\pdfoutput=1
\documentclass[a4paper,11pt]{amsart}
\usepackage[lmargin=60pt,tmargin=56pt,hmarginratio=1:1,vmarginratio=1:1]{geometry}
\usepackage{mathtools,amssymb,amsthm,graphicx,array,booktabs,microtype,xcolor,tikz,needspace,enumitem,etoolbox}
\usepackage[numbers]{natbib}
\IfFormatAtLeastTF{2026-06-01}{}{\usepackage{aliascnt}}

\let\emph\relax
\DeclareTextFontCommand{\emph}{\bfseries\em}
\setlist[enumerate]{itemsep=0.15cm,label=\emph{\upshape(\alph*)}}
\setlist[enumerate,2]{itemsep=0.15cm,label=\emph{\upshape(\roman*)}}
\usetikzlibrary{arrows.meta,decorations.markings,calc,positioning,decorations.pathreplacing,matrix,fit,patterns}
\definecolor{orchid}{RGB}{143,40,194}
\definecolor{mydarkblue}{RGB}{10,10,170}
\definecolor{webink}{gray}{0}
\usepackage[hypertexnames=false,colorlinks=true,linkcolor=mydarkblue,citecolor=teal,urlcolor=orchid,linktocpage=true]{hyperref}
\usepackage{bookmark}
\hypersetup{pdftitle={Positivity in symmetric webs and block permanents},pdfauthor={Daniel Tubbenhauer},pdfstartview=FitH}
\expandafter\def\csname subjclassname@2020\endcsname{2020 Mathematics Subject Classification}

\allowdisplaybreaks
\hypersetup{bookmarksdepth=2}
\numberwithin{equation}{subsection}

\newcommand{\NewTheorem}[2]{%
\IfFormatAtLeastTF{2026-06-01}{%
\newtheorem{#1}[equation]{#2}%
}{%
\newaliascnt{#1}{equation}\newtheorem{#1}[#1]{#2}\aliascntresetthe{#1}%
}%
\expandafter\def\csname #1autorefname\endcsname{#2}}
\theoremstyle{plain}
\NewTheorem{theorem}{Theorem}
\NewTheorem{lemma}{Lemma}
\NewTheorem{proposition}{Proposition}
\NewTheorem{corollary}{Corollary}
\theoremstyle{definition}
\NewTheorem{definition}{Definition}
\NewTheorem{example}{Example}
\AtEndEnvironment{example}{\null\hfill$\Diamond$}
\NewTheorem{remark}{Remark}
\AtEndEnvironment{remark}{\null\hfill$\Diamond$}
\newcommand{\C}{\mathbb C}
\newcommand{\Sym}{\operatorname{Sym}}
\newcommand{\per}{\operatorname{per}}
\newcommand{\tr}{\operatorname{tr}}
\newcommand{\Tr}{\tr}
\newcommand{\id}{\operatorname{id}}
\newcommand{\End}{\operatorname{End}}
\newcommand{\rk}{\operatorname{rank}}
\def\equationautorefname~#1\null{(#1)\null}

\definecolor{wiStrand}{RGB}{32,49,67}
\definecolor{wiActive}{RGB}{0,111,117}
\definecolor{wiFocus}{RGB}{158,82,36}
\definecolor{wiClasp}{RGB}{189,218,239}
\definecolor{wiClaspInk}{RGB}{46,86,119}
\definecolor{wiMapFill}{RGB}{255,229,180}
\definecolor{wiMapInk}{RGB}{134,83,29}
\definecolor{wiHarmonicFill}{RGB}{214,193,230}
\definecolor{wiHarmonicInk}{RGB}{91,61,121}
\definecolor{wiAccentFill}{RGB}{222,240,238}
\definecolor{wiFocusFill}{RGB}{250,235,214}
\colorlet{ink}{wiStrand}
\colorlet{webaccent}{wiActive}
\colorlet{webink}{wiStrand}
\colorlet{webclaspfill}{wiClasp}
\tikzset{
prbase/.style={baseline={([yshift=-.5ex]current bounding box.center)}},
prline/.style={draw=wiStrand,line width=.82pt,line cap=round,line join=round},
prtrack/.style={prline,draw=wiActive,line width=1.13pt},
prcut/.style={draw=wiFocus,line width=.52pt,dash pattern=on 2.25pt off 1.8pt},
prup/.style={prline},
prdown/.style={prline},
prcableup/.style={prup,line width=1.08pt},
prcabledown/.style={prdown,line width=1.08pt},
prvertex/.style={rectangle,draw=black,fill=black,line width=.3pt,inner sep=.93pt},
prsource/.style={circle,draw=black,fill=white,line width=.62pt,inner sep=1.55pt},
prsink/.style={circle,draw=black,fill=black,line width=.62pt,inner sep=1.55pt},
prbox/.style={line width=.58pt,rounded corners=.7pt,minimum height=3.7mm,inner sep=1.6pt,font=\scriptsize,text=black},
prclasp/.style={prbox,draw=wiClaspInk,fill=webclaspfill},
prmorphism/.style={prbox,draw=wiMapInk,fill=wiMapFill},
prharmonic/.style={prbox,draw=wiHarmonicInk,fill=wiHarmonicFill},
web/.style={prline},wa/.style={prup},wr/.style={prdown},clasp/.style={prclasp}
}

\newcommand{\WebBoxSample}[2]{%
\tikz[baseline=(sample.base)]\node[#1,minimum width=7mm] (sample) {$#2$};}
\newcommand{\WebBoxKey}{%
\par\noindent\makebox[\linewidth][c]{\small
\WebBoxSample{prclasp}{P_k}\enspace strand clasp\qquad
\WebBoxSample{prmorphism}{f}\enspace other map\qquad
\WebBoxSample{prharmonic}{J}\enspace trace-free projector}\par}

\definecolor{SWink}{RGB}{38,54,70}

\definecolor{SWteal}{RGB}{6,108,110}

\definecolor{SWgold}{RGB}{167,91,37}

\tikzset{SWwire/.style={draw=SWink,line width=.88pt,line cap=round,line join=round},
SWactive/.style={draw=SWteal,line width=1.02pt,line cap=round,line join=round},
SWmerge/.style={circle,draw=SWink,fill=SWink,inner sep=0pt,minimum size=3.0pt},
SWsplit/.style={circle,draw=SWteal,fill=white,inner sep=0pt,minimum size=3.8pt,line width=.8pt},
SWunit/.style={rectangle,draw=SWink,fill=SWink,inner sep=.82pt},
SWtag/.style={font=\scriptsize,inner sep=1.2pt},
SWseam/.style={draw=SWgold,line width=.60pt,dash pattern=on 1.9pt off 1.7pt}}

\tikzset{trWire/.style={prline},
trCap/.style={prtrack},
trBox/.style={prmorphism,inner sep=1.3pt,minimum height=3.6mm},
trLab/.style={font=\scriptsize,inner sep=.9pt},
trSeam/.style={prcut}}

\tikzset{ftw/.style={draw=black!80,line width=.9pt,line cap=round},
ftc/.style={draw=teal!75!black,line width=1.05pt,line cap=round},
ftbox/.style={prmorphism,inner sep=1.4pt},
ftl/.style={font=\scriptsize,inner sep=.75pt}}

\tikzset{
pdpbase/.style={baseline=(current bounding box.center),
every node/.style={font=\scriptsize,inner sep=1.1pt}},
pdpcable/.style={draw=black!84,line width=1.0pt,line cap=round,line join=round},
pdpcap/.style={pdpcable,draw=teal!78!black,line width=1.15pt},
pdpthin/.style={pdpcable,line width=.65pt},
pdpcut/.style={densely dashed,draw=orange!65!black,line width=.45pt},
pdpbox/.style={prmorphism,inner sep=2pt},
pdpclasp/.style={prharmonic},
pdpstate/.style={prmorphism,inner sep=2pt}
}

\tikzset{lw/.style={draw=black!80,line width=.92pt,line join=round,line cap=round},
lg/.style={draw=teal!75!black,line width=1.08pt,line cap=round},
lt/.style={font=\scriptsize,inner sep=1.2pt},
lp/.style={draw=black!70,fill=blue!7,rounded corners=.7pt,minimum width=9mm,minimum height=3.4mm,inner sep=.5pt,font=\scriptsize}}

\newcommand{\WImatching}[1]{\begin{tikzpicture}[x=7mm,y=6.4mm,prbase]
\node[font=\scriptsize,inner sep=1pt,above] at (0,2.75) {row};
\node[font=\scriptsize,inner sep=1pt,above] at (3,2.75) {column};
\foreach \i/\y in {1/2.55,2/1.92,3/.63,4/0}{
\coordinate (r\i) at (0,\y);\coordinate (c\i) at (3.0,\y);
\node[anchor=east,font=\scriptsize,inner sep=1.2pt] at(r\i){$\i$};
\node[anchor=west,font=\scriptsize,inner sep=1.2pt] at(c\i){$\i$};}
\foreach \pos/\sign in {-.66/1,3.66/-1}{
\foreach \low/\high in {-.13/.76,1.79/2.68}{
\draw[draw=wiStrand!64,line width=.55pt,rounded corners=.65pt]
({\pos+.13*\sign},\low)--(\pos,\low)--(\pos,\high)--({\pos+.13*\sign},\high);}}
\ifnum#1=0
\draw[prtrack] (r1)..controls(.99,2.55)and(2.01,1.92)..(c2);
\draw[prtrack] (r2)..controls(.99,1.92)and(2.01,2.55)..(c1);
\draw[prline] (r3)--(c3);
\else
\draw[prtrack] (r1)--(c1);
\draw[prtrack] (r2)..controls(1.00,1.92)and(2.00,.63)..(c3);
\draw[prline] (r3)..controls(1.00,.63)and(2.00,1.92)..(c2);
\fi
\draw[prline] (r4)--(c4);
\foreach \i in {1,2,3,4}{\fill[wiStrand] (r\i) circle[radius=1pt];
\fill[wiStrand] (c\i) circle[radius=1pt];}
\ifnum#1=0 \node[anchor=north,font=\footnotesize] at(1.50,-.27){whole pairs: $8$ matchings};
\else \node[anchor=north,font=\footnotesize] at(1.50,-.27){split pairs: $16$ matchings};\fi
\end{tikzpicture}}

\newcommand{\ReadCap}[1]{\begin{tikzpicture}[prbase,x=9.5mm,y=7mm,
every node/.style={font=\small,inner sep=1.2pt}]
\draw[prtrack](0,0)--(0,.43)to[out=90,in=90](1.8,.43)--(1.8,0);
\ifnum#1=0 \node[below]at(0,0){$V$};\node[below]at(1.8,0){$V^*$};
\else \node[below]at(0,0){$y$};\node[below]at(1.8,0){$\bar x$};\fi
\end{tikzpicture}}

\newcommand{\ReadTwo}[1]{\begin{tikzpicture}[prbase,x=7.4mm,y=7mm,
every node/.style={font=\small,inner sep=1.2pt}]
\foreach \x/\a in {0/1,2.2/2}{\node[below]at(\x,0){$y_{\a}$};
\node[above]at(\x,1.94){$\bar x_{\a}$};}
\ifnum#1=1
\draw[prline](2.2,0)to[out=90,in=-90](0,1.94);
\draw[prtrack,preaction={draw=white,line width=3pt}](0,0)to[out=90,in=-90](2.2,1.94);
\else
\draw[prtrack](0,0)--(0,1.94);\draw[prline](2.2,0)--(2.2,1.94);
\ifnum#1=2 \node[prclasp,minimum width=22mm]at(1.1,.97){$P_2$};\fi
\fi
\end{tikzpicture}}

\newcommand{\ReadFork}[1]{\begin{tikzpicture}[prbase,x=7mm,y=7mm,
every node/.style={font=\small,inner sep=1.2pt}]
\draw[prline](-1,0)--(-1,.28)to[out=90,in=200](0,1.03);
\draw[prtrack](1,0)--(1,.28)to[out=90,in=-20](0,1.03);
\draw[prline,line width=1.55pt](0,1.03)--(0,2);
\node[rectangle,fill=wiStrand,inner sep=1.25pt]at(0,1.03){};
\ifnum#1=0 \node[below]at(-1,0){$V$};\node[below]at(1,0){$V$};
\node[above]at(0,2){$\Sym^2 V$};
\else \node[below]at(-1,0){$x_1$};\node[below]at(1,0){$x_2$};\node[above]at(0,2){$x_1\odot x_2$};\fi
\end{tikzpicture}}

\newcommand{\ReadCut}[1]{\begin{tikzpicture}[prbase,x=5.5mm,y=5.5mm,
every node/.style={font=\small,inner sep=1.1pt}]
\foreach \x/\lab in {0/F,3.3/{\bar H}}{\draw[prline,line width=1.45pt](\x,0)--(\x,.60);
\node[below]at(\x,0){$\lab$};}
\ifcase#1
\draw[prline,line width=1.45pt](0,.60)--(0,2.5) (3.3,.60)--(3.3,2.5);
\node[above]at(0,2.5){$2$};\node[above]at(3.3,2.5){$2^*$};
\or
\draw[prline](0,.60)to[out=135,in=-90](-.4,1.28)--(-.4,2.5)
(3.3,.60)to[out=45,in=-90](3.7,1.28)--(3.7,2.5);
\draw[prtrack](0,.60)..controls(.58,2.0)and(2.72,2.0)..(3.3,.60);
\foreach \x in {0,3.3}{\node[rectangle,fill=wiStrand,inner sep=1.1pt]at(\x,.60){};}
\node[above]at(-.4,2.5){$1$};\node[above]at(3.7,2.5){$1^*$};
\or
\draw[prtrack,line width=1.45pt](0,.60)to[out=90,in=90](3.3,.60);
\fi
\end{tikzpicture}}

\newcommand{\IntroHalfCut}[1]{%
\begin{tikzpicture}[prbase,x=5.5mm,y=4.6mm,every node/.style={font=\scriptsize,inner sep=1.3pt}]
\draw[prline,line width=1.1pt](0,0)--(0,.72) (3.3,0)--(3.3,.72);
\draw[prline](0,.72)to[out=125,in=-90](-.5,1.37)--(-.5,2.38)
(3.3,.72)to[out=55,in=-90](3.8,1.37)--(3.8,2.38);
\draw[prtrack](0,.72)..controls(.55,1.95)and(2.75,1.95)..(3.3,.72);
\foreach \x in {0,3.3}{\node[rectangle,fill=wiStrand,inner sep=1.1pt]at(\x,.72){};}
\node[below]at(0,0){$F$};\node[below]at(3.3,0){$\bar H$};
\draw[prcut](-1,2.38)--(4.3,2.38);
\ifnum#1=1
\begin{scope}[shift={(0,4.76)},yscale=-1]
\draw[prline,line width=1.1pt](0,0)--(0,.72) (3.3,0)--(3.3,.72);
\draw[prline](0,.72)to[out=125,in=-90](-.5,1.37)--(-.5,2.38)
(3.3,.72)to[out=55,in=-90](3.8,1.37)--(3.8,2.38);
\draw[prtrack](0,.72)..controls(.55,1.95)and(2.75,1.95)..(3.3,.72);
\foreach \x in {0,3.3}{\node[rectangle,fill=wiStrand,inner sep=1.1pt]at(\x,.72){};}
\end{scope}
\node[above]at(0,4.76){$F^\dagger$};\node[above]at(3.3,4.76){$(\bar H)^\dagger$};
\else
\node[above]at(-.50,2.43){$V$};\node[above]at(3.8,2.43){$V^*$};
\fi
\end{tikzpicture}}

\newcommand{\IntroGramStrand}{%
\begin{tikzpicture}[x=7mm,y=7mm,prbase]
\draw[prline] (0,0) -- (0,1.05);
\node[below,inner sep=1pt] at (0,0) {$y$};
\node[above,inner sep=1pt] at (0,1.05) {$\bar x$};
\end{tikzpicture}}

\newcommand{\FigOrganization}{%
\begin{tikzpicture}[x=1cm,y=1cm,prbase,
cell/.style={draw=wiStrand!65,rounded corners=1pt,align=center,text width=4.72cm,minimum height=8mm,inner sep=2.3pt,font=\small},
road arrow/.style={draw=wiStrand!85,-{Stealth[length=1.6mm]},line width=.55pt}]
\node[cell,fill=wiClasp!45,text width=7cm] (start) at (0,0)
{\textbf{\autoref{sec:local}.} The Schur relation\\and the two-block gap};
\node[cell] (a) at (-5.35,-1.30) {\textbf{\autoref{sec:transfer}.} Three blocks\\Transfer and last cap};
\node[cell] (b) at (0,-1.30) {\textbf{\autoref{sec:cubic}.} $m=4$, $k\ge3$\\Two packets on each side};
\node[cell] (c) at (5.35,-1.30) {\textbf{\autoref{sec:threshold}.} Many blocks\\Disjoint two-packet cuts};
\node[cell] (d) at (-5.35,-2.64) {\textbf{\autoref{sec:intersection}--\ref{sec:four}.} $m=4$, $k=2$\\Two caps, one output};
\node[cell] (e) at (5.35,-2.64) {\textbf{\autoref{sec:five}.} $m=5$, $k\ge3$\\A cap with a free cable};
\draw[road arrow](start.south west)--(a.north);\draw[road arrow](start.south)--(b.north);\draw[road arrow](start.south east)--(c.north);
\draw[road arrow](a.south)--(d.north);\draw[road arrow](c.south)--(e.north);
\end{tikzpicture}%
}

\newcommand{\SWforkpair}{\begin{tikzpicture}[x=8.2mm,y=6.5mm,
baseline={([yshift=-.5ex]current bounding box.center)}]
\foreach \SWx/\SWa in {-1.10/p,1.10/q}{
\draw[SWwire] (\SWx,0) to[out=90,in={-90+34*\SWx}] (0,.9);
\draw[SWwire] (0,1.73) to[out={90-34*\SWx},in=-90] (\SWx,2.63);
\node[SWtag,below] at (\SWx,0) {$\SWa$};
\node[SWtag,above] at (\SWx,2.63) {$\SWa$};}
\draw[SWwire,line width=1.20pt] (0,.9)--(0,1.73);
\node[SWtag,right] at (0,1.315) {$p+q$};
\node[SWmerge] at (0,.9) {}; \node[SWsplit] at (0,1.73) {};
\end{tikzpicture}}

\newcommand{\SWexchange}[1]{\begin{tikzpicture}[x=8.2mm,y=6.5mm,
baseline={([yshift=-.5ex]current bounding box.center)}]
\foreach \SWx/\SWa in {0/p,2.6/q}{
\draw[SWwire] (\SWx,0)--(\SWx,2.63);
\node[SWtag,below] at (\SWx,0) {$\SWa$};
\node[SWtag,above] at (\SWx,2.63) {$\SWa$};}
\draw[SWactive] (0,.37) ..controls(.6,.58) and (2,2.05)..(2.6,2.26);
\draw[SWactive,preaction={draw=white,line width=3.0pt}]
(2.6,.37)..controls(2,.58) and (.6,2.05)..(0,2.26);
\foreach \SWx in {0,2.6}{
\ifnum#1=0 \node[SWsplit] at (\SWx,.37) {};\node[SWmerge] at (\SWx,2.26) {};
\else\node[SWunit] at (\SWx,.37) {};\node[SWunit] at (\SWx,2.26) {};\fi}
\node[SWtag,left] at (0,1.315) {$p-r$};
\node[SWtag,right] at (2.6,1.315) {$q-r$};
\node[SWtag,text=SWteal,below] at (.65,.66) {$r$};
\node[SWtag,text=SWteal,below] at (1.95,.66) {$r$};
\end{tikzpicture}}

\newcommand{\SWbent}[1]{\begin{tikzpicture}[x=8.2mm,y=6.5mm,
baseline={([yshift=-.5ex]current bounding box.center)}]
\foreach \SWx/\SWa in {0/p,2.6/{q^*}}{
\draw[SWwire] (\SWx,0)--(\SWx,.37);
\node[SWtag,below] at (\SWx,0) {$\SWa$};}
\draw[SWactive] (0,.37) ..controls(.67,1.17) and (1.93,1.17)..(2.6,.37);
\node[SWtag,text=SWteal] at (1.30,.70) {$r$};
\ifnum#1=0
\foreach \SWx/\SWa in {0/p,2.6/{q^*}}{
\draw[SWwire] (\SWx,.37)--(\SWx,2.63);\node[SWtag,above] at (\SWx,2.63) {$\SWa$};
\node[SWunit] at (\SWx,2.26) {};}
\draw[SWactive] (0,2.26) ..controls(.67,1.46) and (1.93,1.46)..(2.6,2.26);
\draw[SWseam] (-.1,1.315)--(2.7,1.315);
\node[SWtag,left] at (0,1.53) {$p-r$};
\node[SWtag,right] at (2.6,1.53) {$(q-r)^*$};
\else
\foreach \SWx/\SWa in {0/{p-r},2.6/{(q-r)^*}}{
\draw[SWwire] (\SWx,.37)--(\SWx,1.46);\node[SWtag,above] at (\SWx,1.46) {$\SWa$};}
\fi
\foreach \SWx in {0,2.6}{\node[SWunit] at (\SWx,.37) {};}
\end{tikzpicture}}

\newcommand{\SWsmall}[2]{\begin{tikzpicture}[x=9mm,y=6.2mm,
baseline={([yshift=-.5ex]current bounding box.center)}]
\pgfmathsetmacro{\SWleft}{.32*(#1-1)}
\pgfmathsetmacro{\SWright}{\SWleft+.68}
\pgfmathsetmacro{\SWlast}{\SWright+\SWleft}
\foreach \SWj in {1,...,#1}{
\pgfmathsetmacro{\SWx}{\SWleft-.32*(\SWj-1)}
\pgfmathsetmacro{\SWy}{\SWright+.32*(\SWj-1)}
\ifnum\SWj>#2
\draw[SWwire] (\SWx,0)--(\SWx,2.50) (\SWy,0)--(\SWy,2.50);
\else
\pgfmathsetmacro{\SWmid}{(\SWx+\SWy)/2}
\pgfmathsetmacro{\SWrad}{(\SWy-\SWx)/2}
\pgfmathsetmacro{\SWht}{.30+.18*(\SWj-1)}
\draw[SWactive] (\SWx,0)--(\SWx,.46)
arc[start angle=180,end angle=0,x radius=\SWrad,y radius=\SWht]--(\SWy,0);
\draw[SWactive] (\SWx,2.50)--(\SWx,2.04)
arc[start angle=180,end angle=360,x radius=\SWrad,y radius=\SWht]--(\SWy,2.50);
\fi}
\foreach \SWlevel in {.24,2.26}{
\draw[SWwire,fill=blue!14,rounded corners=.45pt]
(-.12,{\SWlevel-.115}) rectangle ({\SWleft+.12},{\SWlevel+.115});
\draw[SWwire,fill=blue!14,rounded corners=.45pt]
({\SWright-.12},{\SWlevel-.115}) rectangle ({\SWlast+.12},{\SWlevel+.115});}
\draw[SWseam] (-.16,1.25)--({\SWlast+.16},1.25);
\end{tikzpicture}}

\newcommand{\LocThreeWords}{\begin{tikzpicture}[x=6mm,y=4.5mm,prbase,every node/.style={SWtag}]
\draw[SWwire,line width=1.10pt](-1.55,0)to[out=90,in=210](0,1.15)
(0,1.15)--(0,2.70) (0,2.70)to[out=145,in=-90](-1.55,3.85);
\draw[SWactive](1.55,0)to[out=90,in=-30](0,1.15)
(0,2.70)to[out=35,in=-90](1.55,3.85);
\node[SWunit]at(0,1.15){};\node[SWunit]at(0,2.70){};
\node[below]at(-1.55,0){$x^2$};\node[below]at(1.55,0){$y$};
\node[right]at(.05,1.94){$3$};
\node[above]at(-1.55,3.85){$2$};\node[above]at(1.55,3.85){$1$};
\end{tikzpicture}}

\newcommand{\LocRecur}[1]{\begin{tikzpicture}[x=5.8mm,y=4.9mm,
baseline={([yshift=-.5ex]current bounding box.center)},every node/.style={SWtag}]
\draw[SWwire,line width=1.15pt] (0,0)--(0,3.5);
\ifnum#1=2
\draw[SWactive](1.6,0)--(1.6,.91)..controls(1.6,1.64)and(3,1.85)..(3,2.59)--(3,3.5);
\draw[SWactive,preaction={draw=white,line width=3.6pt}]
(3,0)--(3,.91)..controls(3,1.64)and(1.6,1.85)..(1.6,2.59)--(1.6,3.5);
\else
\draw[SWwire](1.6,0)--(1.6,3.5) (3,0)--(3,3.5);
\fi
\ifnum#1=0
\draw[SWwire,fill=blue!12,rounded corners=.5pt](-.21,1.40)rectangle(3.21,2.1);
\node at(1.5,1.75){$P_{p+1}$};
\else\ifnum#1=1
\draw[SWwire,fill=blue!12,rounded corners=.5pt](-.21,1.40)rectangle(1.81,2.1);
\node at(.8,1.75){$P_p$};
\else
\foreach \yy in {.55,2.95}{
\draw[SWwire,fill=blue!12,rounded corners=.5pt](-.21,{\yy-.34})rectangle(1.81,{\yy+.34});
\node at(.8,\yy){$P_p$};}
\fi\fi
\foreach \xx/\lab in {0/{p-1},1.6/1,3/1}{\node[below]at(\xx,0){$\lab$};}
\end{tikzpicture}}

\newcommand{\LocTracePair}[1]{\begin{tikzpicture}[x=6mm,y=4.2mm,prbase,every node/.style={SWtag}]
\ifnum#1=0
\draw[SWactive](0,1.45)ellipse[x radius=.85,y radius=1.03];
\node[left]at(-.85,1.45){$W$};\node[right]at(.85,1.45){$W^*$};
\draw[SWseam](-1.26,1.45)--(1.26,1.45);
\else
\draw[SWactive](-.8,0)--(-.8,.42)to[out=90,in=90](.8,.42)--(.8,0)
(-.8,2.90)--(-.8,2.48)to[out=-90,in=-90](.8,2.48)--(.8,2.90);
\node[below]at(-.8,0){$w$};\node[below]at(.8,0){$\phi$};
\node[above]at(-.8,2.90){$W$};\node[above]at(.8,2.90){$W^*$};
\fi\end{tikzpicture}}

\newcommand{\LocClosedBend}[1]{\begin{tikzpicture}[x=6.7mm,y=5.4mm,
baseline=(current bounding box.center),every node/.style={SWtag}]
\draw[SWwire](0,0)--(0,1.17) (0,2.03)--(0,3.26);
\node[below]at(0,0){$F$};\node[above]at(0,3.26){$F^\dagger$};
\ifnum#1=0
\draw[SWactive](1.45,0)--(1.45,1.17) (1.45,2.03)--(1.45,3.26);
\node[below]at(1.45,0){$H$};\node[above]at(1.45,3.26){$H^\dagger$};
\else
\draw[SWactive](1.45,1.17)--(1.45,.81)
..controls(1.45,.18)and(2.37,.18)..(2.37,.81)--(2.37,3.26);
\draw[SWactive,preaction={draw=white,line width=3.2pt}](1.45,2.03)--(1.45,2.37)
..controls(1.45,3.18)and(3.14,3.18)..(3.14,2.37)--(3.14,0);
\node[below]at(3.14,0){$\bar H$};\node[above]at(2.37,3.26){$(\bar H)^\dagger$};
\fi
\draw[draw=wiMapInk,fill=wiMapFill,rounded corners=.6pt](-.22,1.17)rectangle(1.67,2.03);
\node at(.725,1.6){$B$};
\end{tikzpicture}}

\newcommand{\LocTurnExchange}{\begin{tikzpicture}[x=5.45mm,y=5.9mm,prbase,every node/.style={SWtag}]
\draw[SWwire](0,0)--(0,3.13) (1.75,.83)--(1.75,2.29);
\draw[SWactive](0,.83)--(1.75,2.29);
\draw[SWactive,preaction={draw=white,line width=3pt}](1.75,.83)--(0,2.29);
\draw[SWwire](1.75,.83)--(1.75,.66)to[out=-90,in=-90](2.76,.66)--(2.76,3.13);
\draw[SWwire,preaction={draw=white,line width=3pt}](1.75,2.29)--(1.75,2.46)
to[out=90,in=90](3.63,2.46)--(3.63,0);
\foreach \xx in {0,1.75}{\foreach \yy in {.83,2.29}{\node[SWunit]at(\xx,\yy){};}}
\node[below]at(0,0){$p$};\node[below]at(3.63,0){$q^*$};
\node[above]at(0,3.13){$p$};\node[above]at(2.76,3.13){$q^*$};
\end{tikzpicture}}

\newcommand{\LocSeamInsertion}[1]{\begin{tikzpicture}[x=7mm,y=4.1mm,prbase,every node/.style={SWtag}]
\draw[SWwire](0,0)--(0,4.75);
\node[prmorphism]at(0,.66){$h$};\node[prmorphism]at(0,4.09){$h^\dagger$};
\ifnum#1=1\node[prmorphism]at(0,1.65){$J$};\node[prmorphism]at(0,3.10){$J$};\fi
\ifnum#1=2\node[prmorphism]at(0,1.65){$I-J$};\node[prmorphism]at(0,3.10){$I-J$};\fi
\draw[SWseam](-.92,2.375)--(.92,2.375);
\node[right]at(.96,2.375){$Z$};\node[below]at(0,0){$X$};
\end{tikzpicture}}

\newcommand{\LocClos}[2]{\begin{tikzpicture}[x=6.8mm,y=5.5mm,
baseline={([yshift=-.5ex]current bounding box.center)},every node/.style={SWtag}]
\foreach \x/\i in {-2/1,-.8/2,2/m}{
\node[below] at (\x,0) {$F_{\i}$};
\node[above] at (\x,4.35) {$F_{\i}^\dagger$};
\draw[SWwire] (\x,0)--(\x,.50) to[out=90,in={-90+19*\x}](0,1.41);
\draw[SWwire] (0,2.94) to[out={90-19*\x},in=-90](\x,3.85)--(\x,4.35);}
\node at (.62,.32){$\cdots$}; \node at (.62,4.02){$\cdots$};
\draw[SWactive] (0,1.41)--(0,2.94);
\node[SWunit] at (0,1.41){}; \node[SWunit] at (0,2.94){};
\draw[SWseam] (-.87,2.18)--(.85,2.18);
\node[anchor=west,fill=white] at (.15,2.55){$#1$};
\node[left] at (-1.97,.53){$#2$};
\end{tikzpicture}}

\newcommand{\LocPacketPair}{\begin{tikzpicture}[x=6.3mm,y=4.8mm,prbase,every node/.style={SWtag}]
\foreach \x in {0,1.15,3.25}{
\draw[SWwire] (\x,0)--(\x,3.8);}
\foreach \yy in {.65,3.16}{
\node[prclasp,minimum width=25mm] at(1.625,\yy){$P_k$};}
\node[below]at(0,0){$x_{j,1}$};\node[below]at(1.15,0){$x_{j,2}$};
\node[below]at(2.18,0){$\cdots$};\node[below]at(3.25,0){$x_{j,k}$};
\node[above]at(0,3.8){$\bar x_{i,1}$};\node[above]at(1.15,3.8){$\bar x_{i,2}$};
\node[above]at(2.18,3.8){$\cdots$};\node[above]at(3.25,3.8){$\bar x_{i,k}$};
\node[right]at(3.78,.65){$\sqrt{k!}$};\node[right]at(3.78,3.16){$\sqrt{k!}$};
\node at(2.18,2.32){$\cdots$};\draw[SWseam](-.54,1.90)--(3.78,1.9);
\end{tikzpicture}}

\newcommand{\LocBoundary}[1]{\begin{tikzpicture}[x=6.2mm,y=4.4mm,prbase,every node/.style={SWtag}]
\ifnum#1=0
\foreach \x in {0,.85,2.15}{\draw[SWwire](\x,0)--(\x,4);}
\foreach \yy/\ll in {.62/{\sigma},2/{B},3.38/{\sigma^{-1}}}{\node[prmorphism,minimum width=17mm] at(1.075,\yy){$\ll$};}
\else
\draw[SWwire](0,0)--(0,4);\draw[SWactive](2.15,0)--(2.15,4);
\foreach \yy/\ll in {.62/{j},3.38/{j^\dagger}}{\node[prmorphism] at(2.15,\yy){$\ll$};}
\node[prmorphism,minimum width=17mm] at(1.075,2){$B$};
\node[below]at(2.15,0){$K'$};\node[above]at(2.15,4){$K'$};
\fi
\end{tikzpicture}}

\newcommand{\FigDigonRemoval}{%
\begin{tikzpicture}[x=6.8mm,y=5.5mm,baseline=(current bounding box.center),every node/.style={SWtag}]
\draw[SWwire](0,0)--(0,.62) (0,2.17)--(0,2.79);
\draw[SWwire](0,.62)..controls(-1.08,1.02)and(-1.08,1.77)..(0,2.17)
(0,.62)..controls(1.08,1.02)and(1.08,1.77)..(0,2.17);
\node[SWunit]at(0,.62){};\node[SWunit]at(0,2.17){};
\node[left]at(-.81,1.39){$p$};\node[right]at(.81,1.39){$q$};
\end{tikzpicture}%
}

\newcommand{\FigDigonIdentity}{%
\begin{tikzpicture}[x=6.8mm,y=5.5mm,baseline=(current bounding box.center),every node/.style={SWtag}]
\draw[SWwire](0,0)--(0,2.79);\node[right]at(0,1.4){$p+q$};
\end{tikzpicture}%
}

\newcommand{\FigMergeSplit}{%
\begin{tikzpicture}[x=6.8mm,y=5.5mm,baseline=(current bounding box.center),every node/.style={SWtag}]
\foreach \x/\a in {-1/p,1/q}{
\draw[SWwire](\x,0) to[out=90,in={-90+38*\x}](0,.92)
(0,1.87) to[out={90-38*\x},in=-90](\x,2.79);
\node[below]at(\x,0){$\a$};}
\draw[SWwire](0,.92)--(0,1.87);\node[SWunit]at(0,.92){};\node[SWunit]at(0,1.87){};
\end{tikzpicture}%
}

\newcommand{\FigTotalClasp}{%
\begin{tikzpicture}[x=6.8mm,y=5.5mm,baseline=(current bounding box.center),every node/.style={SWtag}]
\foreach \x/\a in {-1/p,1/q}{\draw[SWwire](\x,0)--(\x,2.79);\node[below]at(\x,0){$\a$};}
\draw[draw=SWink,fill=blue!12,rounded corners=.5pt](-1.18,1.02) rectangle(1.18,1.80);
\node at(0,1.40){$P_{p+q}$};
\end{tikzpicture}%
}

\newcommand{\FigSubsetCut}{%
\begin{tikzpicture}[x=6.2mm,y=4.6mm,baseline=(current bounding box.center),
every node/.style={SWtag}]
\foreach \x/\SWinput/\SWoutput in {0/{F}/{2},3.2/{\bar F}/{2^*}}{
\draw[SWwire](\x,0)--(\x,2.5);
\node[below]at(\x,0){$\SWinput$};\node[above]at(\x,2.5){$\SWoutput$};
\node[SWunit]at(\x,.61){};}
\draw[SWactive](0,.61) ..controls(.7,1.65)and(2.5,1.65)..(3.2,.61);
\node[text=SWteal]at(1.6,.97){$2$};
\draw[SWseam](-.28,2.23)--(3.48,2.23);
\end{tikzpicture}%
}

\newcommand{\TrOpen}[1]{\begin{tikzpicture}[x=5.8mm,y=5.3mm,baseline=(current bounding box.center),every node/.style={trLab}]
\foreach \x in {0,1.3,2.9}{\draw[trWire](\x,0)--(\x,.7);}
\draw[trWire](4.95,0)--(4.95,1.45);
\node at(2.08,.26){$\cdots$};
\draw[decorate,decoration={brace,mirror,amplitude=2pt}](0,-.06)--(2.9,-.06);
\node[below]at(1.45,-.23){$n$ packets};\node[below]at(4.95,-.05){$k^*$};
\ifnum#1=0
\node[trBox,minimum width=19mm]at(1.45,.7){$M_n$};
\draw[trWire](1.45,1.06)--(1.45,1.45);
\node[left]at(1.45,1.28){$nk$};
\draw[trWire](1.45,1.45)to[out=130,in=-90](.9,2.12)--(.9,3.05)
(4.95,1.45)to[out=50,in=-90](5.5,2.12)--(5.5,3.05);
\draw[trCap](1.45,1.45)..controls(2.1,2.08)and(2.12,2.48)..(3.20,2.48)
..controls(4.28,2.48)and(4.3,2.08)..(4.95,1.45);
\foreach \x in {1.45,4.95}{\fill(\x,1.45)circle[radius=1.2pt];}
\node[below]at(3.20,2.42){$r$};\node[above]at(.9,3.09){$nk-r$};
\node[above]at(5.5,3.09){$(k-r)^*$};
\draw[trSeam](.53,3.05)--(5.88,3.05);
\else
\node[trBox,minimum width=19mm]at(1.45,.7){$f$};
\draw[trWire](0,1.06)--(0,3.05) (1.3,1.06)--(1.3,3.05) (2.9,1.06)--(2.9,1.45);
\node at(.65,1.64){$\cdots$};
\draw[trCap](2.9,1.45)..controls(2.9,2.83)and(4.95,2.83)..(4.95,1.45);
\node[below]at(3.93,2.39){$k$};
\node[above]at(.65,3.09){$n-1$ packets};
\draw[trSeam](-.4,3.05)--(1.68,3.05);
\fi
\end{tikzpicture}}

\newcommand{\TrFirstMark}[3]{\begin{tikzpicture}[x=5.4mm,y=4.25mm,prbase,every node/.style={font=\scriptsize,inner sep=.8pt}]
\foreach \x/\lab in {0/#1,1.6/#2,4/#3}{
\draw[prline](\x,0)--(\x,2.25);
\node[below]at(\x,0){$\lab$};\node[above]at(\x,2.25){$\lab^\dagger$};}
\draw[prtrack](4,.14)--(4,2.10);
\node[prmorphism,minimum width=13mm]at(.80,1.13){$B$};
\end{tikzpicture}}

\newcommand{\FTexit}[1]{%
\begin{tikzpicture}[x=4.1mm,y=4.1mm,
  baseline=(current bounding box.center)]
\foreach \x/\lab in {0/2,1/3,2/4}{
  \draw[ftw] (\x,0)--(\x,1.10);
  \node[ftl,below] at (\x,0) {$\lab$};
}
\ifnum#1<0
  \foreach \x in {0,1,2}{
    \draw[ftw] (\x,1.10)--(\x,2.7);
  }
  \node[ftbox,minimum width=11.4mm] at (1,1.32) {$f'$};
\else
  \ifnum#1=0
    \foreach \x in {0,1,2}{
      \draw[ftw] (\x,1.10)--(\x,2.7);
    }
  \else
    \def\ftstep##1##2##3{%
      \pgfmathtruncatemacro{\ftleft}{##1}
      \pgfmathtruncatemacro{\ftright}{##1+1}
      \pgfmathtruncatemacro{\ftkeep}{2-2*##1}
      \draw[ftw] (\ftkeep,##2)--(\ftkeep,##3);
      \draw[ftc] (\ftleft,##2)
        to[out=90,in=-90] (\ftright,##3);
      \draw[ftw,preaction={draw=white,line width=2.4pt}]
        (\ftright,##2) to[out=90,in=-90] (\ftleft,##3);
    }
    \ifnum#1=1
      \ftstep{0}{1.10}{1.63}
      \ftstep{1}{1.63}{2.16}
      \ftstep{0}{2.16}{2.70}
    \else
      \ftstep{1}{1.10}{2.70}
    \fi
  \fi
  \node[ftbox,minimum width=7.1mm] at (.5,.53) {$h$};
\fi
\end{tikzpicture}%
}

\newcommand{\FTmixed}[1]{\begin{tikzpicture}[x=6mm,y=4.0mm,baseline=(current bounding box.center)]
\draw[ftw,line width=1.25pt](0,0)--(0,3.85);
\ifnum#1=0
\draw[ftw](1.26,0)--(1.26,.9) to[out=90,in=-90](2.68,2.95)--(2.68,3.85);
\draw[ftc,preaction={draw=white,line width=2.8pt}](2.68,0)--(2.68,.9)
to[out=90,in=-90](1.26,2.95)--(1.26,3.85);
\else
\draw[ftc](1.26,0)--(1.26,.95)..controls(1.26,1.91)and(2.68,1.91)..(2.68,.95)--(2.68,0);
\draw[ftc](1.26,3.85)--(1.26,2.90)..controls(1.26,1.94)and(2.68,1.94)..(2.68,2.90)--(2.68,3.85);
\draw[orange!65!black,densely dashed](-.3,1.925)--(2.94,1.925);
\fi
\node[ftbox,minimum width=11.4mm]at(.63,.45){$f$};
\node[ftbox,minimum width=11.4mm]at(.63,3.40){$B$};
\node[ftl,anchor=north east,xshift=2pt]at(0,0){$1,\dots,n{-}1$};\node[ftl,below]at(1.26,0){$n$};
\ifnum#1=0\node[ftl,below]at(2.68,0){$m$};\else\node[ftl,below]at(2.68,0){$m^*$};\fi
\end{tikzpicture}}

\newcommand{\PDPhalf}[1]{%
\begin{tikzpicture}[x=6.5mm,y=5.8mm,pdpbase]
\foreach \x/\lab in {0/k,1.48/k,3.85/{k^*}}
\node[below] at (\x,0) {$\lab$};
\ifnum#1=2
\draw[pdpcable] (0,0) to[out=90,in=212] (.74,.84);
\draw[pdpcable] (1.48,0) to[out=90,in=-32] (.74,.84)
to[out=66,in=-90] (1.48,1.88);
\fill (.74,.84) circle[radius=1.25pt];
\node[left] at (1.27,1.47) {$2k$};
\else
\ifnum#1=0
\draw[pdpcable](0,0)--(0,1.30) (1.48,0)--(1.48,1.30);
\else
\draw[pdpcable](0,0) to[out=90,in=-90] (1.48,1.30);
\draw[pdpcable,preaction={draw=white,line width=2.6pt}]
(1.48,0) to[out=90,in=-90] (0,1.30);
\fi
\draw[pdpcable](0,1.30)--(0,3.60) (1.48,1.30)--(1.48,1.88);
\node[above left] at (.12,3.60) {$k$};
\fi
\draw[pdpcable](3.85,0)--(3.85,1.88);
\draw[pdpcable](1.48,1.88) to[out=130,in=-90] (.98,2.52)--(.98,3.60);
\draw[pdpcable](3.85,1.88) to[out=50,in=-90] (4.35,2.52)--(4.35,3.60);
\draw[pdpcap](1.48,1.88)..controls(1.96,2.00)and(1.88,2.91)..(2.665,2.91)
..controls(3.45,2.91)and(3.37,2.00)..(3.85,1.88);
\node[below,text=teal!78!black] at (2.665,2.89) {$r$};
\foreach \x in {1.48,3.85}\fill (\x,1.88) circle[radius=1.25pt];
\ifnum#1=2
\node[above] at (.98,3.60) {$2k-r$};
\else
\node[above right] at (.83,3.60) {$k-r$};
\fi
\node[above] at (4.35,3.60) {$(k-r)^*$};
\end{tikzpicture}}

\newcommand{\PDPglue}[1]{%
\begin{tikzpicture}[x=5.8mm,y=5.9mm,pdpbase]
\foreach \x in {0,1.25,2.50}{
\draw[pdpcable](\x,-.32)--(\x,.23) (\x,2.42)--(\x,2.97);}
\draw[pdpcable](0,.23)--(0,2.42);
\ifnum#1=2
\draw[pdpcap](1.25,.23)--(1.25,.60)
..controls(1.25,1.32)and(2.50,1.32)..(2.50,.60)--(2.50,.23);
\draw[pdpcap](1.25,2.42)--(1.25,2.05)
..controls(1.25,1.33)and(2.50,1.33)..(2.50,2.05)--(2.50,2.42);
\else
\foreach \x in {1.25,2.50}\draw[pdpthin](\x,.23)--(\x,2.42);
\ifnum#1=1
\node[pdpclasp,minimum width=10.6mm] at (1.875,.84) {$J$};
\node[pdpclasp,minimum width=10.6mm] at (1.875,1.83) {$J$};
\fi
\fi
\node[pdpbox,minimum width=18mm] at (1.25,.15) {$U_{k-1}$};
\node[pdpbox,minimum width=18mm] at (1.25,2.52) {$U_{k-1}^{\dagger}$};
\draw[pdpcut](-.34,1.33)--(2.88,1.33);
\node[left] at (0,1.60) {$K$};
\ifnum#1=0
\node[left] at (1.25,1.74) {$V$};
\node[right] at (2.50,1.74) {$S^*$};
\fi
\end{tikzpicture}}

\newcommand{\PDPloop}[1]{%
\begin{tikzpicture}[x=5mm,y=4.6mm,pdpbase]
\draw[pdpcap](0,.38)..controls(0,-.22)and(1.55,-.22)..(1.55,.38)
--(1.55,1.65)..controls(1.55,2.25)and(0,2.25)..(0,1.65)--cycle;
\ifnum#1=1
\node[pdpstate] at (0,.56) {$j$};
\node[pdpstate] at (0,1.48) {$j^\dagger$};
\fi
\node at (.87,1.03) {$S$};
\end{tikzpicture}}

\newcommand{\TrFullHalf}[1]{\begin{tikzpicture}[x=5.2mm,y=4.6mm,baseline=(current bounding box.center),every node/.style={trLab}]
\foreach \x/\name in {0/k,1.85/k,4/{k^*}}{\node[below]at(\x,0){$\name$};}
\ifnum#1=0
\foreach \x in {0,1.85}{\draw[trWire](\x,0)--(\x,1.14);}
\node[trBox,minimum width=13.2mm]at(.925,.68){$f$};
\draw[trWire](0,1.14)to[out=90,in=-90](.2,2.65)--(.2,3.22);
\draw[trCap](1.85,1.14)--(1.85,1.55)to[out=90,in=90,looseness=1.25](4,1.55)--(4,0);
\else
\draw[trWire](0,0)to[out=90,in=202](.925,.88) (1.85,0)to[out=90,in=-22](.925,.88)
(.925,.88)--(.925,1.68) (.925,1.68)to[out=140,in=-90](.2,2.48)--(.2,3.22);
\draw[trCap](.925,1.68)to[out=42,in=90](4,1.28)--(4,0);
\foreach \y in {.88,1.68}{\fill(.925,\y)circle(1.3pt);}
\node[left]at(.85,1.3){$2k$};
\fi
\draw[trSeam](-.46,3.15)--(.98,3.15);\node[above]at(.2,3.22){$k$};
\end{tikzpicture}}

\newcommand{\TrCompleteCap}[1]{\begin{tikzpicture}[x=5.7mm,y=4.4mm,prbase,every node/.style={font=\scriptsize,inner sep=.9pt}]
\foreach \x/\lab in {0/k,1.65/k,4/{k^*}}{\node[below]at(\x,0){$\lab$};}
\draw[prline](0,0)--(0,3.65);
\node[above]at(0,3.65){$k$};
\draw[prline](1.65,0)--(1.65,1.11) (4,0)--(4,1.11);
\node[prmorphism,minimum width=13.5mm]at(.825,.55){$f$};
\ifnum#1=0
\draw[prline](1.65,1.11)--(1.65,2.65) (4,1.11)--(4,2.65);
\draw[prtrack](1.65,1.11)to[out=42,in=138](4,1.11);
\draw[prtrack](1.65,2.65)to[out=90,in=90,looseness=.9](4,2.65);
\foreach \x in {1.65,4}{\fill[wiStrand](\x,1.11)circle(1.3pt);}
\draw[prcut](-.4,2.27)--(4.45,2.27);
\node at(2.825,1.15){$k-1$};\node at(2.825,3.13){$1$};
\else
\draw[prtrack](1.65,1.11)to[out=90,in=90,looseness=1.14](4,1.11);
\node at(2.82,1.74){$k$};
\fi
\end{tikzpicture}}

\newcommand{\TrIterExit}[1]{\begin{tikzpicture}[x=6.2mm,y=4.3mm,baseline=(current bounding box.center),
every node/.style={font=\scriptsize,inner sep=1pt}]
\draw[trWire](0,0)--(0,4.15);
\ifnum#1=0
\draw[trWire](1.45,0)--(1.45,4.15);
\draw[trCap](3,0)--(3,4.15);
\node[below]at(1.45,0){$l$};\node[above]at(1.45,4.15){$i$};
\else
\draw[trWire](3,0)--(3,.35)to[out=90,in=-90](1.45,1.38)--(1.45,2.75)to[out=90,in=-90](3,3.82)--(3,4.15);
\draw[trCap,preaction={draw=white,line width=2.8pt}](1.45,0)--(1.45,.35)to[out=90,in=-90](3,1.38)--(3,2.75)to[out=90,in=-90](1.45,3.82)--(1.45,4.15);
\node[below]at(1.45,0){$l$};\node[above]at(1.45,4.15){$i$};
\fi
\node[trBox,minimum width=12.5mm]at(.725,2.05){$D$};
\node[below]at(3,0){$j$};\node[above]at(3,4.15){$j$};
\end{tikzpicture}}

\newcommand{\TrEdgeType}[2]{\begin{tikzpicture}[x=6mm,y=5mm,baseline=(current bounding box.center),every node/.style={trLab}]
\foreach \y in {0,1.62}{
\fill[wiStrand](0,\y)circle(1.55pt);
\draw[draw=wiActive,fill=white,line width=.85pt](1.5,\y)circle(2.2pt);}
\draw[trCap](#1,.11)--(#2,1.51);
\end{tikzpicture}}

\newcommand{\TrCorrectionPath}[1]{\begin{tikzpicture}[x=6mm,y=6mm,prbase]
\foreach\x in{0,1.25,2.5}\draw[prline,line width=1.08pt](\x,0)--(\x,3.5);
\ifcase#1\relax
\node[prmorphism,minimum width=14mm]at(.625,.55){$f$};
\node[prclasp,minimum width=14mm]at(1.875,1.75){$\Pi'$};
\node[prmorphism,minimum width=14mm]at(.625,2.95){$f$};
\or
\node[prmorphism,minimum width=14mm]at(.625,.55){$f$};
\node[prmorphism,minimum width=14mm]at(1.875,1.75){$g$};
\node[prclasp,minimum width=14mm]at(.625,2.95){$\Pi$};
\or
\node[prclasp,minimum width=14mm]at(.625,.55){$\Pi$};
\node[prmorphism,minimum width=14mm]at(1.875,1.75){$g$};
\node[prmorphism,minimum width=14mm]at(.625,2.95){$f$};
\fi
\end{tikzpicture}}

\newcommand{\FigPacketTransfer}{%
\begin{tikzpicture}[x=1cm,y=1cm,baseline=(current bounding box.center),
recut wire/.style={draw=blue!42!black,line width=1.12pt,line cap=round},
recut bracket/.style={draw=black!68,line width=.5pt},
recut map/.style={prmorphism,
line width=.6pt,minimum width=9.5mm,minimum height=4.2mm,
inner xsep=1pt,inner ysep=0pt,font=\small},
recut flag/.style={recut map},
recut label/.style={font=\scriptsize,text=black!70,inner sep=.3pt},
recut sign/.style={font=\normalsize,inner sep=.3pt}]
\def\recutpanel##1##2{%
\begin{scope}[shift={(##1,0)}]
\draw[recut bracket] (-.26,1.98)--(-.35,1.98)--(-.35,-.05)--(-.26,-.05);
\draw[recut bracket] (1.50,1.98)--(1.59,1.98)--(1.59,-.05)--(1.50,-.05);
\foreach \x/\j in {0/1,.62/2,1.24/3}{
\node[recut label,below] at (\x,-.09) {$\j$};
}
\draw[recut wire] (0,0)--(0,1.93);
\ifnum##2<2
\draw[recut wire] (.62,0)--(.62,1.93) (1.24,0)--(1.24,1.93);
\else
\foreach \x in {.62,1.24}
\draw[recut wire] (\x,0)--(\x,.58) (\x,1.35)--(\x,1.93);
\draw[recut wire] (.62,.58)..controls(.62,.93)and(1.24,1.00)..(1.24,1.35);
\draw[recut wire,preaction={draw=white,line width=3.1pt}]
(1.24,.58)..controls(1.24,.93)and(.62,1.00)..(.62,1.35);
\fi
\ifcase##2
\node[recut map] at (.31,.965) {$B$};
\or
\node[recut flag] at (.31,.35) {$f$};
\node[recut map] at (.93,.965) {$B'$};
\node[recut flag] at (.31,1.58) {$f$};
\or
\node[recut flag] at (.31,.35) {$f$};
\node[recut map] at (.31,1.58) {$B$};
\or
\node[recut map] at (.31,.35) {$B$};
\node[recut flag] at (.31,1.58) {$f$};
\fi
\end{scope}}
\recutpanel{0}{0}
\node[recut sign] at (2.02,.965) {$=$};
\node[recut sign] at (2.49,.965) {$2$};
\recutpanel{3.09}{1}
\node[recut sign] at (5.12,.965) {$-$};
\node[recut sign] at (5.62,.965) {$\tfrac12$};
\node[recut sign] at (6.02,.965) {$\Bigg($};
\recutpanel{6.60}{2}
\node[recut sign] at (8.68,.965) {$+$};
\recutpanel{9.48}{3}
\node[recut sign] at (11.36,.965) {$\Bigg)$};
\end{tikzpicture}%
}

\newcommand{\FigThreePacketRouting}{%
\begin{tikzpicture}[x=6mm,y=4.5mm,prbase,every node/.style={font=\scriptsize,inner sep=.8pt}]
\foreach \x/\lab in {0/1,2/2,4.5/3}{
\node[below]at(\x,0){$\lab$};\node[above]at(\x,3.2){$\lab$};
\draw[prtrack,line width=1.3pt](\x,0)--(\x,.55) (\x,2.65)--(\x,3.2);}
\foreach \x in {0,2}{\draw[prline](\x,.55)--(\x,2.65);}
\draw[prline](0,.55)to[out=45,in=225](2,2.65);
\draw[prline,preaction={draw=white,line width=2.6pt}](2,.55)to[out=135,in=-45](0,2.65);
\foreach \x in {0,2}{\foreach \yy in {.55,2.65}{\fill[wiStrand](\x,\yy)circle(1.45pt);}}
\draw[prtrack,line width=1.3pt](4.5,.55)--(4.5,2.65);
\node[right]at(4.53,1.6){$2$};
\end{tikzpicture}%
}

\newcommand{\IntPartner}[1]{\begin{tikzpicture}[prbase,x=5.8mm,y=4.8mm,every node/.style={font=\scriptsize,inner sep=1pt}]
\ifnum#1=0\def\chosen{0}\def\free{2.6}\else\def\chosen{2.6}\def\free{0}\fi
\draw[prline](\free,0)--(\free,2.65);
\draw[prtrack](\chosen,0)--(\chosen,.7)to[out=90,in=90,looseness=1.35](1.3,.7)--(1.3,0);
\node[below]at(0,0){$x$};\node[below]at(1.3,0){$\bar e$};\node[below]at(2.6,0){$y$};
\end{tikzpicture}}

\newcommand{\IntGlue}[1]{\begin{tikzpicture}[prbase,x=4.6mm,y=3.8mm,
every node/.style={font=\scriptsize,inner sep=1pt}]
\ifnum#1=0\def\cupend{0}\def\straight{2.6}\else\def\cupend{2.6}\def\straight{0}\fi
\draw[prline] (\straight,-.35)--(\straight,1.3) (2.6,1.3)--(2.6,3.0);
\draw[prtrack] (\cupend,1.3)--(\cupend,1.02)to[out=-90,in=-90](1.3,1.02)--(1.3,1.30);
\draw[prtrack] (0,1.3)--(0,1.63)to[out=90,in=90](1.3,1.63)--(1.3,1.30);
\node[prmorphism,inner sep=.9pt]at(1.3,1.3){$S$};
\end{tikzpicture}}

\newcommand{\FigMultipleCaps}{%
\begin{tikzpicture}[x=7mm,y=5.5mm,baseline=(current bounding box.center),
every node/.style={font=\scriptsize,inner sep=1.2pt},
v/.style={draw=black!85,line width=.96pt,line cap=round},
d/.style={draw=teal!70!black,line width=1.08pt,line cap=round}]
\draw[v](0,0)--(0,1.0)--(-.5,1.56)--(-.5,2.8)
to[out=90,in=205](2,3.65);
\draw[v](4,0)--(4,1.0)--(4.5,1.56)--(4.5,2.8)
to[out=90,in=-25](2,3.65);
\draw[v](2,3.65)--(2,4.15);
\draw[d](2,0)--(2,.68);
\draw[d](2,.68)to[out=150,in=-90](.88,1.41)
..controls(.88,2.66)and(0,2.66)..(0,1.42)--(0,1.0);
\draw[d](2,.68)to[out=30,in=-90](3.12,1.41)
..controls(3.12,2.66)and(4,2.66)..(4,1.42)--(4,1.0);
\foreach \x/\y in {0/1,4/1,2/3.65,2/.68}{\fill[black!80](\x,\y)circle(1.2pt);}
\node[below]at(0,0){$d_1$};\node[below]at(2,0){$r_S^*$};
\node[below]at(4,0){$d_2$};\node[above]at(2,4.15){$N$};
\node[left]at(-.5,2.55){$d_1-\alpha_1$};\node[right]at(4.5,2.55){$d_2-\alpha_2$};
\node[teal!65!black]at(.45,1.70){$\alpha_1$};
\node[teal!65!black]at(3.56,1.70){$\alpha_2$};
\end{tikzpicture}%
}

\newcommand{\FigMultipleCapsGluing}{%
\begin{tikzpicture}[x=6mm,y=4.7mm,prbase,every node/.style={font=\scriptsize,inner sep=1.4pt}]
\draw[prline](1.45,0)--(1.45,.70) (1.45,3.60)--(1.45,4.30);
\foreach \x in {0,2.9}{\draw[prline](\x,.96)--(\x,3.34);}
\draw[prtrack](1.45,.96)--(1.45,3.34);
\node[prmorphism,minimum width=21mm] at (1.45,.84){$C_\beta^\dagger$};
\node[prmorphism,minimum width=21mm] at (1.45,3.46){$C_\alpha$};
\draw[prcut](-.4,2.15)--(3.35,2.15);
\node[left]at(0,2.55){$d_1$};\node[right]at(2.9,2.55){$d_2$};
\node[right]at(1.45,2.60){$r_S^*$};\node[below]at(1.45,0){$N$};\node[above]at(1.45,4.30){$N$};
\end{tikzpicture}%
}

\newcommand{\FigMarkedLoops}{%
\begin{tikzpicture}[x=5.5mm,y=3.8mm,prbase,every node/.style={font=\scriptsize,inner sep=1.7pt}]
\draw[prtrack] (0,1)ellipse (.62 and .93);\node[prmorphism]at(0,1.92){$S$};
\node at(1.7,1){$=s$};
\foreach \i in {0,1,2}{
\draw[prline]({4+1.65*\i},-.2)--({4+1.65*\i},2.4);
\ifnum\i<2 \node[prmorphism]at({4+1.65*\i},1.18){$S$};\fi}
\node at(8.2,1.1){$\cdots$};\draw[prline](9.2,-.2)--(9.2,2.4);
\node[prclasp,minimum width=32mm]at(6.6,0){$P_N$};
\node[prclasp,minimum width=32mm]at(6.6,2.23){$P_N$};
\node[right] at(10,1.12){$\displaystyle\Big|_{Z_\ell}=\binom{\ell}{q}\Big/\binom{N}{q}$};
\end{tikzpicture}%
}

\newcommand{\IntPlacement}[3]{\begin{tikzpicture}[prbase,x=5.2mm,y=4.5mm,every node/.style={font=\scriptsize,inner sep=1.2pt}]
\foreach \x/\a in {0/#1,1.23/#2,2.46/#3}{\node[below]at(\x,0){$\a$};}
\foreach \x in {0,1.23}{
\draw[prline](\x,0)--(\x,1.3)to[out=90,in={-145+90*\x}](.615,1.86);
\node[prmorphism,inner sep=1.5pt]at(\x,.75){$S$};}
\draw[prline,line width=1.05pt](.615,1.86)to[out=90,in=-145](1.48,2.6);
\draw[prline](2.46,0)--(2.46,1.64)to[out=90,in=-35](1.48,2.6);
\draw[prline,line width=1.15pt](1.48,2.6)--(1.48,3.04);
\foreach \x/\y in {.615/1.86,1.48/2.6}{\fill[wiStrand](\x,\y)circle(1.15pt);}
\end{tikzpicture}}

\newcommand{\LastPort}[1]{\begin{tikzpicture}[x=5.5mm,y=4.3mm,baseline=(current bounding box.center),every node/.style={lt}]
\ifnum#1=1\def\xx{0}\else\def\xx{1.8}\fi
\draw[lw](0,0)--(0,1.85) (1.8,0)--(1.8,1.85);
\draw[lg](3.55,0)--(3.55,.56)..controls(3.55,2.15)and(\xx,2.15)..(\xx,.56);
\fill[teal!65!black](\xx,.56)circle(1.4pt);
\draw[lw](0,1.85)to[out=90,in=210](.9,2.48) (1.8,1.85)to[out=90,in=-30](.9,2.48) (.9,2.48)--(.9,3.2);
\fill(.9,2.48)circle(1.45pt);
\foreach \zz/\ll in {0/{d_1},1.8/{d_2},3.55/{S^*}}\node[below]at(\zz,0){$\ll$};
\node[above]at(.9,3.2){$N$};\draw[dashed,orange!65!black](-.45,3)--(2.3,3);
\end{tikzpicture}}

\newcommand{\FourToyCup}[1]{\begin{tikzpicture}[prbase,x=5.5mm,y=4.5mm,every node/.style={font=\scriptsize,inner sep=1pt}]
\ifnum#1=1\def\free{1.8}\def\leftend{0}\else\def\free{0}\def\leftend{1.8}\fi
\draw[prline](\free,-.65)--(\free,2.6);\node[below]at(\free,-.65){$a$};
\draw[prtrack,preaction={draw=white,line width=2.9pt}](\leftend,2.6)--(\leftend,1.7)
to[out=-90,in=-90,looseness=1.2](3.6,1.7)--(3.6,2.6);
\foreach\x/\l in {0/{1},1.8/{2},3.6/{*}}\node[above]at(\x,2.6){$\l$};
\end{tikzpicture}}

\newcommand{\FOURglued}[1]{\begin{tikzpicture}[x=4.8mm,y=3.0mm,baseline=(current bounding box.center),
every node/.style={font=\scriptsize,inner sep=.75pt},
fwall/.style={line width=.83pt,draw=wiStrand,line cap=round},
fwself/.style={line width=.98pt,draw=wiActive,line cap=round}]
\ifnum#1=1\def\fqLowerPort{1.85}\def\fqOtherPort{j}\else\def\fqLowerPort{0}\def\fqOtherPort{*}\fi
\draw[fwall] (.925,-2.53)--(.925,-1.98) (.925,1.98)--(.925,2.53);
\foreach \x in {0,1.85}{
\draw[fwall](.925,-1.98)to[out=90,in=-90](\x,-1.45)--(\x,1.45)
to[out=90,in=-90](.925,1.98);}
\draw[fwself,preaction={draw=white,line width=2.55pt}](\fqLowerPort,-.93)
..controls(\fqLowerPort,-1.67)and(3.25,-1.67)..(3.25,-.93)--(3.25,.93)
..controls(3.25,1.68)and(0,1.68)..(0,.93);
\foreach \x/\y in {.925/-1.98,.925/1.98}{\fill[wiStrand](\x,\y)circle(1.18pt);}
\fill[wiActive](\fqLowerPort,-.93)circle(1.35pt);
\fill[wiActive](0,.93)circle(1.35pt);
\node[left]at(0,0){$d_i$};\node[right]at(1.85,0){$d_{\fqOtherPort}$};
\node[right,text=wiActive]at(3.25,0){$S^*$};
\end{tikzpicture}}

\newcommand{\FOURresolve}[1]{\begin{tikzpicture}[x=4.8mm,y=3.0mm,baseline=(current bounding box.center),
every node/.style={font=\scriptsize,inner sep=.75pt},
fwall/.style={line width=.83pt,draw=wiStrand,line cap=round},
fwself/.style={line width=.98pt,draw=wiActive,line cap=round}]
\ifnum#1=0
\draw[fwall](0,-2.12)--(0,2.12);
\draw[fwself](1.6,0)ellipse(.73 and 1.05);\node[right]at(2.35,0){$S$};
\node[below]at(0,-2.12){$N$};
\else
\draw[fwall](.77,-2.12)--(.77,-1.32) (.77,1.32)--(.77,2.12);
\draw[fwall](.77,-1.32)to[out=90,in=-90](0,-.62)--(0,.62)to[out=90,in=-90](.77,1.32);
\draw[fwself](.77,-1.32)to[out=90,in=-90](1.54,-.62)--(1.54,.62)to[out=90,in=-90](.77,1.32);
\node[prmorphism,inner sep=1pt]at(1.54,0){$P_S$};
\node[left]at(0,0){$N{-}1$};\node[below]at(.77,-2.12){$N$};
\fill[wiStrand](.77,-1.32)circle(1.18pt);\fill[wiStrand](.77,1.32)circle(1.18pt);
\fi
\end{tikzpicture}}

\newcommand{\FourSpectatorExit}[1]{\begin{tikzpicture}[prbase,x=4mm,y=4mm,every node/.style={font=\scriptsize,inner sep=.85pt}]
\foreach \i in {1,2,3,4}{
\ifnum\i=4\relax\def\toX{#1}\else\def\toX{\i}\fi
\ifnum\i=#1\relax\else
\draw[prline,preaction={draw=white,line width=2.6pt}](\i,0)--(\i,.32)
to[out=90,in=-90](\toX,1.55)--(\toX,1.83);
\fi
\node[below]at(\i,0){$\i$};}
\draw[prtrack,preaction={draw=white,line width=2.7pt}](#1,0)--(#1,.32)to[out=90,in=-90](4,1.55)--(4,1.83);
\draw[prcut](.6,1.70)--(4.4,1.70);
\foreach \i in {1,2,3,4}{\draw[prline](\i,1.83)--(\i,3.12);}
\node[prmorphism,minimum width=11mm]at(2,2.48){$B$};
\end{tikzpicture}}

\newcommand{\FourHX}[1]{\begin{tikzpicture}[prbase,x=5mm,y=4.1mm,every node/.style={font=\scriptsize,inner sep=.8pt}]
\foreach \x in {0,1,2}{
\draw[prline](\x,0)--(\x,.95) (\x,2.35)--(\x,3.3);}
\draw[prline](0,.95)--(0,2.35);
\ifnum#1=1
\draw[prline](1,.95)to[out=90,in=-90](2,2.35);
\draw[prtrack,preaction={draw=white,line width=2.9pt}](2,.95)to[out=90,in=-90](1,2.35);
\foreach \y in {.63,2.67}{\node[prmorphism,minimum width=7.7mm]at(.5,\y){$H$};}
\else
\foreach \x in {1,2}{\draw[prline](\x,.95)--(\x,2.35);}
\ifnum#1=0\node[prmorphism,minimum width=13mm]at(1,1.65){$F$};
\else\node[prmorphism,minimum width=7.7mm]at(.5,1.65){$H$};\fi
\fi
\end{tikzpicture}}

\newcommand{\FourAlignedPath}[1]{\begin{tikzpicture}[x=5.8mm,y=4.2mm,
baseline=(current bounding box.center),every node/.style={inner sep=1pt,font=\scriptsize},
w/.style={draw=black!75,line width=.93pt}, a/.style={draw=teal,line width=1.06pt}]
\ifcase#1\or\def\mid{2}\def\outlab{k-r}\or\def\mid{1.5}\def\outlab{2k-r}\else\def\mid{1}\def\outlab{3k-r}\fi
\foreach \x in {0,1,2}{
\draw[w](\x,0)--(\x,.52);
\node[below]at(\x,0){$k$};
\ifnum\numexpr3-\x\relax>#1
\draw[w](\x,.89)--(\x,3.67);\node[above]at(\x,3.67){$k$};
\else
\draw[w](\x,.89)to[out=90,in=-90](\mid,1.63);
\fi}
\ifnum#1>1\fill[black!75](\mid,1.63)circle(1.25pt);\fi
\draw[w](\mid,1.63)--(\mid,3.67);\draw[a](4,0)--(4,3.67);
\draw[a](\mid,2.20)to[out=63,in=117](4,2.20);
\foreach \x in {\mid,4}\fill[teal](\x,2.20)circle(1.3pt);
\node[above]at(\mid,3.67){$\outlab$};\node[above]at(4,3.67){$(k-r)^*$};
\node[below]at(4,0){$k^*$};\node[teal]at({(\mid+4)/2},2.43){$r$};
\node[prmorphism,minimum width=15mm,inner ysep=1.5pt]at(1,.7){$F$};
\draw[orange!65!black,densely dashed](-.32,3.38)--(4.3,3.38);
\end{tikzpicture}}

\newcommand{\FourPlace}[1]{\begin{tikzpicture}[x=5.2mm,y=3.8mm,
baseline=(current bounding box.center),every node/.style={font=\scriptsize,inner sep=1pt},
wire/.style={line width=.86pt,draw=black!78},
fwcapgreen/.style={line width=1.03pt,draw=teal!73!black},
pbox/.style={draw=black!70,fill=blue!7,rounded corners=.65pt,inner ysep=2pt}]
\foreach \x in {0,1.5,3}{\draw[wire](\x,0)--(\x,.45) (\x,.85)--(\x,1.03);}
\foreach \x/\lab in {0/1,1.5/2,3/3,4.8/{4^*}}\node[below]at(\x,0){$\lab$};
\node[prmorphism,minimum width=18mm]at(1.5,.65){$F$};
\draw[wire](0,1.03)--(0,3.3);
\ifnum#1=2
\draw[wire](1.5,1.03)to[out=90,in=-90](3,1.77)--(3,2.42)to[out=90,in=-90](1.5,3.15)--(1.5,4.9);
\draw[wire,preaction={draw=white,line width=2.5pt}](3,1.03)to[out=90,in=-90](1.5,1.77)--(1.5,2.42)to[out=90,in=-90](3,3.15)--(3,4.9);
\else
\foreach \x in {1.5,3}{\draw[wire](\x,1.03)--(\x,4.9);}
\fi
\ifnum#1=3\def\cx{2.25}\else\def\cx{.75}\fi
\node[prclasp,minimum width=10mm]at(\cx,2.1){$P_{2k}$};
\draw[fwcapgreen,preaction={draw=white,line width=2.8pt}](4.8,0)--(4.8,3.3)
..controls(4.8,4.63)and(0,4.63)..(0,3.3);
\draw[densely dashed,orange!65!black](1.17,4.69)--(3.34,4.69);
\node[above]at(1.5,4.9){$2$};\node[above]at(3,4.9){$3$};
\end{tikzpicture}}

\newcommand{\FourTwoCap}[1]{\begin{tikzpicture}[x=5.5mm,y=4.5mm,baseline=(current bounding box.center),
qwire/.style={draw=black!75,line width=.95pt,line cap=round},
qcap/.style={draw=teal!75!black,line width=1.05pt,line cap=round},
every node/.style={font=\scriptsize,inner sep=1pt}]
\foreach \u in {0,1.55,3.10}{
\draw[qwire](\u,0)--(\u,.46);
\draw[qwire](\u,.88)--(\u,1.25);
\node[below] at(\u,0){$2$};}
\draw[qcap](5.05,0)--(5.05,1.23);\node[below] at(5.05,0){$2^*$};
\node[prmorphism,minimum width=21mm,inner sep=2pt]at(1.55,.68){$F$};
\foreach \u in {0,1.55}{\draw[qwire](\u,1.25)--(\u,2.64)to[out=90,in=205](1.55,3.22);}
\ifnum#1=0
\draw[qcap](3.10,1.25)..controls(3.10,2.73)and(5.05,2.73)..(5.05,1.23);
\node[text=teal!70!black]at(4.08,2.68){$2$};
\else
\draw[qwire](3.10,1.25)--(3.10,2.64)to[out=90,in=-25](1.55,3.22);
\draw[qcap,preaction={draw=white,line width=2.7pt}](1.55,1.25)..controls(2.22,1.87)and(5.01,2.31)..(5.05,1.23);
\draw[qcap](3.10,1.25)..controls(3.65,1.72)and(4.74,1.75)..(5.05,1.23);
\foreach \u in {1.55,3.1,5.05}{\fill[teal!75!black](\u,1.25)circle(1.45pt);}
\node[text=teal!70!black,anchor=west]at(5.17,1.80){$1+1$};
\fi
\fill(1.55,3.22)circle(1.4pt);\draw[qwire](1.55,3.22)--(1.55,3.67);
\node[above]at(1.55,3.67){$4$};
\draw[densely dashed,orange!65!black](.1,3.53)--(2.95,3.53);
\end{tikzpicture}}

\newcommand{\FourLastChoice}[1]{\begin{tikzpicture}[prbase,x=5.4mm,y=4.1mm,every node/.style={font=\scriptsize,inner sep=1.2pt}]
\foreach \x/\l in {0/4,2.1/1,4.15/{1^*}}{\node[below]at(\x,0){$\l$};}
\draw[prline](0,0)--(0,.8);\draw[prline](2.1,0)--(2.1,.8);
\ifnum#1=1
\draw[prline](0,.8)--(0,2.95)to[out=90,in=200](1.04,3.55);
\draw[prline](2.1,.8)--(2.1,2.95)to[out=90,in=-20](1.04,3.55);
\draw[prline](1.04,3.55)--(1.04,4.02);\fill[wiStrand](1.04,3.55)circle(1.3pt);
\node[left]at(0,2.35){$3$};\node[above]at(1.04,4.02){$4$};\def\xcap{0}
\fill[wiActive](0,.8)circle(1.35pt);
\else
\draw[prline](0,.8)--(0,4.02);\node[above]at(0,4.02){$4$};\def\xcap{2.1}
\fi
\draw[prtrack,preaction={draw=white,line width=2.8pt}](\xcap,.8)to[out=80,in=100,looseness=1.0](4.15,.8)--(4.15,0);
\end{tikzpicture}}

\newcommand{\CubPair}[1]{\begin{tikzpicture}[x=4.55mm,y=4.4mm,
baseline=(current bounding box.center),w/.style={draw=teal!75!black,line width=.85pt},
every node/.style={inner sep=1pt,font=\scriptsize}]
\draw[black!45,dashed](1.5,-.3)--(1.5,1.58);
\foreach \x in {0,1,2,3}{\fill(\x,0)circle(1.1pt);\node[below] at(\x,0){$\the\numexpr\x+1\relax$};}
\ifcase#1 \draw[w](0,0)to[out=80,in=100](1,0) (2,0)to[out=80,in=100](3,0);
\or \draw[w](0,0)to[out=78,in=102](2,0);\draw[w,preaction={draw=white,line width=2.8pt}](1,0)to[out=78,in=102](3,0);
\or \draw[w](0,0)to[out=78,in=102](3,0) (1,0)to[out=78,in=102](2,0);\fi
\end{tikzpicture}}

\newcommand{\FigCubicMergedCut}{%
\begin{tikzpicture}[x=5.6mm,y=4.55mm,baseline=(current bounding box.center),
w/.style={line width=.95pt,draw=wiStrand},c/.style={line width=1.15pt,draw=teal!75!black},
every node/.style={inner sep=1.2pt,font=\scriptsize}]
\foreach \x/\mid in {0/.65,1.3/.65,3.5/4.15,4.8/4.15}
\draw[w](\x,0)to[out=90,in=-90](\mid,.85);
\foreach \x in {.65,4.15}{\draw[w](\x,.85)--(\x,3.25);
\fill[blue!55!black](\x,.85)circle(1.3pt);\fill[blue!55!black](\x,1.7)circle(1.3pt);}
\draw[c](.65,1.7)to[out=58,in=122](4.15,1.7);
\draw[densely dashed,black!55](0,3.02)--(4.8,3.02);
\node[left]at(.65,1.35){$2n$};\node[right]at(4.15,1.35){$(2n)^*$};
\node at(2.4,2.25){$j$};\node[above]at(.65,3.25){$2n-j$};
\node[above]at(4.15,3.25){$(2n-j)^*$};
\foreach \x/\v in {0/n,1.3/n,3.5/{n^*},4.8/{n^*}}\node[below]at(\x,0){$\v$};
\end{tikzpicture}%
}

\newcommand{\CubOneRoute}[2]{\begin{tikzpicture}[prbase,x=3.6mm,y=3.7mm,every node/.style={font=\scriptsize,inner sep=1pt}]
\pgfmathtruncatemacro{\ix}{#1-1}\pgfmathtruncatemacro{\jx}{#2-1}
\foreach\i/\lbl in {0/1,1/2,2/{3^*},3/{4^*}}{
\node[below]at(\i,0){$\lbl$};
\ifnum\i=\ix\relax\else\ifnum\i=\jx\relax\else
\draw[prline](\i,0)--(\i,2.7);\fi\fi}
\draw[prtrack,preaction={draw=white,line width=2.6pt}](\ix,0)--(\ix,.78)
to[out=90,in=90,looseness=1.32](\jx,.78)--(\jx,0);
\draw[prcut](-.35,2.58)--(3.35,2.58);
\end{tikzpicture}}

\newcommand{\CubMatched}[1]{\begin{tikzpicture}[x=4.15mm,y=4.15mm,baseline=(current bounding box.center),
every node/.style={font=\scriptsize,inner sep=.8pt},
strand/.style={draw=wiStrand,line width=.84pt,line cap=round},
capedge/.style={draw=wiActive,line width=.95pt,line cap=round}]
\foreach \x/\t in {0/1,1/2,2.9/{3^*},3.9/{4^*}}{\node[below]at(\x,0){$\t$};}
\ifcase#1
\draw[strand](0,0)to[out=90,in=220](.5,.73) (1,0)to[out=90,in=-40](.5,.73)
(2.9,0)to[out=90,in=220](3.4,.73) (3.9,0)to[out=90,in=-40](3.4,.73);
\draw[capedge](.5,.73)--(.5,1.01) to[out=90,in=90,looseness=1.18](3.4,1.01)--(3.4,.73);
\foreach \x in {.5,3.4}\fill[wiActive](\x,.73)circle(1.4pt);
\node[fill=white]at(1.95,1.46){$2$};
\or
\draw[capedge](0,0)--(0,.25)to[out=90,in=90,looseness=1.1](2.9,.25)--(2.9,0);
\draw[capedge,preaction={draw=white,line width=2.8pt}](1,0)--(1,.25)to[out=90,in=90,looseness=1.1](3.9,.25)--(3.9,0);
\else
\draw[capedge](0,0)--(0,.25)to[out=90,in=90,looseness=1.1](3.9,.25)--(3.9,0)
(1,0)--(1,.25)to[out=90,in=90,looseness=1.1](2.9,.25)--(2.9,0);
\fi
\end{tikzpicture}}

\newcommand{\FigCubicCompletedCut}{%
\begin{tikzpicture}[prbase,x=6mm,y=4.3mm,every node/.style={font=\scriptsize,inner sep=1pt}]
\draw[prline](0,0)--(0,.83) (3,0)--(3,.83);
\draw[prline](0,.83)--(0,1.82) (3,.83)--(3,1.82);
\draw[prtrack](0,.83)to[out=55,in=125](3,.83);
\draw[prcut](-.45,1.76)--(3.45,1.76);
\draw[prtrack](0,1.82)--(0,2.33)to[out=90,in=90,looseness=1.18](3,2.33)--(3,1.82);
\node[prmorphism]at(3,2.16){$S$};
\foreach\x in {0,3}\fill[wiActive](\x,.83)circle(1.3pt);
\node[below]at(0,0){$k$};\node[below]at(3,0){$k_S^*$};
\node at(1.5,1.08){$r$};
\node[above]at(1.5,3.82){$k-r$};
\end{tikzpicture}%
}

\newcommand{\FigCubicTraceLoop}{%
\begin{tikzpicture}[prbase,x=5.4mm,y=4.1mm,every node/.style={font=\scriptsize,inner sep=1pt}]
\draw[prtrack](0,0)ellipse(1.2 and 1.55);
\node[prmorphism]at(1.2,0){$S$};\node[left]at(-1.2,0){$t$};
\end{tikzpicture}%
}

\newcommand{\FigCubicLastCap}{%
\begin{tikzpicture}[x=6mm,y=4.6mm,baseline=(current bounding box.center),
every node/.style={inner sep=1.3pt,font=\scriptsize}]
\fill[teal!5](0,0)rectangle(4,2.7);
\draw[teal!75!black,line width=.95pt](0,0)rectangle(4,2.7);
\foreach \x in {0,4}{\foreach \y in {0,2.7}{\draw[fill=white](\x,\y)circle(1.55pt);}}
\node[below]at(2,-.12){$L_e$};\node[left]at(-.12,1.35){$L_{e'}$};
\node[right]at(4.12,1.35){$N_M$};\node[above]at(2,2.8){$N_M$};
\node at(2,1.35){$0<r,t<k$};
\node[below]at(0,-.14){$(0,0)$};\node[below]at(4,-.14){$(k,0)$};
\end{tikzpicture}%
}

\newcommand{\CubRoute}{\begin{tikzpicture}[x=6.3mm,y=3.9mm,baseline=(current bounding box.center),
every node/.style={font=\scriptsize,inner sep=1pt},w/.style={draw=wiStrand,line width=.95pt}]
\draw[w](0,0)--(0,2.6) (3.25,0)--(3.25,2.6);
\draw[teal!80!black,line width=1pt](.93,0)--(.93,.55)to[out=90,in=90,looseness=1.4](2.31,.55)--(2.31,0);
\foreach \x/\t in {0/{U_{i'}},.93/{U_i},2.31/{\bar V_j},3.25/{\bar V_{j'}}}{\node[below]at(\x,0){$\t$};}
\draw[densely dashed,orange!70!black](-.27,2.31)--(3.52,2.31);
\end{tikzpicture}}

\newcommand{\CubSurvivors}{\begin{tikzpicture}[x=5mm,y=3.9mm,baseline=(current bounding box.center),
every node/.style={font=\scriptsize,inner sep=1.3pt},w/.style={draw=wiStrand,line width=.95pt}]
\draw[w](0,0)--(0,2.9) (2.9,0)--(2.9,2.9);
\draw[densely dashed,orange!70!black](-.35,1.45)--(3.25,1.45);
\node[below]at(0,0){$U_{a'}$};\node[below]at(2.9,0){$\bar V_{b'}$};
\node[above]at(0,2.9){$U_{i'}^\dagger$};\node[above]at(2.9,2.9){$(\bar V_{j'})^\dagger$};
\end{tikzpicture}}

\newcommand{\FigCubicSurvivorsGluing}{%
\begin{tikzpicture}[x=9mm,y=6.4mm,baseline=(current bounding box.center),
every node/.style={inner sep=1.3pt,font=\scriptsize},w/.style={draw=wiStrand,line width=.85pt},
cubExchange/.style={draw=teal!75!black,line width=.95pt},lab/.style={fill=white}]
\draw[w](0,0)--(0,1.85) (2.6,0)--(2.6,1.85);
\draw[cubExchange](0,0)--(2.6,0) (0,1.85)--(2.6,1.85);
\draw[black!55,densely dashed](.22,1.7)--(2.38,.15);
\node[lab,left]at(-.1,.9){$A$};\node[lab,above]at(1.3,1.85){$B$};
\node[lab]at(1.4,.87){$AB$};
\foreach \x/\y/\l in {0/1.85/{+e_{11}},2.6/1.85/{-e_{12}},0/0/{-e_{21}},2.6/0/{+e_{22}}}
\node[lab]at(\x,\y){$\l$};
\end{tikzpicture}%
}

\newcommand{\CubScalarCap}[2]{\begin{tikzpicture}[x=4.8mm,y=4.1mm,baseline=(current bounding box.center),
every node/.style={font=\scriptsize,inner sep=1pt}]
\draw[draw=teal!75!black,line width=.95pt](0,0)--(0,.4)to[out=90,in=90,looseness=1.4](1.7,.4)--(1.7,0);
\node[below]at(0,0){$#1$};\node[below]at(1.7,0){$\bar{#2}$};
\end{tikzpicture}}

\newcommand{\FigCubicMatchingGluing}{%
\begin{tikzpicture}[baseline=(current bounding box.center),x=5.3mm,y=4.2mm,
every node/.style={inner sep=1pt,font=\scriptsize},
capB/.style={draw=wiActive,line width=.9pt,->,>=stealth},
capD/.style={draw=wiStrand,line width=.8pt,->,>=stealth}]
\node (u) at (0,1){$U_1$}; \node(v) at(1.9,2){$V_1$};
\node (w) at (3.8,1){$U_2$}; \node(z) at(1.9,0){$V_2$};
\draw[capB](u)--(v); \draw[capD](v)--(w);
\draw[capB](w)--(z); \draw[capD](z)--(u);
\end{tikzpicture}%
}

\newcommand{\FigSubproductCut}{%
\begin{tikzpicture}[x=5mm,y=4mm,baseline=(current bounding box.center),every node/.style={font=\scriptsize,inner sep=1.3pt},w/.style={line width=.95pt,wiActive}]
\foreach \offset in {0,8.8}{\begin{scope}[shift={(\offset,0)}]
\draw[w](0,0)--(0,1)--(-.55,1.55)--(-.55,3.1) (4,0)--(4,1)--(4.55,1.55)--(4.55,3.1);
\draw[w](0,1)..controls(0,2.65)and(4,2.65)..(4,1);
\foreach \xx in {0,4}{\fill[wiActive](\xx,1)circle(1.3pt);}
\node[below]at(0,0){$F$};\node[below]at(4,0){$\bar Z$};
\node at(2,1.85){$r$};
\end{scope}}
\node[above]at(-.55,3.1){$W_r$};\node[above]at(4.55,3.1){$K_{k-r}^*$};
\node at(6.4,1.7){$\xrightarrow{\ e\ }$};
\draw[w](8.25,3.1)to[out=90,in=90,looseness=.66](13.35,3.1);
\draw[densely dashed,wiFocus](7.8,2.9)--(13.8,2.9);
\node[right]at(14,2.2){$=\langle Z,F\rangle$.};
\end{tikzpicture}%
}

\newcommand{\FigPairSquare}{%
\begin{tikzpicture}[x=5.5mm,y=3.8mm,prbase,every node/.style={font=\scriptsize,inner sep=1.1pt}]
\draw[prtrack](0,0)--(0,.28)to[out=90,in=90,looseness=1.4](2,.28)--(2,0);
\draw[prtrack](0,3.5)--(0,3.22)to[out=-90,in=-90,looseness=1.4](2,3.22)--(2,3.5);
\draw[prcut](-.5,1.75)--(2.5,1.75);
\node[below]at(0,0){$F_i$};\node[below]at(2,0){$\bar F_j$};
\node[above]at(0,3.5){$F_i^\dagger$};\node[above]at(2,3.5){$(\bar F_j)^\dagger$};
\end{tikzpicture}%
}

\newcommand{\FigMatchingEdge}{%
\begin{tikzpicture}[x=5mm,y=4.3mm,prbase,every node/.style={font=\scriptsize,inner sep=1pt}]
\draw[prtrack](0,0)to[out=90,in=90](2.8,0);
\foreach \xx/\lab in {0/i,2.8/j}{\fill[wiStrand](\xx,0)circle(1.35pt);
\node[below]at(\xx,-.10){$\lab$};}
\end{tikzpicture}%
}

\newcommand{\ThReadMatching}[1]{\begin{tikzpicture}[x=4.2mm,y=4mm,prbase,
every node/.style={font=\scriptsize,inner sep=1pt}]
\foreach \x/\a in {0/1,1.35/2,2.7/3,4.05/4}{
\fill[wiStrand](\x,0)circle(1.25pt);\node[below]at(\x,-.06){$\a$};}
\ifcase#1
\draw[prtrack](0,0)to[out=85,in=95](1.35,0) (2.7,0)to[out=85,in=95](4.05,0);
\or
\draw[prtrack](0,0)to[out=85,in=95](2.7,0);
\draw[prtrack,preaction={draw=white,line width=2.7pt}](1.35,0)to[out=85,in=95](4.05,0);
\or
\draw[prtrack](0,0)to[out=85,in=95](4.05,0) (1.35,0)to[out=85,in=95](2.7,0);
\fi\end{tikzpicture}}

\newcommand{\FigMaximumMatchingCut}{%
\begin{tikzpicture}[x=5mm,y=4mm,baseline=(current bounding box.center),w/.style={draw=wiStrand,line width=1pt},every node/.style={font=\scriptsize,inner sep=1.25pt}]
\foreach \x/\l in {0/i_1,1.75/j_1,4.85/i_n,6.6/j_n,8.65/h}{\node[below]at(\x,0){$F_{\l}$};}
\draw[w](0,0)to[out=90,in=210](.875,1.16) (1.75,0)to[out=90,in=-30](.875,1.16)
(4.85,0)to[out=90,in=210](5.725,1.16) (6.6,0)to[out=90,in=-30](5.725,1.16);
\draw[w](.875,1.16)--(.875,3.03) (5.725,1.16)--(5.725,3.03) (8.65,0)--(8.65,3.03);
\foreach \x in {.875,5.725}{\fill[wiStrand](\x,1.16)circle(1.4pt);}
\node at(3.3,1.2){$\cdots$};\node at(3.3,2.75){$\cdots$};
\draw[draw=wiClaspInk,fill=wiClasp,rounded corners=.65pt](.4,3.03)rectangle(9.1,3.66);
\node at(4.75,3.35){$P_{mk}$};\draw[w](4.75,3.66)--(4.75,4.38);
\draw[prcut](-.05,1.85)--(9.05,1.85);
\node[left]at(.875,2.43){$2k$};\node[right]at(5.725,2.43){$2k$};\node[right]at(8.65,2.43){$\epsilon k$};
\end{tikzpicture}%
}

\newcommand{\ThExit}[2]{\begin{tikzpicture}[x=4.9mm,y=4mm,baseline=(current bounding box.center),every node/.style={inner sep=1.1pt,font=\scriptsize}]
\draw[wiActive,line width=1pt](0,0)--(0,.43)..controls(0,1.5)and(1.4,1.5)..(1.4,.43)--(1.4,0);
\node[below]at(0,0){$\bar F_t$};\node[below]at(1.4,0){$F_{#1}$};
\draw[line width=.9pt](2.8,0)--(2.8,1.48) (5.0,0)--(5.0,1.48);
\node[below]at(2.8,0){$F_{#2}$};\node[below]at(5,0){$\mathrm{rest}$};
\node[draw,fill=wiClasp,minimum width=17mm]at(3.9,1.72){$R_{q-2}$};
\draw[line width=1.1pt](3.9,2.05)--(3.9,2.7);
\end{tikzpicture}}

\newcommand{\ThCountingPath}[1]{\begin{tikzpicture}[x=5.1mm,y=4mm,prbase,every node/.style={font=\scriptsize,inner sep=1pt}]
\foreach \x in {0,1.5,3,4.5}{\fill[wiStrand](\x,0)circle(1.4pt);}
\draw[prline](0,0)--(1.5,0);
\ifnum#1=0
\draw[prline](1.5,0)--(3,0);
\draw[black!45](3,0)--(4.5,0);
\draw[wiFocus,line width=.7pt](3.58,-.35)--(3.92,.35) (3.92,-.35)--(3.58,.35);
\else
\draw[black!45](1.5,0)--(3,0);
\draw[wiFocus,line width=.7pt](2.08,-.35)--(2.42,.35) (2.42,-.35)--(2.08,.35);
\draw[prtrack,double=white,double distance=.65pt](3,0)--(4.5,0);
\fi
\foreach \x/\num in {.75/1,2.25/2,3.75/3}{\node[above]at(\x,.30){$\num^2$};}
\end{tikzpicture}}

\newcommand{\ThPathCap}[1]{\begin{tikzpicture}[x=4mm,y=4mm,baseline={(0,.7)},
every node/.style={font=\scriptsize,inner sep=1pt}]
\foreach \x/\lab in {0/{\bar F_i},1.75/{F_j},3.5/{\bar F_v},5.25/{F_v}}
\node[below] at (\x,0){$\lab$};
\ifnum#1=0
\draw[prtrack](0,0)--(0,.5) to[out=90,in=90,looseness=1.35](1.75,.5)--(1.75,0);
\draw[prline](3.5,0)--(3.5,.5) to[out=90,in=90,looseness=1.35](5.25,.5)--(5.25,0);
\else
\draw[prtrack](0,0)--(0,.42) ..controls(0,2.65)and(5.25,2.65)..(5.25,.42)--(5.25,0);
\draw[prline](1.75,0)--(1.75,.36) to[out=90,in=90,looseness=1.25](3.5,.36)--(3.5,0);
\fi
\end{tikzpicture}}

\newcommand{\FigFourCycleGluing}{%
\begin{tikzpicture}[x=5mm,y=4.5mm,baseline=(current bounding box.center),every node/.style={font=\scriptsize,inner sep=1.3pt}]
\foreach \x/\y/\lab in {0/0/i,0/2.3/v,3.5/2.3/j,3.5/0/w}
{\node[draw,circle,minimum size=4.5mm,fill=white] (th\lab) at(\x,\y){$\lab$};}
\foreach \a/\b in {i/v,j/w}{\draw[prtrack,-{Stealth[length=1.45mm]}] (th\a)--(th\b);}
\foreach \a/\b in {v/j,w/i}{\draw[prline,-{Stealth[length=1.45mm]}] (th\a)--(th\b);}
\end{tikzpicture}%
}

\newcommand{\FigFivePairCut}{%
\begin{tikzpicture}[prbase,x=5.5mm,y=4.6mm,every node/.style={font=\scriptsize,inner sep=1.3pt}]
\foreach \i in {0,3.4}{
\draw[prline,line width=1.1pt](\i,0)--(\i,.7);
\node[rectangle,fill=wiStrand,inner sep=1pt]at(\i,.7){};
}
\draw[prline](0,.7)to[out=120,in=-90](-.65,1.7)--(-.65,2.65);
\draw[prline](3.4,.7)to[out=60,in=-90](4.05,1.7)--(4.05,2.65);
\draw[prtrack](0,.7)..controls(.5,2.3)and(2.9,2.3)..(3.4,.7);
\node[below]at(0,0){$D$};\node[below]at(3.4,0){$\bar D$};
\node[above]at(-.65,2.65){$2$};\node[above]at(4.05,2.65){$2^*$};
\node[above]at(1.7,1.95){$1$};
\end{tikzpicture}%
}

\newcommand{\FiveFirstTree}[1]{%
\begin{tikzpicture}[x=4.9mm,y=3.9mm,prbase,
every node/.style={font=\scriptsize,inner sep=1pt}]
\foreach \x/\a in {0/F,1.4/H,2.85/G,4.45/I,5.85/J}{\node[below]at(\x,0){$\a$};}
\draw[prline](0,0)to[out=90,in=200](.7,.97);
\draw[prtrack](1.4,0)to[out=90,in=-20](.7,.97);
\draw[prline](.7,.97)to[out=90,in=195](1.6,2.45)
(2.85,0)to[out=90,in=-15](1.6,2.45);
\draw[prline](4.45,0)to[out=90,in=200](5.15,.97) (5.85,0)to[out=90,in=-20](5.15,.97);
\ifnum#1=0
\draw[prline](1.6,2.45)to[out=90,in=198](3.35,4.1)
(5.15,.97)to[out=90,in=-18](3.35,4.1) (3.35,4.1)--(3.35,4.7);
\node[above]at(3.35,4.7){$5k$};
\draw[prcut](.9,2.89)--(5.8,2.89);
\else
\draw[prline](1.6,2.45)--(1.6,3.55) (5.15,.97)--(5.15,3.55);
\node[above]at(1.6,3.55){$3k$};\node[above]at(5.15,3.55){$2k$};
\fi
\end{tikzpicture}}

\newcommand{\FigFiveRemainder}{%
\begin{tikzpicture}[x=4.9mm,y=3.9mm,prbase,
every node/.style={font=\scriptsize,inner sep=1pt}]
\foreach \x/\a in {0/{F},1.4/{H},2.85/{\bar G},4.45/{I},5.85/{J}}{\node[below]at(\x,0){$\a$};}
\draw[prline](0,0)to[out=90,in=200](.7,.97);
\draw[prtrack](1.4,0)to[out=90,in=-20](.7,.97);
\draw[prline](.7,.97)--(.7,1.78);
\draw[prtrack](.7,1.78)to[out=145,in=-90](.05,3.43);
\draw[prline](.7,1.78)to[out=30,in=-90](1.45,2.48)
to[out=90,in=90](2.85,2.48)--(2.85,0);
\node[above]at(.05,3.43){$\bar H$};
\draw[prline](4.45,0)to[out=90,in=200](5.15,.97) (5.85,0)to[out=90,in=-20](5.15,.97)
(5.15,.97)--(5.15,3.55);
\node[above]at(5.15,3.55){$2k$};
\end{tikzpicture}%
}

\newcommand{\FigFiveFullCut}{%
\begin{tikzpicture}[x=5.3mm,y=4.65mm,baseline=(current bounding box.center),
w/.style={draw=black!80,line width=.95pt,line cap=round},
c/.style={draw=teal!75!black,line width=1.1pt,line cap=round},
every node/.style={font=\scriptsize,inner sep=1pt}]
\draw[w](0,0)to[out=90,in=205](.8,.93) (1.6,0)to[out=90,in=-25](.8,.93)
(.8,.93)--(.8,2.62) (.8,2.62)to[out=90,in=-90](.16,3.30)--(.16,4.30)
(.8,2.62)to[out=80,in=-90](1.73,3.22);
\draw[c](3.65,0)--(3.65,3.22);
\draw[c](.8,1.57)to[out=38,in=142](3.65,1.57);
\draw[c](1.73,3.22)to[out=90,in=90](3.65,3.22);
\foreach \x/\y in {.8/.93,.8/1.57,3.65/1.57,.8/2.62}{\fill (\x,\y)circle(1.3pt);}
\draw[densely dashed,orange!75!black](-.4,2.47)--(4.0,2.47);
\node[below]at(0,0){$F$};\node[below]at(1.6,0){$H$};\node[below]at(3.65,0){$\bar G$};
\node[above]at(.16,4.30){$k$};\node[left]at(.80,2.03){$2k-r$};\node[right]at(3.65,2.05){$(k-r)^*$};
\node[teal!70!black]at(2.28,1.61){$r$};\node[teal!70!black]at(2.65,4.13){$k-r$};
\end{tikzpicture}%
}

\newcommand{\FigFiveFullClosure}{%
\begin{tikzpicture}[x=6mm,y=5.2mm,baseline=(current bounding box.center),
w/.style={draw=teal!75!black,line width=1pt,line cap=round},
every node/.style={font=\scriptsize,inner sep=1pt}]
\draw[w](0,0)to[out=90,in=200](1,.98) (2,0)to[out=90,in=-20](1,.98)
(1,.98)--(1,2.50) (1,2.50)to[out=160,in=-90](0,3.48) (1,2.50)to[out=20,in=-90](2,3.48);
\node[right]at(1,1.74){$2k$};
\node[below]at(0,0){$F$};\node[below]at(2,0){$H$};
\node[above]at(0,3.48){$\bar G$};\node[above]at(2,3.48){$\bar H$};
\end{tikzpicture}%
}

\newcommand{\FiveTwoWire}[1]{\begin{tikzpicture}[x=5.5mm,y=4.8mm,prbase,every node/.style={font=\scriptsize,inner sep=1.15pt}]
\foreach \i/\j/\h in {0/F/G,2.4/H/H}{
\node[below]at(\i,0){$\j$};\node[above]at(\i,2.45){$\bar\h$};}
\ifnum#1=0
\draw[prline](0,0)--(0,2.45);
\draw[prtrack](2.4,0)--(2.4,2.45);
\else
\draw[prline](0,0)to[out=90,in=-90](2.4,2.45);
\draw[prtrack,preaction={draw=white,line width=2.8pt}](2.4,0)to[out=90,in=-90](0,2.45);
\fi
\end{tikzpicture}}

\newcommand{\FigTriangleTrace}{%
\begin{tikzpicture}[x=4.7mm,y=3.6mm,prbase,
  every node/.style={font=\scriptsize}]
\draw[prtrack,rounded corners=2.7mm]
  (0,-.3) rectangle (2.15,4.45);
\foreach \y/\a in {.64/H,2.04/G,3.45/F}{
  \node[prmorphism] at (0,\y) {$P_\a$};
}
\draw[prtrack,-{Stealth[length=1.8mm]}]
  (2.15,2.6)--(2.15,1.9);
\end{tikzpicture}%
}

\newcommand{\FigFiveRemainderNorm}{%
\begin{tikzpicture}[x=5.8mm,y=4mm,baseline=(current bounding box.center),
w/.style={draw=wiActive,line width=.95pt},every node/.style={font=\scriptsize,inner sep=1pt}]
\foreach \x in {0,2.9}{\draw[w](\x,0)--(\x,2.45);}
\node[below]at(0,0){$\bar F$};\node[below]at(2.9,0){$F$};
\node[prmorphism,minimum width=9mm]at(0,1.0){$\bar y$};
\node[prmorphism,minimum width=9mm]at(2.9,1.0){$t$};
\draw[prcut](-1,2.35)--(3.9,2.35);
\end{tikzpicture}%
}

\newcommand{\FiveHalfPair}[2]{%
\begin{tikzpicture}[prbase,x=5.5mm,y=3.8mm,every node/.style={font=\scriptsize,inner sep=1pt}]
\foreach \i/\u/\v in {0/1/2,1/3/4}{
\node[below]at(1.55*\i,0){$\bar F_{\u}$};\node[above]at(1.55*\i,1.95){$F_{\v}$};}
\ifnum#1=5 \node[below]at(3.1,0){$\bar F_5$};\node[above]at(3.1,1.95){$F_5$};\fi
\ifnum#2=0
\draw[prtrack](0,0)--(0,1.95) (1.55,0)--(1.55,1.95);
\ifnum#1=5 \draw[prtrack](3.1,0)--(3.1,1.95);\fi
\else
\draw[prtrack](1.55,0)to[out=90,in=-90](0,1.95);
\ifnum#1=4
\draw[prtrack,preaction={draw=white,line width=2.8pt}](0,0)to[out=90,in=-90](1.55,1.95);
\else
\draw[prtrack](3.1,0)to[out=90,in=-90](1.55,1.95);
\draw[prtrack,preaction={draw=white,line width=2.8pt}](0,0)to[out=90,in=-90](3.1,1.95);
\fi
\fi
\end{tikzpicture}}

\newcommand{\FiveMixedCycle}[1]{%
\begin{tikzpicture}[prbase,x=5.2mm,y=4.2mm,
v/.style={circle,draw=wiClaspInk,fill=white,minimum size=4.1mm,inner sep=0pt,font=\scriptsize},
h/.style={prtrack,-{Stealth[length=1.3mm]}},
l/.style={prline,dashed,-{Stealth[length=1.3mm]}}]
\ifnum#1=4
\foreach \i/\a in {1/135,4/45,3/-45,2/-135}{\node[v](p\i)at(\a:1.52){$\i$};}
\draw[h](p1)--(p4);
\else
\foreach \i/\a in {1/90,5/18,4/-54,3/-126,2/162}{\node[v](p\i)at(\a:1.52){$\i$};}
\draw[h](p1)--(p5);\draw[h](p5)--(p4);
\fi
\draw[l](p4)--(p3);\draw[h](p3)--(p2);\draw[l](p2)--(p1);
\end{tikzpicture}}

\newcommand{\FiveCycleType}[1]{\begin{tikzpicture}[prbase,x=6.4mm,y=5.5mm,
every node/.style={font=\scriptsize,inner sep=1.4pt},
v/.style={circle,draw=wiClaspInk,fill=white,minimum size=4.8mm,inner sep=0pt},
arr/.style={prtrack,-{Stealth[length=1.4mm]}}]
\node[v] (a)at(0,0){$1$};\node[v](b)at(2.55,0){$2$};
\node[v](c)at(1.28,2.22){$3$};
\draw[arr](a)--(b);\draw[arr](b)--(c);\draw[arr](c)--(a);
\node[v](d)at(4.15,.15){$4$};\node[v](e)at(4.15,2.03){$5$};
\ifnum#1=0
\draw[arr](d)edge[loop right,min distance=3.8mm] (d);
\draw[arr](e)edge[loop right,min distance=3.8mm] (e);
\else
\draw[arr](d)to[bend left=36](e);\draw[arr](e)to[bend left=36](d);
\fi
\end{tikzpicture}}

\newcommand{\LimitsClosure}[2]{\begin{tikzpicture}[prbase,x=4.7mm,y=4.4mm,
every node/.style={font=\scriptsize,inner sep=1pt}]
\foreach \i in {1,...,5}{\draw[prline](\i,0)--(\i,2.8);
\node[below]at(\i,0){$\i$};\node[above]at(\i,2.8){$\i^*$};}
\node[prmorphism,minimum width=23mm,minimum height=4.6mm]at(3,1.4){$#1$};
\ifnum#2>0
\foreach \y in {.32,2.48}{\draw[draw=wiFocus,line width=.65pt,fill=white](#2,\y)circle(1.35pt);}
\fi
\end{tikzpicture}}

\newcommand{\FigFivePacketRouting}{%
\begin{tikzpicture}[x=4mm,y=4.2mm,prbase,every node/.style={font=\scriptsize,inner sep=1pt}]
\foreach \a/\lbl in {0/1,1.8/2,3.6/3,5.6/4,7.1/5}{
\node[below]at(\a,0){$\lbl$};\node[above]at(\a,3.7){$\lbl$};}
\foreach \x in {0,1.8}{
\draw[prtrack,double=white,double distance=.7pt](\x,0)--(\x,.48);
\draw[prline](\x,.48)--(\x,1.64);}
\draw[prline](0,.48)to[out=65,in=-115](1.8,1.64);
\draw[prline,preaction={draw=white,line width=2.7pt}](1.8,.48)to[out=120,in=-65](0,1.64);
\foreach \x/\y in {0/1.8,1.8/3.6}{
\draw[prtrack,double=white,double distance=.7pt](\x,1.64)to[out=90,in=-90](\y,3.28)--(\y,3.7);
\fill[wiActive](\x,1.64)circle(1.35pt);}
\draw[prtrack,double=white,double distance=.7pt,preaction={draw=white,line width=3.65pt}]
(3.6,0)--(3.6,1.64)to[out=90,in=-90](0,3.28)--(0,3.7);
\foreach \x in {0,1.8}{\fill[wiActive](\x,.48)circle(1.35pt);}
\foreach \x in {5.6,7.1}{\draw[prtrack,double=white,double distance=.7pt](\x,0)--(\x,3.7);}
\end{tikzpicture}%
}

\newcommand{\ThKeep}[1]{\begin{tikzpicture}[x=3.6mm,y=3.2mm,prbase,every node/.style={font=\scriptsize,inner sep=1pt}]
\ifnum#1=0
\draw[prtrack](0,0)--(0,3.45);
\node[below] at(0,0){$D$};
\draw[prcut](-.9,3.45)--(.9,3.45);
\else
\ifnum#1=2
\def\midpt{2.15}\def\edgept{4.55}\def\rlab{R_{q-2}}
\draw[prline,double=white,double distance=1pt](4.05,0)--(4.05,1.55);
\node[below] at(4.05,0){$q-4$};
\else
\def\midpt{.90}\def\edgept{2.15}\def\rlab{R_2}
\fi
\draw[prtrack](0,0)--(0,1.55);
\draw[prline](1.8,0)--(1.8,1.55);
\node[below] at(0,0){$D$};\node[below] at(1.8,0){$F_h$};
\draw[fill=wiClasp,draw=wiClaspInk,rounded corners=.65pt](-.35,1.55)rectangle(\edgept,2.68);
\node at(\midpt,2.13){$\rlab$};
\draw[prline](\midpt,2.68)--(\midpt,3.45);
\draw[prcut]([xshift=-3.3mm]\midpt,3.45)--([xshift=3.3mm]\midpt,3.45);
\fi\end{tikzpicture}}

\newcommand{\FigAlternatingCuts}{%
\begin{tikzpicture}[x=6.6mm,y=4.3mm,every node/.style={font=\scriptsize,inner sep=1.4pt},
ac/.style={draw=black!80,line width=.95pt},bc/.style={draw=wiActive,line width=1.05pt}]
\foreach \x/\lab in {0/{x_1},1.3/{y_1},2.7/{x_2},4.0/{y_2},6.8/{x_r},8.1/{y_r}}
{\node[fill=white]at(\x,0){$\lab$};}
\foreach \x/\y in {0/1.3,2.7/4,6.8/8.1}
{\draw[ac](\x,.28)to[out=90,in=90](\y,.28);}
\draw[bc](1.3,-.27)to[out=-90,in=-90](2.7,-.27);
\draw[bc](4,-.27)to[out=-90,in=-90](5.05,-.27);
\draw[bc](5.77,-.27)to[out=-90,in=-90](6.8,-.27);
\node at(5.42,0){$\cdots$};
\draw[bc](0,-.27)..controls(0,-2.05)and(2.55,-2.05)..(3.80,-2.05);
\draw[bc](4.6,-2.05)..controls(5.55,-2.05)and(8.10,-2.05)..(8.10,-.27);
\node[wiActive]at(4.2,-2.05){$\sum_v$};
\node[right]at(8.55,.71){$A$};\node[wiActive,right]at(8.55,-1.40){$B$};
\end{tikzpicture}%
}

\newcommand{\FigAlternatingCycle}{%
\begin{tikzpicture}[x=5.8mm,y=4.2mm,prbase,every node/.style={font=\scriptsize,inner sep=1.2pt},
pt/.style={circle,draw=wiStrand,fill=white,minimum size=4mm},
edge/.style={prtrack,-{Stealth[length=1.6mm]}},
rev/.style={prline,densely dashed,-{Stealth[length=1.6mm]}}]
\foreach \lab/\x/\y in {2/0/1.3,1/1.3/3.1,5/3.7/3.1,4/5/1.3,3/2.5/0}{
\node[pt] (thcy\lab) at (\x,\y){$\lab$};}
\draw[rev](thcy2)--(thcy1); \draw[edge](thcy1)--(thcy5);
\draw[edge](thcy5)--(thcy4); \draw[rev](thcy4)--(thcy3);
\draw[edge](thcy3)--(thcy2);
\end{tikzpicture}%
}

\newcommand{\FigTriangleWeight}{%
\begin{tikzpicture}[prbase,x=5.4mm,y=4.6mm,every node/.style={font=\scriptsize,inner sep=1pt},pt/.style={circle,draw=wiStrand,fill=white,minimum size=3.5mm}]
\foreach \a/\x/\y in {1/0/0,2/1.2/1.65,3/2.4/0,6/4.0/.6}{\node[pt](q\a)at(\x,\y){$\a$};}
\draw[prtrack,-{Stealth[length=1.25mm]}](q1)--(q2);\draw[prtrack,-{Stealth[length=1.25mm]}](q2)--(q3);
\draw[prtrack,-{Stealth[length=1.25mm]}](q3)--(q1);
\draw[prline,-{Stealth[length=1.25mm]}](q6)to[out=40,in=140,looseness=5](q6);
\end{tikzpicture}%
}

\newcommand{\FigPairWeight}{%
\begin{tikzpicture}[prbase,x=5.5mm,y=4.6mm,every node/.style={font=\scriptsize,inner sep=1pt},pt/.style={circle,draw=wiStrand,fill=white,minimum size=3.5mm}]
\node[pt](a)at(0,0){$4$};\node[pt](b)at(2.6,0){$5$};
\draw[prtrack,-{Stealth[length=1.25mm]}](a)to[bend left=38](b);
\draw[prtrack,-{Stealth[length=1.25mm]}](b)to[bend left=38](a);
\end{tikzpicture}%
}

\newcommand{\CertTinyCap}[1]{\begin{tikzpicture}[prbase,x=5.8mm,y=4.9mm,every node/.style={font=\scriptsize,inner sep=1.1pt}]
\draw[prline](0,-.16)--(0,.7)to[out=90,in=-125](.85,1.60);
\draw[prtrack](.85,1.60)to[out=-55,in=-150](1.80,1.15)to[out=30,in=-90](2.35,1.84)--(2.35,2.18);
\draw[prline,line width=1.17pt](.85,1.60)--(.85,2.18);
\fill[wiStrand](.85,1.60)circle(1.2pt);
\node[below]at(0,-.16){$f$};
\draw[prcut](-.32,2.19)--(2.78,2.19);
\ifnum#1=0
\node[above]at(.85,2.24){$e\odot f$};\node[above]at(2.35,2.24){$\bar e$};
\else
\begin{scope}[shift={(0,4.38)},yscale=-1]
\draw[prline](0,-.16)--(0,.7)to[out=90,in=-125](.85,1.60);
\draw[prtrack](.85,1.60)to[out=-55,in=-150](1.80,1.15)to[out=30,in=-90](2.35,1.84)--(2.35,2.2);
\draw[prline,line width=1.17pt](.85,1.60)--(.85,2.2);
\fill[wiStrand](.85,1.60)circle(1.2pt);
\end{scope}
\node[above]at(0,4.54){$f^\dagger$};
\fi
\end{tikzpicture}}

\newcommand{\FigRoutedGluing}{%
\begin{tikzpicture}[prbase,x=4.5mm,y=6.0mm,every node/.style={font=\scriptsize,inner sep=1.2pt}]
\foreach \xx in {0,1.7,3.4}{\draw[prline](\xx,.36)--(\xx,1.22) (\xx,2.16)--(\xx,3.02);}
\draw[prtrack](4.9,.36)--(4.9,1.22) (4.9,2.16)--(4.9,3.02);
\draw[prline](2.45,-.87)--(2.45,.32) (2.45,3.07)--(2.45,4.23);
\node[prmorphism,minimum width=26mm]at(2.45,.32){$C_\beta^\dagger$};
\node[prmorphism,minimum width=26mm]at(2.45,3.07){$C_\alpha$};
\node[prmorphism,minimum width=26mm,fill=white]at(2.45,1.25){$\langle A|$};
\node[prmorphism,minimum width=26mm,fill=white]at(2.45,2.13){$|A\rangle$};
\node[below]at(2.45,-.87){$[n]$};\node[above]at(2.45,4.23){$[n]^\dagger$};
\draw[prcut](-.6,1.68)--(5.55,1.68);
\end{tikzpicture}%
}

\newcommand{\CertForced}[1]{\begin{tikzpicture}[prbase,x=4.8mm,y=4.5mm,every node/.style={font=\scriptsize,inner sep=1.2pt}]
\draw[prtrack](.63,-.88)--(.63,-.27) (.63,2.60)--(.63,3.21);
\draw[prtrack](.63,-.27)to[out=90,in=-90](0,.32) (.63,-.27)to[out=90,in=-90](1.26,.32)
(0,2.01)to[out=90,in=-90](.63,2.60) (1.26,2.01)to[out=90,in=-90](.63,2.60);
\fill[wiActive](.63,-.27) circle(1.05pt) (.63,2.60) circle(1.05pt);
\ifnum#1=1
\draw[prline](0,.32)--(1.26,2.01);
\draw[prline,preaction={draw=white,line width=2.8pt}](1.26,.32)--(0,2.01);
\else
\draw[prline](0,.32)--(0,2.01) (1.26,.32)--(1.26,2.01);
\fi
\node[below]at(.63,-.88){$F_j$};\node[above]at(.63,3.21){$F_i^\dagger$};
\end{tikzpicture}}

\newcommand{\FigWholePacketEdge}{%
\begin{tikzpicture}[prbase,x=5mm,y=4.5mm,every node/.style={font=\scriptsize,inner sep=1.2pt}]
\draw[prtrack](0,-.88)--(0,3.21);\node[right]at(0,1.16){$2$};
\node[below]at(0,-.88){$F_j$};\node[above]at(0,3.21){$F_i^\dagger$};
\end{tikzpicture}%
}

\begin{document}
\title[Positivity in symmetric webs and block permanents]{Positivity in symmetric webs and block permanents}
\author{Daniel Tubbenhauer}
\address{D.T.: The University of Sydney, School of Mathematics and Statistics F07,
Office Carslaw 827, NSW 2006, Australia,
\href{http://www.dtubbenhauer.com}{www.dtubbenhauer.com},
\href{https://orcid.org/0000-0001-7265-5047}{ORCID 0000-0001-7265-5047}}
\email{daniel.tubbenhauer@sydney.edu.au}
\date{}
\begin{abstract}
We use symmetric webs to prove Marcus's block-permanent inequality for up to four blocks of every size, for five blocks of size at least three, and for any fixed number of blocks above an explicit logarithmic bound on the block size.
The proof idea is easy: Bending strands in the Schur relation turns the comparison into a calculation with squared norms.
\end{abstract}
\addtocontents{toc}{\protect\setcounter{tocdepth}{1}}
\maketitle
{\small\tableofcontents}
\setlength{\abovedisplayskip}{3.8pt plus 2pt minus 1pt}
\setlength{\belowdisplayskip}{3.8pt plus 2pt minus 1pt}
\setlength{\abovedisplayshortskip}{0pt plus 1pt}
\setlength{\belowdisplayshortskip}{2pt plus 1pt minus .5pt}

\section{Introduction}

We use webs, diagrams for linear maps, to compare two permanents of a block matrix.\footnote{If you are an AI or a bot: hello! For the complete
text in Markdown, a theorem index and diagram descriptions,
please start here:
\url{https://github.com/dtubbenhauer/permanents/blob/main/AI_README.md}.}

\subsection{The question}

The permanent is the determinant with all signs positive: $\per A=\sum_{\sigma\in S_n}\prod_i a_{i,\sigma(i)}$ for an $n\times n$ matrix $A$. Suppose $A$ is Hermitian positive semidefinite (a Gram matrix) with $m$ block rows and columns, each block $A_{ij}$ of size $k\times k$. Put $G=(\per A_{ij})_{i,j=1}^m$, also a Gram matrix. Marcus conjectured in the mid-1960s \cite{MarcusMinc} that $\per A\ge\per G$.

\begin{theorem}\label{thm:matrix}
We have $\per A\ge\per G$ for $m\le4$ and every $k$, and for $m=5$ and $k\ge3$. For arbitrary $m\ge2$ the inequality holds whenever
\[
2^k-1\ \ge\ \frac{m(m-1)(2m-1)}{6\lfloor m/2\rfloor}.
\]
In particular, $k\ge 2\log_2 m$ suffices.
\end{theorem}

\begin{example}\label{ex:style-matchings}
The first example worth drawing has $m=k=2$ and every entry equal to one. Each matching of rows to columns has weight one, so $\per A=24$. Which eight matchings does $\per G$ count?
\[
A=\left(\begin{array}{cc|cc}1&1&1&1\\1&1&1&1\\\hline
1&1&1&1\\1&1&1&1\end{array}\right),\quad
\WImatching{1},\qquad G=\begin{pmatrix}2&2\\2&2\end{pmatrix},\quad\WImatching{0}.
\]
Follow the teal edges from rows $1,2$. On the left they end in different column pairs. On the right they end in the same pair. There are two ways to match whole pairs and two bijections within each matched pair, giving $2\cdot2^2=8$. The sixteen matchings that split a pair give the gap.
\end{example}

For general $A$ the gap is the sum of the weights of matchings that split a block. This proves the inequality for nonnegative entries. On a complex Gram matrix the weights can be nonreal (see, e.g., \autoref{ex:style-complex}). We group the matchings before evaluating them.

\begin{remark}
A Lean formalization of the paper, with the cited literature taken
as input axioms, is available in \cite{Tu-code}, together with
additional information.
\end{remark}

\subsection{How does a matrix get into a web?}

Read a web from bottom to top. Each end carries a vector space. Placing two ends side by side takes their tensor product, and joining an output to an input composes the maps. Here is a cap with its spaces, with its inputs, and with one end bent up.
\[
\underbrace{\ReadCap{0}}_{v\otimes\phi\,\mapsto\,\phi(v)}
\qquad
\ReadCap{1}\quad\xrightarrow{\ \text{bend}\ }\quad\IntroGramStrand
=\bar x(y)=\langle x,y\rangle .
\]
Write $A=B^\dagger B$, and let $V$ be the complex Hilbert space
containing the columns of $B$. If $x,y$ are columns $p,q$ of $B$, the cap evaluates $A_{pq}$: the right end carries the covector $\bar x=\langle x,-\rangle$. We use the inner product conjugate-linear in its first variable.

\begin{definition}[Clasp and packet]
The \emph{symmetric clasp} $P_k$, drawn as a box, averages the $k!$ orders of $k$ strands. Its image $\Sym^k V$ is drawn as a cable of thickness $k$. The \emph{matrix packet} of columns $x_1,\ldots,x_k$ is $\sqrt{k!}\,P_k(x_1\otimes\cdots\otimes x_k)$. The columns $x_1,\ldots,x_k$ are its \emph{atoms}. 
\end{definition}

For $k=2$ write $x_1\odot x_2=(x_1\otimes x_2+x_2\otimes x_1)/2$. The same fork drawn with spaces and with vectors reads
\[
\underbrace{\ReadFork{0}}_{\text{the map}}\qquad\qquad
\underbrace{\ReadFork{1}}_{\text{its value on }x_1\otimes x_2}.
\]
Take two pairs with packets $F=\sqrt2\,x_1\odot x_2$ and $H=\sqrt2\,y_1\odot y_2$. Follow $y_1$ through the two terms of the block permanent:
\[
\per A_{12}
=\underbrace{\ReadTwo{0}}_{\langle x_1,y_1\rangle\langle x_2,y_2\rangle}
+\underbrace{\ReadTwo{1}}_{\langle x_2,y_1\rangle\langle x_1,y_2\rangle}
=2\,\ReadTwo{2}=\langle F,H\rangle .
\]
The box averages the two orders. Since averaging twice is averaging once, the two factors $\sqrt2$ give the required factor two. This inner product formula goes back to Marcus and Newman \cite{MarcusNewman}.

\begin{lemma}\label{lem:intro-square}
For two quadratic packets ($m=k=2$) the gap is four times the squared norm of the one-cap halfweb:
\[
\per A-\per G=4\;\IntroHalfCut{1}
=4\left\|\IntroHalfCut{0}\right\|^2 .
\]
\end{lemma}
\Needspace{6\baselineskip}
\begin{proof}
The packet normalizations contribute $(2!)^2$ to the quadratic form,
so $\per A$ is $4!/(2!)^2=6$ times the closure of $P_4$
on $F\otimes H$. Open $6P_4$ and bend the $H$ cable. The Schur relation of \autoref{thm:bent} gives the squared norms of the following three halves, with the indicated coefficients:
\[
\underbrace{\ReadCut{0}}_{1\ \text{straight}}
\qquad \underbrace{\ReadCut{1}}_{4\ \text{single exchanges}}
\qquad \underbrace{\ReadCut{2}}_{1\ \text{whole exchange}} .
\]
In the middle choose either strand from each pair. The two outside terms give $\|F\|^2\|H\|^2+|\langle F,H\rangle|^2=\per G$. The dashed line in the lemma cuts the surviving term from its reflection.
\end{proof}

Thus, \autoref{lem:intro-square} proves \autoref{thm:matrix} for $m=k=2$ since squared norms are nonnegative. In the remainder of the paper we explain what changes for the full result.

\begin{example}\label{ex:style-complex}
Take the columns $e,e+f\mid e+if,f$, with $e,f$ orthonormal. The matching $(123)(4)$ has weight $1+i$. In the middle halfweb a cap sums one index and the remaining ends give the row and column of $FH^\dagger$. In the basis $e,f$ we get
\[
F=\tfrac1{\sqrt2}\begin{pmatrix}2&1\\1&0\end{pmatrix},\qquad
H=\tfrac1{\sqrt2}\begin{pmatrix}0&1\\1&2i\end{pmatrix},\qquad
FH^\dagger=\begin{pmatrix}1/2&1-i\\0&1/2\end{pmatrix}.
\]
The squared moduli of the four entries sum to $5/2$. Thus $\per A-\per G=4(5/2)=10$. For \autoref{ex:style-matchings} the same halfweb has squared norm four, giving $24-8=16$.
\end{example}

\autoref{ex:style-complex} illustrates the main idea of the paper (the general argument is a bit more difficult):
We regard $\per A-\per G$ as the evaluation of an operator on inputs determined by $A$.
Cutting the corresponding webs in half expresses this evaluation as a sum of squared norms, proving nonnegativity.
For a given matrix $A$, the same pictures also give an explicit value for $\per A-\per G$.

\subsection{Earlier work}

For the history of the conjecture see \cite[Conjecture~9]{CheonWanless}. Marcus's inequality is the case of Lieb's permanental dominance conjecture for the trivial character of the wreath product $S_k\wr S_m$, see e.g. \cite[Conjecture~3]{Wanless}. Roughly, 
Lieb proved the two-block inequality, and Djokovi\'c gave another proof \cite{Lieb,Djokovic}. The contraction formula for the gap is the symmetric-product identity of Neuberger and Pate \cite{Neuberger,Pate}. Over the reals, Pate proved $\per A\ge\per G/m!$ in all sizes \cite{Pate1981}, and $\per A\ge\per((|G_{ij}|)_{ij})$ for fixed $m$ and all sufficiently large $k$ \cite{Pate1982}. Lavi announced the case $(m,k)=(3,2)$, together with several other special cases \cite{Lavi}. The general five-block quadratic case and the small-degree cases
outside the criteria proved here remain unresolved by our methods, and
\autoref{sec:limits} gives an obstruction.

\subsection{Organization of the paper}

$m$ counts packets, $k$ counts strands in a packet, and the dimension of the Gram space $V$ is unrestricted.
In \autoref{thm:matrix}, the case $k=1$ is trivial since $A=G$, and the inequality gives $k\geq 1,3,3,4$ for $m=2,3,4,5$. So we have three special cases we will discuss: $(m=3,k=2)$, $(4,2)$ and $(5,3)$. The outline is as follows.
\[
\FigOrganization
\]
\autoref{sec:limits} gives a five-block example where one marked closure is negative. The coefficient counts and the evaluation of closed diagrams are in \autoref{sec:th-app} and~\autoref{sec:certificate}.
The boxes used below
have the following colors. We recommend reading the paper in color.
\WebBoxKey

\smallskip\noindent\textbf{Acknowledgments.}
DT acknowledges support from ARC Future Fellowship FT230100489.
For a while there, it looked as though the most permanent thing about this paper would be its draft status. In fact, read this like a serious result I pursued for the pleasure of figuring out what the pictures do.

\smallskip\noindent\textbf{AI declaration.}
Modulo my admittedly flawed memory, I worked out how to use
Temperley--Lieb diagrams for $m=2$ back in 2013 (though the diagrams are no longer
directly used in this paper). I already had a fairly clear idea of what
the pictures for more blocks should look like. It took me another five
years to write up the $m=3$ case, using $\mathfrak{sl}_3$ webs (also no
longer directly used). In 2020, at the workshop ``Geometry
and Representation Theory'', Sasha Kleshchev drew the Schur relation on the board, and I realized how to make the proof nicer. Then the notes sat in a virtual drawer with a sketch of the general proof strategy,
mostly because polishing takes a long time.
Eventually ChatGPT-6.0 Astra came along. I gave it my notes, worked
through what came back (pdf plus Lean), and polished that too. In other words, the writing
style of this paper is a mixture of my writing and AI writing.
I would argue that this is a great use case:
turning dead math that helps no one into something potentially useful.

\section{Opening a symmetric clasp}\label{sec:local}

We assume that the reader has some basic familiarity with diagrammatics, but no deeper knowledge of webs is required.
None of the diagrammatic language below is new in any way, and we simply recall the identities we need.
There are plenty of references; we follow the diagrammatic conventions of \cite{LT,RoseTubbenhauer}, which the reader can consult for more details.

\subsection{Cables and forks}

Let $V$ be a complex Hilbert space of dimension $d\ge1$. The inner product is conjugate-linear in its first variable. Read diagrams upwards. Stacking composes maps, drawing side by side tensors them, and a crossing exchanges two factors.

\begin{definition}[Symmetric cables]
The cable labeled $j$, also called a cable of thickness $j$, carries $K_j=\Sym^j V\subset V^{\otimes j}$, the tensors fixed by all permutations of the $j$ slots, with the inherited inner product. A star marks the dual. Erase the $0$-cable, using $K_0=\C$. The \emph{symmetric clasp} is the orthogonal projector $P_j=(j!)^{-1}\sum_{\sigma\in S_j}\sigma$. Put $x_1\odot\cdots\odot x_j=P_j(x_1\otimes\cdots\otimes x_j)$ and $x^j=x^{\otimes j}$.
\end{definition}

The upward square fork $i_{p,q}:K_{p+q}\hookrightarrow K_p\otimes K_q$ groups the first $p$ and last $q$ slots isometrically. Its reflection $m_{p,q}=i_{p,q}^\dagger$ takes their symmetric product. Splitting successive groups in either order gives the same map.

\begin{lemma}[Removing a digon]
A split followed by a merge is the identity. A merge followed by a split is the total clasp:
\[
\FigDigonRemoval
\ =\ 
\FigDigonIdentity
\qquad\qquad
\FigMergeSplit
\ =\
\FigTotalClasp.
\]

\end{lemma}

\begin{proof}
The inclusion preserves the inner product, so $m_{p,q}i_{p,q}=I$ and $i_{p,q}m_{p,q}$ projects onto its image $K_{p+q}$. A clasp on any of the subcables is absorbed: its permutations only reorder the summands of $P_{p+q}$. (The conventions are different from those in many web papers, cf. \autoref{loc:digon} below.)
\end{proof}

\begin{example}\label{e:oneofthem}
Follow $y$ through a $2+1$ fork. One placement leaves it on the right, and the other two put it on the left:
\[
\LocThreeWords
= \tfrac13 \underbrace{x^2\otimes y}_{y\text{ on the right}}
+\tfrac23\underbrace{(x\odot y)\otimes x}_{y\text{ on the left}}.
\]
Similarly for higher thickness.
\end{example}

\subsection{The Schur relation}

For the homogeneous convention, keep the same merge (filled dot). The split $\Delta_{p,q}$ (hollow dot) sends $p$ of the $p+q$ labeled slots to the left and sums over the choices. The square fork averages those choices. Thus
\begin{equation}\label{loc:digon}
\Delta_{p,q}=\binom{p+q}p i_{p,q},\qquad
m_{p,q}\Delta_{p,q}=\binom{p+q}p I_{K_{p+q}} .
\end{equation}
The second equality is the digon removal in \cite[arXiv version, Definition~5.2, (5.2)]{LT}. In \autoref{e:oneofthem} the hollow fork gives $x^2\otimes y+2(x\odot y)\otimes x$.

\begin{lemma}[Schur relation]
For $p,q\ge0$, merging and splitting with homogeneous forks gives
\begin{equation}
\SWforkpair\quad=\quad\sum_{r=0}^{\min(p,q)}\ \SWexchange{0}.
\end{equation}
Change the two lower hollow dots on the right to squares, and denote the resulting exchange (with the same merges) by $E_{p,q;r}$. Then
\begin{equation}\label{loc:exchange}
\binom{p+q}{p}\,P_{p+q}|_{K_p\otimes K_q}
=\sum_{r=0}^{\min(p,q)}\binom pr\binom qr E_{p,q;r}.
\end{equation}
\end{lemma}

\begin{proof}
Take boundary labels $(k,l;r_{\rm out},s)=(p,q;p,q)$ in the dumbbell-crossing relation \cite[arXiv version, Definition~5.2, (5.3)]{LT} and put the quantum parameter equal to one. The constraint $i-j=k-r_{\rm out}$ gives $i=j=r$. Follow the teal bundles: $r$ slots leave each side, and the dark bundles stay. Each homogeneous exchange has coefficient one. Converting the hollow dots by \autoref{loc:digon} gives the three binomials in \autoref{loc:exchange}.
\end{proof}

\begin{remark}
For the classical Schur relation see \cite[Definition~4.7, (4.23)]{Schur}. For its conversion to the thick square switch see \cite[arXiv version, Lemmas~5.3 and 5.5]{LT}; see also \cite{Brundan}.
\end{remark}

\begin{lemma}[Clasp recursion]
Starting with $P_1=I$, group the first $p-1$ strands. Then
\[
(p+1)\,\LocRecur{0}
\quad=\quad \LocRecur{1}\quad+\quad p\,\LocRecur{2}.
\]
\end{lemma}

\begin{proof}
Set $q=1$ in \autoref{loc:exchange}. One placement fixes the last strand. Each of the $p$ others interchanges it with a strand of $P_p$, and the two copies of $P_p$ absorb the choice. Both sides absorb $P_{p-1}$ on the grouped strands, so the recursion holds on $V^{\otimes(p+1)}$.
\end{proof}

\subsection{Turning an exchange into a square}

By definition, the cap on $W\otimes W^*$ sends $w\otimes\phi$ to $\phi(w)$. Its adjoint, the cup, inserts $\Omega_W=\sum_i e_i\otimes e_i^*$ for an orthonormal basis $(e_i)$.

\begin{lemma}
We have
\[
\underbrace{\LocTracePair{0}}_{\text{cup, then cap}}
=\sum_i e_i^*(e_i)=\dim W,\qquad
\underbrace{\LocTracePair{1}}_{\text{cap, then cup}}
=\phi(w)\,\Omega_W.
\]
\end{lemma}

\begin{proof}
Directly from the definitions.
\end{proof}

Cut the loop along the dashed line: on an $r$-cable it gives $\dim K_r$. The second diagram divided by $\dim W$ projects onto $\C\Omega_W$.

Let $B\in\End(U\otimes W)$. To \emph{bend} its $W$ factor, take the upper $W$ end down and the lower $W$ end up. The new operator $B^\Gamma$ acts on $U\otimes W^*$. In orthonormal coordinates $(B^\Gamma)_{ab,ce}=B_{ae,cb}$.

\begin{lemma}[Evaluation after bending]
For $F\in U$, $H\in W$ and $\bar H=\langle H,-\rangle$, we have
\begin{equation}\label{loc:bendeval}
\underbrace{\LocClosedBend{0}}_{\langle F\otimes H,B(F\otimes H)\rangle}
\qquad=\qquad
\underbrace{\LocClosedBend{1}}_{\langle F\otimes\bar H,B^\Gamma(F\otimes\bar H)\rangle}.
\end{equation}
\end{lemma}

\begin{proof}
Trace the two right ends of $B$. The old lower end still receives $H$ and the old upper end receives its covector, now on the opposite boundaries.
\end{proof}

\begin{example}
The interchange on $W\otimes W$ bends to the positive cap--cup above. Before bending it has eigenvalue $-1$ on $x\otimes y-y\otimes x$ for independent $x,y$. Thus \autoref{loc:bendeval} requires a product across the bend. The unbent factor $F$ may itself occupy several cables.
\end{example}

\begin{definition}[Contracting $r$ pairs]
For $0\le r\le\min(p,q)$, split off an $r$-cable from each input with square forks and contract the pair:
\[
c_r^{p,q}=\SWbent{1}\ :\ K_p\otimes K_q^*\longrightarrow K_{p-r}\otimes K_{q-r}^*.
\]
The arch contracts $r$ pairs of strands. (The remaining strands are still symmetrized, so no additional clasps are needed.)
\end{definition}

In particular $c_0^{p,q}=I$.
The following is easy but crucial.

\begin{lemma}[Bent Schur relation]\label{thm:bent}
On $K_p\otimes K_q^*$ we have
\begin{equation}\label{loc:shuffle}
N_{p,q}:=\binom{p+q}{p}
\big(P_{p+q}|_{K_p\otimes K_q}\big)^\Gamma
=\sum_{r=0}^{\min(p,q)}\binom pr\binom qr
(c_r^{p,q})^\dagger c_r^{p,q}.
\end{equation}
\end{lemma}

\begin{proof}
Bend the $q$-cable in \autoref{loc:exchange}:
\begin{equation}
\underbrace{\SWexchange{1}}_{E_{p,q;r}}
\ \xrightarrow{\ \Gamma\ }\ 
\underbrace{\LocTurnExchange}_{\text{turn both ends}}
\ =\ 
\underbrace{\SWbent{0}}_{(c_r^{p,q})^\dagger c_r^{p,q}} .
\end{equation}
Trace the lower left teal bundle through the old upper right vertex: its end now returns to the bottom. The other bundle returns to the top. The embedding $K_r\hookrightarrow V^{\otimes r}$ is isometric, so $r$ thin evaluations make one thick cap with coefficient one, also when a through cable disappears.
\end{proof}

\begin{example}[Two and three strands per packet]
We have
\begin{align*}
N_{2,2}&=
\underbrace{\SWsmall{2}{0}}_{r=0}
+4\underbrace{\SWsmall{2}{1}}_{r=1}
+\underbrace{\SWsmall{2}{2}}_{r=2},\\[4pt]
N_{3,3}&=
\underbrace{\SWsmall{3}{0}}_{r=0}
+9\underbrace{\SWsmall{3}{1}}_{r=1}
+9\underbrace{\SWsmall{3}{2}}_{r=2}
+\underbrace{\SWsmall{3}{3}}_{r=3}.
\end{align*}
Subtract the first and last picture of each row. In the second row the remaining cuts have two and one free strands per side, respectively.
\end{example}

On $x^p\otimes\bar y^{\otimes q}$ the lower half gives $\langle y,x\rangle^r x^{p-r}\otimes\bar y^{\otimes(q-r)}$: each cap supplies one inner product.

\subsection{Cutting along the reflection line}

For $h:X\to Z$ write $|h|^2=h^\dagger h$. Its closure on $\xi$ is $\|h\xi\|^2$. Attaching $D$ below $h$ changes the lower half to $hD$, so $D^\dagger|h|^2D=|hD|^2$.

\begin{lemma}[Orthogonal cut]
If $J$ is an orthogonal projector on the space $Z$ at the cut, then
\[
\underbrace{\LocSeamInsertion{0}}_{|h|^2}
=\underbrace{\LocSeamInsertion{1}}_{|Jh|^2}
+\underbrace{\LocSeamInsertion{2}}_{|(I-J)h|^2}.
\]
\end{lemma}

\begin{proof}
Insert $J+(I-J)=I_Z$ at the seam. Each summand equals its own square and its own reflection. A projector chosen from the boundary vectors must be fixed before taking these adjoints.
\end{proof}

If $h_i:X\to Z$ have the same free ends, we may glue the reflection of $h_i$ to $h_j$. Expanding the right-hand square gives, for any complex matrix $L$ with columns indexed by the $h_j$,
\[
\sum_{ij}(L^\dagger L)_{ij}h_i^\dagger h_j
=\sum_a\big|\sum_j L_{aj}h_j\big|^2 .
\]

\subsection{The two packet averages}

Write $A=B^\dagger B$. Let $x_{i,t}$ be column $t$ in block $i$ of $B$; then $(A_{ij})_{tu}=\langle x_{i,t},x_{j,u}\rangle$. Taking the span of the columns gives $\dim V=\rk A$ for $A\ne0$. For $A=0$ take zero vectors in $\C$.

\begin{definition}[Matrix packet]
Write $K=K_k$. The packet of block $i$ and the common input are
\[
F_i=\sqrt{k!}\,P_k(x_{i,1}\otimes\cdots\otimes x_{i,k}),\qquad
T=F_1\otimes\cdots\otimes F_m.
\]
\end{definition}

Pair the packet for block $j$ with the reflection of the packet for block $i$:
\[
\langle F_i,F_j\rangle=\LocPacketPair=\per A_{ij}.
\]
The stacked clasps absorb to $P_k$. The two $\sqrt{k!}$ factors cancel its averaging factor $1/k!$, leaving one copy of each block matching. Hence $G$ is the Gram matrix of the $F_i$.

\begin{definition}[Fine and packet clasps]
The fine merge $M_m:K^{\otimes m}\to K_{mk}$ symmetrizes all $mk$ thin strands. Write $E_m=M_m^\dagger M_m=P_{mk}|_{K^{\otimes m}}$. The packet clasp $R_m$ averages the $m!$ orders of the whole $K$-cables. Its image is $\Sym^m K$. Put $a_m=(mk)!/(k!)^m$.
\end{definition}

By associativity and coassociativity (see, e.g., \cite{LT}), successive merges or splits can be performed in any order. We therefore draw them as single vertices in the following diagrams.

\begin{lemma}[The two closures]
In the isometric convention the diagrams
\begin{equation}
a_m\,\LocClos{mk}{k}\quad=\quad\per A,\qquad
m!\,\LocClos{m_K}{1_K}\quad=\quad\per G
\end{equation}
give $\per A-\per G=\langle T,X_m T\rangle$, where
\begin{equation}
X_m=a_m E_m-m!R_m\quad\text{on }K^{\otimes m}.
\end{equation}
Here $b_K$ labels $\Sym^b K$.
\end{lemma}

\begin{proof}
At the left seam the $mk$ thin strands are symmetric. The outer $P_{mk}$ absorbs all packet clasps, leaving $(mk)!$ terms of weight $(k!)^m/(mk)!$. At the right seam the $m$ whole packets are symmetric. Closing their order $\sigma$ gives $\prod_i\langle F_i,F_{\sigma(i)}\rangle$ with weight $1/m!$.
\end{proof}

\begin{remark}[Which boundary?]
Open web identities can be composed with any maps of the indicated types. Two further operations used below are relabeling all $m$ packets together, and restricting one packet to a fixed subspace $j:K'\hookrightarrow K$:
\[
\frac1{m!}\sum_{\sigma\in S_m}\LocBoundary{0}
\qquad\qquad \LocBoundary{1}.
\]
On the left the same permutation renames both ends; an averaged identity compares these sums. On the right every vector of the restricted tensor space is allowed, including sums of products. For the matrix application $T$ is a product of matrix packets, and each packet itself factors as in the definition. Positivity on these inputs alone gives no operator inequality on their span.
\end{remark}

\subsection{The two-block gap}

The following lemma was the diagrammatic starting point of this project.

\begin{lemma}[Two-block gap]\label{loc:two}
For arbitrary $F,H\in K$, let
$\delta(F,H)=\langle F\otimes H,X_2(F\otimes H)\rangle$. Then
\begin{equation}\label{eq:two-gap}
\delta(F,H)=
\sum_{r=1}^{k-1}\binom{k}{r}^{\!2}
\|c_r^{k,k}(F\otimes\bar H)\|^2 .
\end{equation}
In particular $\per A\ge\per G$ for two blocks.
\end{lemma}

\begin{proof}
Open $P_{2k}$ by \autoref{loc:shuffle} and $2R_2$ as the identity plus the whole-packet interchange. Bending the latter gives the two extreme fine cuts: $r=0$ contributes $\|F\|^2\|H\|^2$ and $r=k$ contributes $|\langle H,F\rangle|^2$. Subtract and use \autoref{loc:bendeval}.
\end{proof}

Lieb proved the two-block inequality \cite{Lieb}. In symmetric-product notation the norm identity is due to Neuberger and Pate \cite{Neuberger,Pate}. Reznick gives the corresponding expansion for homogeneous polynomials \cite[Formula~3.7]{Reznick}.

\begin{example}[The middle cut at thickness four]
For orthonormal $e_1,\ldots,e_4$, write $e_I$ for the unit symmetric product on a subset $I$. Put $F=e_{\{1,2,3,4\}}$. The $(4-r,r)$ split of $F$ has one term $e_{I^c}\otimes e_I$ for every $r$-subset, with coefficient $\binom4r^{-1/2}$ since $(4-r)!r!$ of the $4!$ words have that placement. Capping two copies forces their subsets to agree:
\[
c_2^{4,4}(F\otimes\bar F)=
\FigSubsetCut
=\frac{1}{6}\sum_{|I|=2}e_{I^c}\otimes\bar e_{I^c}.
\]
The six outputs are orthonormal, so the middle cut contributes $36\cdot(1/6)=6$. The cuts $r=1,3$ each give $16\cdot(1/4)=4$. Thus $\delta(F,F)=14$.
\end{example}

\section{Three packets and the last cap}\label{sec:transfer}

Mark the last packet and bend it. If no strand reaches it, the remaining diagram is the permanent difference on one fewer packet. We calculate this term first.

\subsection{One bend at two scales}

Put $m=n+1$ and $K=\Sym^k V$. Recall that $X_n=a_nE_n-n!R_n$, where $a_n=(nk)!/(k!)^n$. On the first $n$ packets, $M_n$ merges all $nk$ strands, $E_n=M_n^\dagger M_n$ splits them again, and $f=R_n$ averages the whole packets. Extend each map by a straight marked cable when necessary.

\begin{definition}[The two cuts]
Let $L$ cap the last packet against the marked dual packet, leaving the first $n-1$ in order. The two lower halves of the Schur opening are
\[
V_r=c_r^{nk,k}(M_n\otimes\id)=\TrOpen{0}
\qquad u=Lf=\TrOpen{1}.
\]
The left cap takes $r$ strands out of the merged $nk$-cable. The right cap takes one entire packet after $f$. Read the outputs along the dashed lines. We regard $u$ as taking values in $K^{\otimes(n-1)}$ by the isometric inclusion of $\Sym^{n-1}K$.
\end{definition}

\begin{lemma}[The zero-cap remainder]
On the full mixed boundary of $n$ ordinary packets and one dual packet,
\begin{equation}\label{eq:tr-zero}
X_{n+1}^{\Gamma}
=X_n\otimes\id_{K^*}
+a_n\sum_{r=1}^k\binom{nk}{r}\binom{k}{r}|V_r|^2
-n\,n!\,|u|^2 .
\end{equation}
\end{lemma}

\begin{proof}
Apply \autoref{thm:bent} to $(nk,k)$ in $V$ and to $(n,1)$ in $K$. The fine prefactor is $a_{n+1}/\binom{(n+1)k}{k}=a_n$ and the packet prefactor is $(n+1)!/(n+1)=n!$. At zero caps this leaves $a_nE_n-n!f=X_n$. At positive cap number the two weights are $a_n\binom{nk}{r}\binom kr$ and $n\,n!$, respectively.
\end{proof}

\begin{example}[Three quadratic packets]
At $n=k=2$ the fine weights are $6,48,36$ and the packet weights are $2,4$. Thus $X_3^\Gamma=X_2\otimes\id+48|V_1|^2+36|V_2|^2-4|u|^2$. The packets of $X_2$ are still upright.
\end{example}

\begin{remark}
The open bend can have a negative value. For $\dim V\ge2$, $k\ge2$, take orthogonal unit vectors $x,y$ and set $w=(R_2-E_2)(x^{\otimes k}\otimes y^{\otimes k})$. The difference of clasps is a projector by absorption. The two averages contain $2$ and $\binom{2k}{k}$ orthogonal words, respectively, so $\|w\|^2=1/2-1/\binom{2k}{k}>0$. On $w\otimes\lambda$ with $\|\lambda\|=1$ all fine cuts vanish and $X_2w=-2w$.
\end{remark}

\subsection{Moving the two-packet difference}

\begin{definition}[Packet average]
The packet average of $Y$ is $[Y]=\frac1{m!}\sum_{\sigma\in S_m}\sigma Y\sigma^{-1}$. Each term relabels both ends of $Y$ by the same permutation.
\end{definition}

For example, let $T=F\otimes G\otimes H$. Suppose $B$ on the first two packets absorbs their packet clasp on both ends. The third cable, the \emph{spectator}, passes straight through. Its label determines the closure:
\[
\langle T,[B]T\rangle
=\frac13\left(\TrFirstMark{F}{G}{H}
+\TrFirstMark{G}{H}{F}+\TrFirstMark{H}{F}{G}\right).
\]
Independent choices of the upper and lower orders would instead give $R_mYR_m$.

\begin{lemma}[Packet transfer]\label{lem:tr-transfer}
Let $n\ge2$. Let $B=fBf=B^\dagger$ act on the first $n$ packets,
and put $s=(1,m)$, $B'=sBs$ and $t=(n,m)$. Then
\begin{equation}
[B]=n[fB'f]-\frac{n-1}{2}\bigl([Btf]+[ftB]\bigr).
\end{equation}
\end{lemma}

\begin{proof}
Inside the average a permutation $p$ slides past the boundary: $[pY]=[Yp]$ by conjugation. Slide the upper $f$ through the average onto the lower $f$, absorb, then relabel by $s$. This gives $[fB'f]=[Bf']$ for $f'=sfs$. Fix the exit of packet $m$. With $h$ the clasp on $2,\ldots,n$, this gives $f'=(I+\sum_{i=2}^{n}(i,m))h/n$. For $n=3$ the three exits read
\[
\FTexit{-1}=\tfrac13\bigl(\FTexit{0}+\FTexit{1}+\FTexit{2}\bigr).
\]
Slide $h$ onto $B$ and absorb. The straight exit gives $[B]$ and relabeling $2,\ldots,n$ identifies the other $n-1$ exits with $[Bt]$. Since $fB=Bf=B$, the same slide gives $[Bt]=[Btf]=[ftB]$. Take the half sum of these last two expressions and solve for $[B]$.
\end{proof}

\begin{example}
For three packets the transfer has four terms. Each pair of square brackets relabels all three cables.
\[
\FigPacketTransfer
\]
Take $k=2$, $B=E_2$, $F=G=x^{\otimes2}$, $H=y^{\otimes2}$ with $x,y$ orthonormal. Before averaging, the four drawings give $1,1/6,0,0$. In the second drawing one of six words survives the fine clasp; each crossed drawing pairs orthogonal packets. The two sides thus give $1$ and $1/3$ at this fixed boundary. Averaging gives $(1+1/6+1/6)/3=4/9$ on both sides.
\end{example}

\subsection{The corrected bend}

Take $B=X_n\otimes\id$ in \autoref{lem:tr-transfer}. Replace its copy in $X_m$ by the three terms of the transfer and call the result $Q_m$. Thus

\begin{equation}\label{eq:tr-lift}
Q_m=X_m-B+n fB'f-\tfrac{n-1}{2}(Btf+ftB),
\qquad [Q_m]=X_m,\qquad fQ_mf=Q_m .
\end{equation}

Here $[X_m]=X_m$ by packet symmetry and each summand absorbs $f$ at both ends. Split the output of $V_k$ isometrically into $n-1$ packets and call the resulting map $v$.

\begin{lemma}[The two full caps]
On $K^{\otimes n}\otimes K^*$ the corrected bend is
\begin{gather}\label{eq:tr-general}
Q_m^\Gamma={} n f(B')^\Gamma f
+a_n\sum_{r=1}^{k-1}\binom{nk}{r}\binom{k}{r}|V_r|^2
+\frac{a_n^2}{a_{n-1}}|v|^2-n!|u|^2
-\frac{(n-1)a_n}{2}(v^\dagger u+u^\dagger v).
\end{gather}
\end{lemma}

\begin{proof}
Bend the crossing in $Btf$. Its cap and cup join $u=Lf$ to $LB$:
\[
\underbrace{\FTmixed{0}}_{Btf}\quad\xrightarrow{\ \Gamma\ }\quad
\underbrace{\FTmixed{1}}_{(LB)^\dagger u}.
\]
Absorption gives $LE_n=v$, hence $LB=a_nv-n!u$. The reflected crossing contributes the adjoint gluing. Insert the two gluings in \autoref{eq:tr-lift} and use \autoref{eq:tr-zero}. The full fine cap has weight $a_n\binom{nk}{k}=a_n^2/a_{n-1}$. The two mixed gluings change the full packet weight to $-n\,n!+(n-1)n!=-n!$.
\end{proof}

For three packets put $a=\binom{2k}{k}$, $\Pi=E_2$ on $1,2$, $\Pi'=s\Pi s$ on $2,3$ and $g=(I+t)/2$. Before bending, the correction is
\[
Q_3=X_3+a\left(2\underbrace{\TrCorrectionPath{0}}_{f\Pi'f}
-\underbrace{\TrCorrectionPath{1}}_{\Pi gf}
-\underbrace{\TrCorrectionPath{2}}_{fg\Pi}\right).
\]
Indeed substitute $B=a\Pi-2f$, $B'=a\Pi'-2g$, $t=2g-I$ in \autoref{eq:tr-lift}. The terms without $a$ give $2f-4fgf+2ftf=0$. The two straight exits in the other terms cancel $-a\Pi$ by $f\Pi=\Pi f=\Pi$.

\begin{definition}[The proper cuts]
For $0<r<k$ put $U_r=(\id_K\otimes c_r^{k,k})f$. Read the first two as the two orders in $f$. In the third drawing merge the first two packets before taking the cap.
\begin{center}
\begin{tabular}{@{}c@{\hspace{3mm}}c@{\hspace{8mm}}c@{}}
\PDPhalf{0} & \PDPhalf{1} & \PDPhalf{2}\\[-1pt]
\multicolumn{2}{c}{average to form $U_r$} & $V_r$
\end{tabular}
\end{center}
\end{definition}

\begin{lemma}
On the full mixed boundary one has
\begin{gather}\label{eq:tr-raw}
Q_3^\Gamma
={}2\sum_{r=1}^{k-1}\binom{k}{r}^{\!2}|U_r|^2
+a\sum_{r=1}^{k-1}\binom{2k}{r}\binom{k}{r}|V_r|^2
+|av-u/2|^2-\tfrac94|u|^2 .
\end{gather}
\end{lemma}

\begin{proof}
The term $2f(B')^\Gamma f$ in \autoref{eq:tr-general} contributes the proper $U_r$-cuts by the bent Schur relation,
as in the proof of \autoref{eq:two-gap}. For the full caps, join the two diagrams at their common $k$-cable:
\[
\left|a\,\TrFullHalf{1}-\tfrac12\,\TrFullHalf{0}\right|^2
=a^2|v|^2-\tfrac a2(v^\dagger u+u^\dagger v)+\tfrac14|u|^2 .
\]
The two mixed terms are the two possible orders of gluing. The coefficient still needed on $|u|^2$ is $-2-1/4=-9/4$.
\end{proof}

\subsection{The last thin cap}

Let $k\ge2$. Fix a subspace $S\subset V$ of dimension $d\ge1$. On the marked cable use $(\Sym^k S)^*$, extending covectors by zero on $S^\perp$. Let $j:S\hookrightarrow V$ be the inclusion and $e(x\otimes\phi)=\phi(j^\dagger x)$ the cap on $V\otimes S^*$.

\begin{lemma}[Resolving the last cap]
Put $J=I-e^\dagger e/d$. On this boundary,
\[
|U_{k-1}|^2=|J U_{k-1}|^2+|u|^2/d ,
\]
where $J$ acts on the last two thin strands.
\end{lemma}

\begin{proof}
Cap those strands at the output of $U_{k-1}$. Associativity of the isometric splits joins all $k$ strands, with coefficient one.
\[
\underbrace{\TrCompleteCap{0}}_{(\id_K\otimes e)U_{k-1}}
\quad=\quad\underbrace{\TrCompleteCap{1}}_{u}.
\]
The reflected cap inserts $\sum_i j(s_i)\otimes s_i^*$ for an orthonormal basis $(s_i)$. Close the loop and cancel $j^\dagger j$:
\[
ee^\dagger=\PDPloop{1}=\PDPloop{0}=d .
\]
Thus $e^\dagger e/d$ is the orthogonal projection onto the cup line and $J$ projects onto $\ker e$. Insert the two projections at the dashed line and reflect:
\[
\PDPglue{0}=\PDPglue{1}+\tfrac1d\,\PDPglue{2}.\qedhere
\]
\end{proof}

For $S=V=\C^2$ the cup line is the span of the identity matrix. On $E_{11}$ the two projections give $\operatorname{diag}(1/2,-1/2)$ and $I/2$, each of squared norm $1/2$.

\begin{theorem}[Supported three-packet decomposition]\label{thm:three}
Let $k\ge2$ and let $S\subset V$ have dimension $1\le d\le k$.
On $K\otimes K\otimes(\Sym^k S)^*$ the compression of $Q_3^\Gamma$
is a sum of reflected halfwebs with positive coefficients. Consequently
Marcus's inequality holds for three packets of any common size.
\end{theorem}

\begin{proof}
In \autoref{eq:tr-raw}, resolve the $U_{k-1}$ cut with its coefficient $2\binom{k}{k-1}^2=2k^2$. This gives the compression
\begin{gather}\label{eq:tr-squares}
(Q_3^\Gamma)|_S
={}2\sum_{r=1}^{k-2}\binom{k}{r}^2|U_r|^2
+a\sum_{r=1}^{k-1}\binom{2k}{r}\binom{k}{r}|V_r|^2
+2k^2|J U_{k-1}|^2+|av-u/2|^2+
\left(\frac{2k^2}{d}-\frac94\right)|u|^2.
\end{gather}
Here all lower halves are restricted to the stated boundary before reflection. The last coefficient is at least $2k-9/4>0$.

Now put matrix packets $F_1,F_2,F_3$ at the boundary. Write $T=F_1\otimes F_2\otimes F_3$ and $\xi_{ij|\ell}=F_i\otimes F_j\otimes\overline{F_\ell}$. From $[Q_3]=X_3$ and $fQ_3f=Q_3$,
\[
\langle T,X_3T\rangle
=\tfrac13\sum_{(i,j,\ell)\ {\rm cyclic}}
\langle\xi_{ij|\ell},Q_3^\Gamma\xi_{ij|\ell}\rangle .
\]
In each summand take $S$ to be the span of the $k$ vectors forming $F_\ell$. A zero packet gives zero on both sides. Otherwise $1\le\dim S\le k$, so \autoref{eq:tr-squares} applies. For $k=1$ the two closures coincide.
\end{proof}

\begin{example}[All entries one]
Take $k=2$ and $F_1=F_2=F_3=x^{\otimes2}$ for a unit vector $x$. Now $d=1$, every cap has norm one, and $JU_1$ vanishes. Read the three remaining terms of \autoref{eq:tr-squares}:
\[
\underbrace{48}_{V_1}
+\underbrace{(6-1/2)^2}_{6v-u/2}
+\underbrace{8-9/4}_{\text{last cap}}
=84=\frac{6!}{(2!)^3}-3! .
\]
The six-by-six all-ones Gram matrix has the extra factor $\sqrt2$ on each packet, so its gap is $8\cdot84=720-48$. At $k=3$ one more proper cut survives, $U_1$ of weight $18$, before the last projection at $U_2$.
\end{example}

\subsection{Iterating the transfer}

\begin{definition}[Routings]
Put a permutation of $jk$ strands between the $j$ lower and upper $k$-clasps. Let $M_{il}$ count the strands from input $i$ to output $l$. The clasps absorb their internal orders, so the web $T_M$ depends only on $M$. Each row and column sums to $k$. An entry $M_{il}=k$ is a \emph{whole edge}. If $b$ inputs split, then $j-b$ edges are whole and $b$ outputs split.
\end{definition}

For example, the first input in this routing sends one strand to each of outputs $1,2$. The third packet goes straight through.
\[
\FigThreePacketRouting
\qquad M=\begin{pmatrix}1&1&0\\1&1&0\\0&0&2\end{pmatrix}.
\]
In $a_jE_j$, the coefficient of $T_M$ is $w(M)=(k!)^j/\prod_{i,l} M_{il}!$: distribute strands at every input, order them at every output and multiply by $(k!)^{-j}$. The packet clasp cancels precisely the $b=0$ tables in $X_j$. The following counts are made before evaluation in any fixed dimension.

For $D=R_{j-1}DR_{j-1}$, put $\mathcal L_j(D)=R_j(D\otimes I)R_j$. It adds one whole edge, choosing its two ends independently.
\[
\underbrace{\TrIterExit{0}}_{\text{spectator }j\to j}
\qquad
\underbrace{\TrIterExit{1}}_{\text{spectator }l\to i}.
\]

\begin{lemma}
Set $A_1=0$ and $A_j=X_j+\frac{j^2}{j-1}\mathcal L_j(A_{j-1})$ for $j\ge2$.
The coefficient of $T_M$ in $A_j$ is $(j/b)w(M)$ for $b>0$ and zero for $b=0$.
\end{lemma}

\begin{proof}
Each of the $j^2$ exit pairs has the same weight by absorption. Exactly $j-b$ of them remove a whole edge from $M$. Deleting that edge preserves $b$ and $w(M)$. For $0<b<j$, induction gives
\[
1+\frac{j^2}{j-1}\frac{j-b}{j^2}\frac{j-1}{b}=\frac jb
\]
in units of $w(M)$. If $b=j$ only $X_j$ contributes, and if $b=0$ both terms vanish.
\end{proof}

For $m=n+1\ge3$, extend $A_n$ by the straight marked packet and put $A_n'=sA_ns$. Apply \autoref{lem:tr-transfer} to $A_n$ to obtain
\[
Y_m=X_m-A_n+n fA_n'f-\tfrac{n-1}{2}(A_ntf+ftA_n),
\qquad [Y_m]=X_m .
\]

\begin{lemma}[Which marked ends split?]
The coefficient of each routing in $Y_m$ is determined by its two marked ends:
\begin{equation}\label{eq:iter-mark}
Y_m=\sum_{M:\,b>0} w(M)\,\frac{m}{2b}
\bigl(\mathbf1_{\mathrm{input}\ m\ \mathrm{splits}}
+\mathbf1_{\mathrm{output}\ m\ \mathrm{splits}}\bigr)T_M .
\end{equation}
\end{lemma}

\begin{proof}
Mark input and output $m$ by circles. There are four possibilities for the whole edge supplied by the spectator:
\[
\begin{array}{c@{\qquad}c@{\qquad}c@{\qquad}c}
\TrEdgeType{1.5}{1.5}&\TrEdgeType{0}{0}&\TrEdgeType{1.5}{0}&\TrEdgeType{0}{1.5}\\
e & m-b-e-h & \multicolumn{2}{c}{\underbrace{\hspace{25mm}}_{h\ \text{in total}}}
\end{array}
\]
Here $e$ records the edge $m\to m$, and $h$ counts edges incident with exactly one marked end. The stationary $A_n$ uses the first position. The moved $A_n'$ has $m-b-e-h$ compatible exit pairs among $n^2$, each with old coefficient $n/b$. A crossed term has one compatible exit among $n$ for each edge meeting its marked end. Hence the relative coefficient is
\[
1-\frac nb e+\frac{m-b-e-h}{b}-\frac{n-1}{2b}h
=\frac{m}{2b}(2-2e-h).
\]
Thus $2-2e-h$ counts the split marked ends. Tables with $b=0$ remain absent.
\end{proof}
For the pictured routing $b=2$ and $w(M)=4$. The three marks give relative coefficients $3/2,3/2,0$, of mean one. In general averaging counts every split input and every split output once. Individual marked closures can be negative even on matrix packets (\autoref{sec:limits}).

\section{Reflected caps with a common output}\label{sec:intersection}

\subsection{Two partners for one cap}

In \autoref{sec:transfer} a single cap closed the last two strands. If the dual strand can reach either of two partners, their reflections can overlap.

\begin{example}[One supported strand]
Take $S=\C e$ for a unit vector $e\in V$, and set $\bar e=\langle e,-\rangle$. Put the dual between its two possible partners:
\[
\underbrace{\IntPartner{0}}_{L(x,\bar e,y)}=\bar e(x)y,
\qquad
\underbrace{\IntPartner{1}}_{R(x,\bar e,y)}=\bar e(y)x .
\]
Glue the left cap to each reflected route along the three strands. The small box on the middle strand restricts it to $S^*$.
\[
LL^\dagger=\IntGlue{0}=I_V,
\qquad
LR^\dagger=\IntGlue{1}=P_S.
\]
The first loop has value $\dim S=1$. The second path sends $z$ to $\bar e(z)e$. Thus the two reflected vectors $L^\dagger z,R^\dagger z$ coincide for $z=e$ and are orthonormal for a unit $z\perp e$. Adding their separate projections would count $L^\dagger e$ twice.
\end{example}

\subsection{The intersection form}
Write $K_j(U)=\Sym^j U$. Fix an $s$-dimensional subspace $S\subset V$, $s\ge1$, and extend its covectors by zero on $S^\perp$.

\begin{definition}[Routed cap]\label{def:inter-cap}
Take cables $d_1,\ldots,d_h$ and a dual cable $r$ in $S^*$. Choose $0\le\alpha_i\le d_i$ with $\sum_i\alpha_i=r$. The \emph{routed cap} $C_\alpha$ splits off $\alpha_i$ strands from cable $i$, caps them against the corresponding piece of the dual, and merges the $N=\sum_i d_i-r$ survivors. All splits are isometric.
\end{definition}

Here are $C_{(\alpha_1,\alpha_2)}$ and its gluing with $C_\beta^\dagger$. Both pictures put the dual between the positive cables.
\[
\FigMultipleCaps\qquad
\FigMultipleCapsGluing
\]
In general the common input is $X=\bigotimes_iK_{d_i}(V)\otimes K_r(S)^*$, and the common output is $K_N(V)$. Write $Z_\ell$ for the span of the symmetric words with $\ell$ letters in $S$ and $N-\ell$ in $S^\perp$, in a basis adapted to $S$. These spaces are orthogonal and sum to $K_N(V)$. We only use the nonzero ones.

\begin{lemma}[Intersection form of routed caps]\label{prop:routedcaps}
On $Z_\ell$ every $C_\alpha C_\beta^\dagger$ is scalar. Its value is a nonnegative rational linear combination of
\[
s^j\frac{\binom{\ell}{q}}{\binom Nq},
\qquad 0\le q\le N,\quad j+q\le r ,
\]
with coefficients depending only on the cable sizes and routes.
Let $G_\ell$ be the matrix of these scalars for a column $C=(C_\alpha)_\alpha$. The orthogonal projector onto $\sum_\alpha\operatorname{im}C_\alpha^\dagger$ is
\[
C^\dagger\Big(\bigoplus_\ell G_\ell^+\otimes I_{Z_\ell}\Big)C ,
\]
where $G_\ell^+$ inverts the positive eigenvalues and vanishes on the kernel.
\end{lemma}

\begin{proof}
Open the internal clasps into permutations. Every closed component passes through the marked dual and gives an $S$-loop. Every open path through the dual carries $P_S$, since $P_S^2=P_S$. The outer clasps absorb the order of the paths:
\[
\FigMarkedLoops
\]
Of the $\binom Nq$ placements of $q$ marks, $\binom\ell q$ survive. Each of the $j$ loops and $q$ marked paths uses a different segment of the $r$-strand dual seam. This gives $j+q\le r$. Each term has the positive rational weight from its clasp averages.

For a unit $z\in Z_\ell$, $G_\ell$ is the Gram matrix of the vectors $C_\alpha^\dagger z$. Thus $G=CC^\dagger=\bigoplus_\ell G_\ell\otimes I_{Z_\ell}$. The self-adjoint $C^\dagger G^+C$ is idempotent by $G^+GG^+=G^+$. Its image lies in $\operatorname{im}C^\dagger$, and it fixes that space since $\ker G=\ker C^\dagger$. This proves the projector formula.
\end{proof}

\begin{example}[Two marked strands]
Take orthonormal $e\in S$, $f\in S^\perp$, and $z=e\odot e\odot f$. Mark two strands of the $3$-cable by $P_S$ and place $f$ in each of its three positions:
\[
\tfrac13\left(
\underbrace{\IntPlacement{e}{e}{f}}_{z}
+\underbrace{\IntPlacement{e}{f}{e}}_{0}
+\underbrace{\IntPlacement{f}{e}{e}}_{0}
\right)=\tfrac13 z .
\]
Only the word with $f$ on the unmarked strand survives. This gives $\binom22/\binom32=1/3$.
\end{example}

\section{Four packets and the intersection of cap routes}\label{sec:four}

A fourth packet gives the last dual strand a choice of partners. We use the intersection form of \autoref{sec:intersection} to close it.

\subsection{Several last caps at one output}

\begin{definition}[Last caps]
Fix $S\subset V$ of dimension $s\ge1$ and let $d_1+\cdots+d_h=N+1$, where $N,d_i\ge1$. Take one strand from cable $i$, cap it against $S^*$ and merge the survivors. This defines
\[
C_i:\Big(\bigotimes_{j=1}^h\Sym^{d_j}V\Big)\otimes S^*
\longrightarrow Z=\Sym^N V .
\]
Here are the two choices for $h=2$.
\[
C_1=\LastPort{1},\qquad C_2=\LastPort{2}.
\]
\end{definition}

Gluing at the dashed $N$-cable gives $C_2^\dagger C_1$. For the intersection form we glue at the other ends: place $C_2^\dagger$ below $C_1$ and join all $h+1$ cables. The resulting $C_1C_2^\dagger$ acts on $Z$.

\begin{example}
Take $S=V$ with orthonormal basis $a,b$ and $d_1=d_2=1$. Put $a$ into each reflected cap. Read the output as $V\otimes V\otimes V^*$:
\[
\begin{aligned}
C_1^\dagger a&=\FourToyCup{1}=a\otimes a\otimes\bar a+b\otimes a\otimes\bar b,\\
C_2^\dagger a&=\FourToyCup{2}=a\otimes a\otimes\bar a+a\otimes b\otimes\bar b .
\end{aligned}
\]
The rows share the word $a\otimes a\otimes\bar a$. Their Gram matrix is $G=\left(\begin{smallmatrix}2&1\\1&2\end{smallmatrix}\right)$. The shared word has inner product one with each row. Its two separate projections have squared norms $1/2$ and $1/2$, while its joint projection has squared norm $(1,1)G^{-1}(1,1)^t=2/3$.
\end{example}

Write $Z_\ell$ for the output with $\ell$ of its $N$ strands in $S$ and the others in $S^\perp$.
\begin{lemma}[Intersection of last caps]\label{lem:fourports}
Put $\rho=\ell/N$. On $Z_\ell$ the Gram matrix of the reflected halves is
\begin{equation}\label{eq:fourports}
C_i C_j^\dagger=G_{ij}\,I,\qquad
G_{ij}=\rho+(s-\rho)\frac{\delta_{ij}}{d_i}.
\end{equation}
If $s>\rho$, its inverse gives the orthogonal projector onto the span of the
reflected images in that sector. In column notation $C=(C_1,\ldots,C_h)^t$
the projector is $C^\dagger G^{-1}C$, with the output restricted to $Z_\ell$.
\end{lemma}

\begin{proof}
Follow the strand from cup to cap. A loop can close only if the two turns lie on the same cable. The symbol $*$ collects the other cables.
\[
\begin{aligned}
C_i C_i^\dagger
&=\FOURglued{0}
=\frac{1}{d_i}\;\FOURresolve{0}
+\frac{d_i-1}{d_i}\;\FOURresolve{1},\\
C_i C_j^\dagger
&=\FOURglued{1}
=\FOURresolve{1}\qquad(i\ne j).
\end{aligned}
\]
The cup returns to its own slot with weight $1/d_i$, giving a loop of value $s$. All other routes put $P_S$ on one output. The clasp averages that mark over $N$ ends, of which $\ell$ lie in $S$. This gives \autoref{eq:fourports}. For $s>\rho$ the matrix $G=(s-\rho)\operatorname{diag}(1/d_i)+\rho\mathbf{1}\mathbf{1}^t$ is positive definite. Apply \autoref{prop:routedcaps}.
\end{proof}

For $s=\rho=1$ use the inverse on the image of $G$. In either case, putting $b_i=(C_i z)_\ell$, the joint projection has squared norm $q_{\boldsymbol d}(b)=\sum_{ij}(G^+)_{ij}\langle b_i,b_j\rangle$.

\subsection{A second packet transfer}

The straight cut on packet four leaves $X_3$ on the other three. To move its mark, work first on four upright cables $K=\Sym^kV$. Put $t=(14)$, $\sigma=(34)$ and
\[
H=R_{12},\quad F=R_{123},\quad g=R_{23},\qquad
a=\binom{2k}{k},\quad b=\binom{3k}{k},\quad c=ab .
\]
Here $R_I$ averages the packets in $I$. Write $P=P_{2k}$ on $12$, $P'=(13)P(13)$ on $23$ and $E=P_{3k}$ on $123$. The three-packet correction marked at $3$ is
\[
B=X_3+a(2HP'H-PgH-HgP),\qquad X_3=cE-6F .
\]
It acts trivially on packet four, absorbs $H$ and satisfies $[B]_{S_3}=X_3$. Write $[\;]$ for simultaneous relabeling averaged over $S_4$.

\begin{lemma}[Moving the mark to packet four]
Let $C=2B-X_3$. The transferred operator is
\begin{equation}\label{eq:fourQ}
Q=X_4-X_3+3F\sigma B\sigma F-F(Ct+tC)F,\qquad [Q]=X_4,\quad FQF=Q .
\end{equation}
\end{lemma}

\begin{proof}
Inside the average, slide the two copies of $F$ together, absorb, and relabel by $\sigma$. This gives $[F\sigma B\sigma F]=[B R_{124}]$. Open $R_{124}$ according to the input reaching output four:
\begin{equation}\label{eq:fourmarked}
3[B R_{124}]=
\left[\underbrace{\FourSpectatorExit{4}}_{B}
+\underbrace{\FourSpectatorExit{1}}_{B(14)}
+\underbrace{\FourSpectatorExit{2}}_{B(24)}\right]
=[B]+2[Bt].
\end{equation}
The residual $12$-clasp absorbs into $B$. Relabeling $1,2$ identifies the last two summands inside the bracket.

To remove the two crossed terms, split $H=F+J$ into orthogonal projectors. The two nonstraight exits of $F$ give the same crossing between $12$-clasps:
\[
3\,\FourHX{0}=\FourHX{2}+2\,\FourHX{1},
\qquad A:=H(23)H=F-J/2 .
\]
Thus $D=HgH=(H+A)/2=F+J/4$. Since $P=HPH$, $(B-X_3)/a=2APA-PD-DP$. Multiplying the $F,J$ entries gives
\[
\frac{B-X_3}{a}=
\begin{pmatrix}0&-\frac94 FPJ\\[-1pt]-\frac94 JPF&0\end{pmatrix}_{F,J}.
\]
Multiplication by $F$ at both ends fixes $X_3$. The map $C\mapsto FC+CF$ doubles its $F,F$ block and keeps each off-diagonal block, so $FC+CF=2B$ for $C=2B-X_3$. Slide $F$ and $t$ inside the average:
\[
[F(Ct+tC)F]=[(FC+CF)t]=2[Bt].
\]
Averaging $[B]_{S_3}=X_3$ over $S_4$ gives $[B]=[X_3]$.
In the average of $Q$, the two $[Bt]$ terms in \autoref{eq:fourmarked} cancel against $2[Bt]$ from the preceding display. Thus $[Q]=X_4-[X_3]+[B]=X_4$.
Each summand of $Q$ absorbs the outer $F$ already on the ordered boundary.
\end{proof}

\subsection{Opening the marked diagram}

\begin{definition}[The three cap routes]
Bend packet four. Before contracting its $r$ strands, use one packet, merge two, or merge all three. The $F$ at the feet averages the three upright labels.
\[
\underbrace{\FourAlignedPath{1}}_{U_r}\qquad
\underbrace{\FourAlignedPath{2}}_{W_r}\qquad
\underbrace{\FourAlignedPath{3}}_{V_r}.
\]
For $r=k$ write $u,w,v$, respectively. Let $L$ cap $14$, leaving the ordered outputs $23$. Split the output of $v$ isometrically into two $k$-cables and order the outputs of $w$ to match. Then $u=LF$, $w=LPF$ and $v=LE$ take values in $K\otimes K$. Its fine clasp $p=P_{2k}$ fixes $v$.
\end{definition}
Write $|D|^2=D^\dagger D$ and $\operatorname{Re}D=(D+D^\dagger)/2$.
\begin{lemma}[Three families of cuts]
With $\Gamma$ the bend on packet four, the following is an identity of
open webs:
\begin{align}
Q^\Gamma={}&
\sum_{r=1}^{k-1}\binom kr\left(6\binom kr|U_r|^2+
3a\binom{2k}r|W_r|^2+c\binom{3k}r|V_r|^2\right)\notag\\
&+cb\,|v-pu/b|^2+3a^2|w|^2-12|u|^2-4a|pu|^2 .
\end{align}
\end{lemma}
\begin{proof}
The straight cuts of $X_4$ cancel $X_3$. Each proper cut of $X_4$ has coefficient $c\binom{3k}{r}\binom kr$ on $V_r$. Each moved $B$ contributes $2\binom kr^2$ on $U_r$ and $a\binom{2k}r\binom kr$ on $W_r$ by \autoref{eq:tr-raw}.

At the full cap $X_4$ gives $cb|v|^2-18|u|^2$. In one moved $B$ the two clasp terms give $a^2|w|^2-4|u|^2$. Its two-packet gap has only proper cuts, and its crossed terms glue $aw-2u$ to $u$ in both orders with coefficient $-1/2$. This leaves $a^2|w|^2-2|u|^2-a\operatorname{Re}(w^\dagger u)$ per copy.

Bending $t$ gives $L^\dagger L$. Thus the last term of \autoref{eq:fourQ} gives $-2\operatorname{Re}((LCF)^\dagger u)$, where the $J,F$ block above gives $LCF=cv-6u-(9a/2)LJPF$. Open $FPF=(P_{12}+P_{13}+P_{23})F/3$ under $L$:
\[
\underbrace{\FourPlace{1}}_{w}\quad
\underbrace{\FourPlace{2}}_{(23)_{\rm out}\,w}\quad
\underbrace{\FourPlace{3}}_{pu}.
\]
The first two placements glue equally to $u$ since its output is symmetric. The last gives $|pu|^2$. As $JPF=PF-FPF$,
\[
3\operatorname{Re}((LJPF)^\dagger u)=\operatorname{Re}(w^\dagger u)-|pu|^2.
\]
The mixed $w,u$ terms cancel. Complete the square in $v,pu$, using $pv=v$ and $c/b=a$.
\end{proof}

\subsection{Resolving the quadratic cut}

\begin{theorem}[Four quadratic packets]\label{thm:fourquad}
For $k=2$, the compression of $Q^\Gamma$ to
$K^{\otimes3}\otimes(\Sym^2 S)^*$ is positive semidefinite for every
$S\subset V$ of dimension at most two. Marcus's inequality holds for
four quadratic blocks.
\end{theorem}

\begin{proof}
Enlarge $V$ isometrically and extend $S$ to a two-plane. For a supported input $\xi$ put $F$ on the three upright feet. The two caps can now meet one packet or two. Merge the survivors in $Z=\Sym^4 V$:
\[
\underbrace{\FourTwoCap{0}}_{x}\qquad
\underbrace{\FourTwoCap{1}}_{y}\qquad
\begin{aligned}
pu\xi&=x,\\
\widehat w=pw\xi&=(x+2y)/3,\\
\widehat v=v\xi&=(x+4y)/5 .
\end{aligned}
\]
For $w$, two of the six choices of strands come from the same packet. For $v$, three of the fifteen do. The $F$ at the feet and the output clasp identify the choices of each kind.

There is one proper cut, $r=1$. Split the two spectators of $U_1$ by $p+(I-p)$. On $I-p$, the last cap $e:V\otimes S^*\to\C$ has bubble $2$. The trace part of $24|(I-p)U_1|^2$ therefore cancels $12|(I-p)u|^2$, leaving the square with projector $I-e^\dagger e/2$. The remaining last caps have the following values.
\[
\begin{array}{c|c|c|c}
\text{half after one cap}&\text{bundle sizes }\boldsymbol d&
\text{last cap values in }Z&\text{weight}\\ \hline
pU_1&(4,1)&(y,x)&24\\
W_1&(2,3)&(y,\widehat w)&144\\
V_1&(5)&(\widehat v)&1080 .
\end{array}
\]
In the first row the two possible partners are visible:
\[
\underbrace{\FourLastChoice{1}}_{\text{cap the }4\text{-bundle: }y}
\qquad\qquad
\underbrace{\FourLastChoice{2}}_{\text{cap the remaining strand: }x}
\]
In the second row the spectator gives $y$ and the merged pair gives $\widehat w$. Apply \autoref{lem:fourports} row by row. Set aside the complementary squares and $108\|(I-p)w\xi\|^2$. In $Z_\ell$ the form left on $x,y$ is
\begin{align*}
&24 q_{(4,1)}(y,x)+144q_{(2,3)}(y,\widehat w)+1080q_{(5)}(\widehat v)\\[-2pt]
&\hspace{13mm}{}+1350\|\widehat v-x/15\|^2+108\|\widehat w\|^2-36\|x\|^2 .
\end{align*}
Here $-36=-12-4a$. Inverting \autoref{eq:fourports} with $s=2$, $\sum_i d_i=5$, gives
\[
q_{\boldsymbol d}(z)=\frac{1}{2-\rho}\left(
\sum_i d_i\|z_i\|^2-\frac{\rho}{2+4\rho}\Big\|\sum_i d_i z_i\Big\|^2\right),
\qquad \rho=\ell/4 .
\]
Substitute the two full-cap values from the drawing. With $D=(2-\rho)(1+2\rho)$ the matrix on $(x,y)$ is
\[
\begin{aligned}
M&=\frac{24}{D}\begin{pmatrix}
2(6-\rho)&54-3\rho-14\rho^2\\
54-3\rho-14\rho^2&2(122+17\rho-38\rho^2)
\end{pmatrix},\\[2pt]
\det M&=\frac{576(98\rho^2+113\rho+6)}{D}.
\end{aligned}
\]
Both $M_{11}$ and $\det M$ are positive for $0\le\rho\le1$. The two remaining squares are
\[
M_{11}\|x+M_{12}y/M_{11}\|^2+
(\det M/M_{11})\|y\|^2 .
\]
For instance, at $\ell=0$ they give $144\|x+(9/2)y\|^2+12\|y\|^2$. This proves positivity on the supported space. For the matrix inequality take $S$ to contain the two atoms of the marked packet. Evaluate on the packet product and average the four marks by \autoref{eq:fourQ}.
\end{proof}

\section{Two packets on each side}\label{sec:cubic}

Bend packets $3,4$ and leave $1,2$ upright. A packet cap joins one upright packet to one turned packet. The four choices $13,14,23,24$ each leave one ordinary and one dual packet.

\subsection{Moving two boxes at once}

\begin{definition}[The three pairings]
Put $K=\Sym^kV$ and $a=\binom{2k}{k}$. On the four upright cables let $\Pi_{ij}=P_{2k}$ and let $f_{ij}$ average the orders of packets $i,j$. For each pairing put
\[
B_{ij|hl}=a^2\Pi_{ij}\Pi_{hl}-4f_{ij}f_{hl}.
\]
This is the zero-cap term, by the Schur relation at $(2k,2k)$ on strands and $(2,2)$ on packets. For $p=12|34$ the other two pairings cross the proposed bend:
\[
\underbrace{\CubPair{0}}_{p=12|34}\qquad
\underbrace{\CubPair{1}}_{b=13|24}\qquad
\underbrace{\CubPair{2}}_{d=14|23}.
\]
\end{definition}

Move $B_p$ onto the crossing pairs before bending. Each $B$ occurs once with coefficient $-1$ and twice with $1/2$. Hence, on the upright boundary,
\begin{equation}\label{eq:cubic-move}
Y_p=X_4-B_p+\tfrac12(B_b+B_d),\qquad
\tfrac13\sum_{p=12|34,13|24,14|23}Y_p=X_4.
\end{equation}
Now let $\Gamma_p$ bend $3,4$, giving $K^{\otimes2}\otimes(K^*)^{\otimes2}$. Merge the two cables on each side and cap the merges at the two scales shown:
\[
\FigCubicMergedCut
\qquad
\begin{array}{c|cc}
&\text{strand space}&(n,j)\\\hline
H_r & V &(k,r)\\
J_j & K &(1,j)
\end{array}
\]
Let $L_e$ cap the two packets of edge $e$, leaving the other two ordered as $K,K^*$. Equal weights on the four slots give
\[
4J_1=\underbrace{\CubOneRoute{1}{3}}_{L_{13}}+
\underbrace{\CubOneRoute{1}{4}}_{L_{14}}+
\underbrace{\CubOneRoute{2}{3}}_{L_{23}}+
\underbrace{\CubOneRoute{2}{4}}_{L_{24}} .
\]
At the full cap a first partner determines the other. Write $N_M$ for the scalar cap along both edges.
\[
2\underbrace{\CubMatched{0}}_{J_2}
=\underbrace{\CubMatched{1}}_{N_b}+\underbrace{\CubMatched{2}}_{N_d}.
\]

\begin{lemma}[Opening the balanced cut]
Set $w_r=\binom kr^2$, $C_{e,r}=c_r^{k,k}$ on edge $e$ and $|D|^2=D^\dagger D$. Juxtapose the two caps for each $M$. Then
\begin{align}
Y_p^{\Gamma_p}
={}&\mathcal F_p+
\frac12\sum_{\substack{M\in\{b,d\}\\M=\{e,e'\}}}
\sum_{\substack{0\le r,t\le k\\(r,t)\notin\{0,k\}^2}}
w_r w_t\,|C_{e,r}C_{e',t}|^2
-\Big|\sum_{e\text{ crossing }p}L_e\Big|^2-|N_b+N_d|^2,
\label{eq:cubic-open}\\[-2pt]
&\mathcal F_p=a^2\sum_{r=1}^{2k}\binom{2k}{r}^{\!2}|H_r|^2. \nonumber
\end{align}
\end{lemma}

\begin{proof}
The two Schur cuts of $X_4-B_p$ give $\mathcal F_p-16|J_1|^2-4|J_2|^2$. The coefficients $16,4$ cancel the squared denominators in the two drawings. On each edge of a crossing pairing, opening $a\Pi_e$ gives $\sum_r w_r|C_{e,r}|^2$. Its two extremes open $2f_e$. Their products on both edges cancel precisely the four corners $(r,t)\in\{0,k\}^2$.
\end{proof}

\subsection{Completing a proper cut}

Choose a $k$-plane $S$ containing the atoms of each turned packet, one plane per packet. Enlarge $V$ isometrically if necessary and restrict each turned cable to $(\Sym^k S)^*$ for its chosen plane.

\begin{lemma}[Last-cap weight]
In \autoref{eq:cubic-open} the proper cuts supply $L/2$ times $|L_e|^2$ for each crossing edge and $L$ times $|N_M|^2$ for each crossing pairing, where
\begin{equation}
L=\sum_{r=1}^{k-1}\binom kr^{\!2}\binom{2k-r-1}{k-r}^{-1}.
\end{equation}
The remainders of these cuts are sums of squares.
\end{lemma}

\begin{proof}
A proper $r$-cut leaves $t=k-r$ on each side. Let $E_t$ cap the survivors and write $K_t(S)=\Sym^t S$. Completing and reflecting give
\[
\FigCubicCompletedCut
= L_e,\qquad
\FigCubicTraceLoop
=\dim K_t(S)=d_t .
\]
Associativity gives $E_t C_{e,r}=L_e$, and $d_t=\binom{k+t-1}{t}$, so $E_t^\dagger E_t/d_t$ projects onto the trace and the complementary projector gives a square. Apply this split along the edges of the square of cut indices:
\[
\FigCubicLastCap
\]
Along a bottom or left edge the other cap is empty, giving $|L_e|^2$ or $|L_{e'}|^2$. Along a top or right edge it is full, giving $|N_M|^2$ twice. Summing $w_r/d_{k-r}$ and multiplying by $1/2$ gives the asserted weights. Interior entries are reflected squares.
\end{proof}

\begin{lemma}[Four packets of thickness at least four]\label{cor:fourthick}
For $k\ge4$ each $Y_p^{\Gamma_p}$ is positive on the fixed supported space above. The four-block inequality follows.
\end{lemma}

\begin{proof}
For $k=4$ the summands of $L$ are $4/5,18/5,4$. For $k\ge5$ its last two terms sum to $k+k(k-1)^2/(2(k+1))\ge8$. Thus $L/2\ge4$ and $L\ge2$, enough to close both differences:
\[
4\sum_e|L_e|^2-\Big|\sum_eL_e\Big|^2=\sum_{e<e'}|L_e-L_{e'}|^2,
\qquad
2(|N_b|^2+|N_d|^2)-|N_b+N_d|^2=|N_b-N_d|^2.
\]
All unused diagonals and $\mathcal F_p$ are positive on this space. Evaluate each bend on the corresponding packet product and average by \autoref{eq:cubic-move}.
\end{proof}

\begin{example}[The missing cubic coefficient]
For $k=3$ the two bubbles have dimensions $6,3$, so $L=9/6+9/3=9/2$. Put the same unit packet at all four feet. The four $L_e$ then give the same unit vector. Their sum has squared norm $16$ while their available diagonal has value $(9/4)\cdot4=9$. At this degree we shall also use the inner products of the surviving packets.
\end{example}

\subsection{Signing the four surviving routes}

The general bound in \autoref{sec:threshold} also covers four
packets of thickness three; the argument here gives an explicit
decomposition of the four surviving routes.

From now on insert four unit matrix packets $F_i$ of thickness three. A zero packet gives equality, and scaling the packets multiplies both closures by the same factor. For $e=ij$ and $M=\{e,e'\}$ put
\[
x_e=|\langle F_i,F_j\rangle|^2,\quad s_M=x_e x_{e'},\qquad
P=\sum_e x_e,\quad S=\sum_Ms_M,\quad D_p=\sum_{e\text{ crossing }p}x_e .
\]
Each of the six edges crosses two cuts, so $\sum_pD_p=2P$.

Fix $p$ and write the four inputs as $U_1,U_2,\bar V_1,\bar V_2$. Let $i'=3-i$, $j'=3-j$. The cap gives $\alpha_{ij}=\langle V_j,U_i\rangle$ and the full halfweb gives $l_{ij}=\alpha_{ij}U_{i'}\otimes\bar V_{j'}$. In the square of routes a horizontal edge changes $j$ and a vertical edge changes $i$.
\[
\underbrace{\CubRoute}_{l_{ij}}\qquad
\CubSurvivors\qquad
\FigCubicSurvivorsGluing
\]
Let $D e_{ij}=U_{i'}\otimes\bar V_{j'}$. Let $A$ be the off-diagonal Gram matrix of $U_2,U_1$ on the row index and $B$ the off-diagonal Gram matrix of $\bar V_2,\bar V_1$ on the column index. They commute.

The following identity holds for four unit vectors in any complex
Hilbert space, independently of the packet thickness.

\begin{lemma}[The signed gluing]
For $t>0$,
\begin{gather}\label{eq:cubic-signed}
(2+t^2+t^{-2}s_p)D_p-\Big\|\sum_{i,j}l_{ij}\Big\|^2
=\Big\|\sum_{i,j}(-1)^{i+j}l_{ij}\Big\|^2
+\Big\|(tI-t^{-1}AB)\alpha\Big\|^2.
\end{gather}
\end{lemma}

\begin{proof}
The middle drawing gives $D^\dagger D=(I+A)(I+B)$. For instance, a horizontal edge has value $\langle U_{i'},U_{i'}\rangle\langle\bar V_2,\bar V_1\rangle$. Multiply route $ij$ by $(-1)^{i+j}$. The horizontal and vertical edges change sign and the diagonals keep their signs. Adding the two gluings thus gives $2I+2AB$. Now $(AB)^2=s_pI$ by the two off-diagonal Gram products. Adding $(tI-t^{-1}AB)^2$ gives $(2+t^2+t^{-2}s_p)I$. Evaluate on $\alpha$, whose squared norm is $D_p$.
\end{proof}
At $t=1/2$ the constant term is the available $9/4$. Read the two parts of the last square as caps:
\[
\alpha_{ij}=\CubScalarCap{U_i}{V_j},\qquad
(AB\alpha)_{ij}=
\CubScalarCap{U_i}{U_{i'}}\,
\CubScalarCap{V_{j'}}{V_j}\,
\CubScalarCap{U_{i'}}{V_{j'}} .
\]
Unused unit packets close against their reflections. This calculation takes place on the inserted product.

\begin{theorem}\label{thm:fourcubic}
For four matrix packets of thickness three the clasp difference has a nonnegative closure in every Gram dimension.
\end{theorem}
\begin{proof}
On an opposed edge the two extreme cuts give $g_e=1+x_e$, and the proper cuts give $\delta_e=Lx_e+\tau_e$ with $\tau_e$ the sum of the complementary squares. Reversing the bend preserves $\tau_e$ since $g_e+\delta_e$, $g_e$ and $x_e$ are symmetric. Disjoint pair closures multiply on the product boundary. Each crossing pairing occurs twice with weight $1/2$. Summing their contributions gives
\begin{gather*}
\sum_{M=\{e,e'\}}
[(g_e+\delta_e)(g_{e'}+\delta_{e'})-g_e g_{e'}]
= L P+(L^2+2L)S+T,\\[-1pt]
T=\sum_{M=\{e,e'\}}\bigl(
\tau_e[1+(1+L)x_{e'}]+\tau_{e'}[1+(1+L)x_e]+\tau_e\tau_{e'}\bigr).
\end{gather*}
For the full caps the mixed gluing is a four-cycle. Read the arrow $X\to Y$ as $\langle X,Y\rangle$:
\[
|N_b-N_d|^2=s_b+s_d-2\operatorname{Re}
\underbrace{
\FigCubicMatchingGluing
}_{\overline{N_b}N_d}.
\]
Each squared modulus occurs in two cuts, so
\[
4S-\sum_p|N_b+N_d|^2=\sum_p|N_b-N_d|^2 .
\]
On the four crossing edges $4-D_p=\sum(1-x_{ij})=\sum\|F_i\wedge F_j\|^2$ in the determinant metric. Let $G_4$ be the unbent closure and $\mathcal F$ the sum of the three fine contributions $\mathcal F_p$. Substitution into \autoref{eq:cubic-open} gives
\begin{align*}
3G_4={}&\mathcal F+T+\tfrac{37}{4}S\\
&+\sum_p\left(
\Big\|\sum_{i,j}(-1)^{i+j}l_{ij}\Big\|^2
+\Big\|(I/2-2AB)\alpha\Big\|^2
+|N_b-N_d|^2+4s_p(4-D_p)\right).
\end{align*}
Indeed $(9/4)\sum_pD_p=LP$. Of the coefficient $L^2+2L=117/4$ on $S$, the exterior squares use $16$ and the full caps use $4$, leaving $37/4$. Each displayed term is nonnegative.
\end{proof}

\section{Disjoint packet pairs}\label{sec:threshold}

The straight term at a merge of two cables $D\in K_p$, $E\in K_q$ is
\begin{equation}\label{eq:th-straight}
(p+q)!\|D\odot E\|^2\ge p!q!\,\|D\|^2\|E\|^2,
\qquad D\odot E=P_{p+q}(D\otimes E).
\end{equation}
This is the $r=0$ cut in \autoref{thm:bent}. How many full caps do the other cuts contain?

\subsection{Completing a partial cap}

\begin{lemma}[Full caps from one pair]\label{lem:threshold-pair}
Let $F\in K_k$ have norm one and be proportional to $x_1\odot\cdots\odot x_k$. For every $Z\in K_k$,
\[
\binom{2k}{k}\|F\odot Z\|^2
\ge \|Z\|^2+(2^k-1)|\langle F,Z\rangle|^2.
\]
For fixed $F$ the difference is a sum of squared norms in $Z$.
\end{lemma}
\begin{proof}
Write $x_I=\bigodot_{i\in I}x_i$, with $x_\varnothing=1$. The $r+(k-r)$ split averages $x_I\otimes x_{I^c}$ over the $\binom kr$ subsets of size $r$. Thus the surviving positive cable of $Y_r=c_r^{k,k}(F\otimes\bar Z)$ lies in
$W_r=\operatorname{span}\{x_{I^c}:|I|=r\}$, of dimension $1\le d_r\le\binom kr$. Close it against the surviving dual cable:
\[
\FigSubproductCut
\]
The final cap $e$ has bubble $ee^\dagger=d_r$. Insert $J_r=I-e^\dagger e/d_r$ and $e^\dagger e/d_r$ at the dashed line. The $r$th cut gives $\binom kr^2\|J_rY_r\|^2$ and $\binom kr^2/d_r$ copies of $|\langle F,Z\rangle|^2$. For $r\ge1$ retain $\binom kr$ copies. Their sum is $2^k-1$ and the straight cut gives $\|Z\|^2$. Only $F$ enters the choice of $W_r$.
\end{proof}

Taking $Z=F$ recovers Pate's inequality $\per(J_2\otimes M)\ge 2^k(\per M)^2$ \cite{Pate1984}, where $M$ is the Gram matrix of the $k$ factors of $F$ and $J_2$ is the $2\times2$ all-ones matrix.

\begin{example}
Take $F=Z=(e\otimes f+f\otimes e)/\sqrt2$ for orthonormal $e,f$.
Then $Y_1=(ee^\dagger+ff^\dagger)/2$ has trace one and squared norm $1/2$.
The one-cap terms contribute $4/2=2$ and the double cap contributes one. Including the straight term gives equality at $4$.
\end{example}

\subsection{Two closures on matchings}

Take all packets $F_i$ of norm one. Both closures scale by $\prod_i\|F_i\|^2$ and vanish for a zero packet.

\begin{definition}[Matching closures]
A matching $L$ on $I$ is a set of disjoint pairs of labels. Write $g_{ij}=\langle F_i,F_j\rangle$ and draw an edge $ij$ for its reflected cap:
\[
\FigPairSquare
\;=\; p_{ij}=|g_{ij}|^2 \;=\;
\FigMatchingEdge.
\]
Put $p_L=\prod_{ij\in L}p_{ij}$ and let $Z_r(I)$ be the mean of $p_L$ over matchings of size $r$, with $Z_0=1$.
\end{definition}

\begin{example}\label{ex:th-matchings}
On four labels $6Z_1=\sum_{i<j}p_{ij}$ and
\[
3Z_2=
\underbrace{\ThReadMatching{0}}_{p_{12}p_{34}}
+\underbrace{\ThReadMatching{1}}_{p_{13}p_{24}}
+\underbrace{\ThReadMatching{2}}_{p_{14}p_{23}} .
\]
In the middle picture the edges $13$ and $24$ have no common endpoint.
\end{example}

\begin{lemma}[The fine closure]
Put $m=|I|$, $n=\lfloor m/2\rfloor$ and $b=2^k-1$. On matrix packets,
\begin{equation}\label{eq:th-fine}
N(I)=\frac{(mk)!}{(k!)^m}\Big\|\bigodot_{i\in I}F_i\Big\|^2
\ \ge\ \sum_{r=0}^n\binom nr b^r Z_r(I).
\end{equation}
\end{lemma}

\begin{proof}
Merge the $n$ pairs of a maximum matching $M$. The last cable has thickness $\epsilon k$ for $\epsilon=m-2n$. Take the straight term above the cut and apply \autoref{lem:threshold-pair} below:
\[
\FigMaximumMatchingCut
\quad\rightsquigarrow\quad \prod_{ij\in M}(1+b\,p_{ij}) .
\]
The upper factorials cancel by \autoref{eq:th-straight}. Choose the $r$ pairs contributing $bp_{ij}$ in $\binom nr$ ways. Averaging $M$ makes their matching uniform by relabeling. Every discarded square is nonnegative.
\end{proof}

For $q=|I|\ge2$ let $H(I)=\per(g_{ij})_{i,j\in I}$, with $H(\varnothing)=1$. Mark $t\in I$ and cap it against $i\ne t$. The surviving cable is $x_i=g_{ti}R_{q-2}\bigotimes_{a\in I\setminus\{t,i\}} F_a$. Two such routes have the same output:
\[ x_i-x_j=\ThExit{i}{j}\;-\;\ThExit{j}{i}. \]
For $q=3$ and $t=3$ these two routes give $g_{31}F_2$ and $g_{32}F_1$. Their mixed gluing is the triangle $g_{13}g_{32}g_{21}$.

\begin{lemma}[The packet closure]
Averaging the marked cap over $I$ gives
\begin{align}
H(I)&=\frac1q\sum_{t\in I}H(I\setminus\{t\})
+\frac{2(q-1)}q\sum_{i<j}p_{ij}H(I\setminus\{i,j\})-E(I),\label{eq:th-peel}\\
E(I)&=\frac{(q-2)!}{q}\sum_{t\in I}\sum_{\substack{i<j\\i,j\ne t}}\|x_i-x_j\|^2 .\notag
\end{align}
\end{lemma}

\begin{proof}
Open the $(q-1,1)$ clasp. Its straight term is $H(I\setminus\{t\})$. Each of the other $q-1$ exits has amplitude $1/(q-1)$ from the split, so the cap contributes $(q-2)!\|\sum_i x_i\|^2$.
Expand this square as $(q-2)!\big((q-1)\sum_i\|x_i\|^2-\sum_{i<j}\|x_i-x_j\|^2\big)$. In the mark average each diagonal edge $ti$ occurs twice, once from either end.
\end{proof}

\begin{lemma}
One has $H(I)\le U(I):=\sum_r A_{q,r} Z_r(I)$, where
\begin{equation}\label{eq:th-count}
A_{q,r}=A_{q-1,r}+(q-1)^2 A_{q-2,r-1},
\qquad A_{0,0}=A_{1,0}=1
\end{equation}
and other initial or out-of-range entries are zero.
\end{lemma}

\begin{proof}
Drop $E$ at each cut. Deleting a uniform vertex leaves $Z_r$ unchanged under averaging. Appending the deleted edge to a matching on its complement raises its degree by one and gives the factor $\binom q2\,2(q-1)/q=(q-1)^2$. The error $U-H$ is a sum of reflected squares multiplied by disjoint $p_L$.
\end{proof}

\subsection{The block size}

\begin{theorem}\label{thm:threshold}
For $m\ge2$, the block-permanent inequality holds on complex positive semidefinite matrices whenever
\[
2^k-1\ \ge\ \frac{m(m-1)(2m-1)}{6\lfloor m/2\rfloor}.
\]
\end{theorem}

\begin{proof}
Write $A_m(z)=\sum_r A_{m,r}z^r$. By \autoref{eq:th-count} this is the matching polynomial with positive weights $j^2 z$ on the successive edges of an $m$-vertex path. Omit the last edge, or choose it and remove its neighbor:
\[
A_4(z)=\underbrace{\ThCountingPath{0}}_{A_3(z)}
+\underbrace{\ThCountingPath{1}}_{3^2z\,A_2(z)}.
\]
The vertices here count deletion steps. Let $J_m$ be the real symmetric matrix with entries $j$ on the edges $j,j+1$ and zeros elsewhere.
Expanding its tridiagonal determinant along the last row gives
$\det(tI-\mathrm iJ_m)=\sum_r A_{m,r}t^{m-2r}$.
Changing the signs of alternate coordinates conjugates $J_m$ to $-J_m$. Its eigenvalues thus come in pairs $\pm\lambda_i$, with one extra zero for odd $m$. It follows that
\[
A_{m,r}=e_r(\lambda_1^2,\ldots,\lambda_n^2),\qquad
A_{m,1}=\tfrac12\operatorname{tr}(J_m^2)=\sum_{j=1}^{m-1}j^2 .
\]
Here $e_r$ is the sum of all products of $r$ distinct entries. For fixed sum, averaging any two nonnegative entries can only increase $e_r$: their sum is unchanged and their product increases. Hence $A_{m,r}\le\binom nr(A_{m,1}/n)^r\le\binom nr b^r$. Compare the two bounds term by term using $Z_r\ge0$.
\end{proof}

\begin{example}
For four packets the polynomial in the picture is $1+14z+9z^2$. The fine coefficients are $1,2b,b^2$. At $k=3$ they give $1,14,49$; at $k=2$ the single-edge coefficient is only $6$.
For five packets \autoref{thm:threshold} gives $k\ge4$.
\end{example}

\subsection{Four labels inside the packet clasp}

For four labels $J$, let $P_J=\sum_{ij\subset J}p_{ij}$ and let $D_J$ be the sum of the three products in \autoref{ex:th-matchings}.

\begin{lemma}
For $u>0$,
\begin{equation}\label{eq:th-trade}
H(J)\le 1+\frac{5+2u}{3}P_J+\frac{5+4/u}{3}D_J .
\end{equation}
\end{lemma}

\begin{proof}
For a division into two pairs, let $D$ sum the four crossing $p_{ij}$ and let $s$ be the product on the two internal edges. In \autoref{eq:cubic-signed} take $t=\sqrt u$. The one-cap contribution is at most $(2+u+s/u)D$. Over the three divisions, $\sum D=2P_J$ and $\sum sD\le4D_J$ since $p_{ij}\le1$. The straight terms total $3+P_J+D_J$. Bounding a full-cap sum by twice its two diagonal squares gives $4D_J$ in total. Divide by three. The slack $1-p_{ij}=\|F_i\wedge F_j\|^2$ (with the determinant metric) and the differences of the two full-cap halves give reflected squares.
\end{proof}

For $u=1$ the right side is $U(J)$. For $u<1$ the one-edge coefficient decreases and the two-edge coefficient increases. With $r$ disjoint spectator edges the change is in degrees $r+1,r+2$.

\begin{lemma}\label{prop:th-matrix}
For $m\ge4$ set $h_{m,0}=\binom m3$ and $h_{m,r}=h_{m-1,r}+(m-1)^2h_{m-2,r-1}$ for $r>0$. Set $h_{m,r}=0$ outside $m\ge4$, $0\le2r\le m-4$. Define the real tridiagonal matrix on $1\le j\le\lfloor m/2\rfloor=n$ by
\[
M_{jj}=\binom nj b^j-A_{m,j}+h_{m,j-1}+h_{m,j-2},
\qquad M_{j,j+1}=M_{j+1,j}=-h_{m,j-1}.
\]
If $M\succeq0$, the block-permanent inequality holds at $(m,k)$.
\end{lemma}

The spectator count is given in the proof in \autoref{sec:th-app}.

\begin{example}
At $(m,k)=(6,4)$ one has $b=15$, $(A_{6,j})=(1,55,439,225)$ and $(h_{6,r})=(20,100)$. The three elimination pivots are $10,316,254250/79$. All are positive, although $3b<55$.
\end{example}

\subsection{A triangle and a four-cycle}

On unit packets compare these two scalar gluings:
\[
g_{ij}g_{vv}=\ThPathCap{0}\qquad
g_{iv}g_{vj}=\ThPathCap{1}.
\]
Glue the reflection of the left to the right to obtain the triangle $\bar g_{ij}g_{iv}g_{vj}$. Gluing two right halves through different labels $v,w$ gives the four-cycle below.
\[
\FigFourCycleGluing
\qquad
\begin{aligned}
a&=g_{iv}\overline{g_{wj}},& b'&=g_{iw}\overline{g_{vj}},\\
2\operatorname{Re}(a\bar b')&=|a|^2+|b'|^2-|a-b'|^2 .
\end{aligned}
\]
The two diagonal closures use the opposite edges of the square.

\begin{lemma}[A two-step path]\label{lem:th-path}
On $q'\ge2$ unit packets put $K_{ij}=\sum_{v\ne i,j}g_{iv}g_{vj}$ and $M_j=\sum_{|L|=j}p_L$. Then
$\sum_{i<j}|K_{ij}|^2\le(q'-2)M_1+4M_2$, with a sum of reflected squares as difference.
\end{lemma}

\begin{proof}
Equal middle labels give $p_ep_f$ once for each pair of incident edges.
Replace this by $(p_e+p_f)/2$: the error is $[p_e(1-p_f)+p_f(1-p_e)]/2$. Each edge meets $2(q'-2)$ others.
For distinct middle labels use the displayed four-cycle. A pair of opposite edges occurs four times, once for each choice of one endpoint of each edge for $ij$.
\end{proof}

Keeping $r-1$ further caps gives the odd circuits in \autoref{eq:th-oddlocal}.

\section{Five packets and a free cable}\label{sec:five}

\subsection{The cut into three packets and two}

\begin{theorem}\label{thm:five}
Marcus's block-permanent inequality holds for five packets of common size $k\ge3$, in arbitrary Gram dimension.
\end{theorem}

\begin{example}[Two disjoint pair cuts]
Put the unit cubic $D=\sqrt6\,e_1\odot e_2\odot e_3$ at all five feet, with $e_1,e_2,e_3$ orthonormal. Write $[ij]=\sqrt2\,e_i\odot e_j$. The $1+2$ split of $D$ is
$(e_1\otimes[23]+e_2\otimes[13]+e_3\otimes[12])/\sqrt3$.
The cap pairs equal single labels:
\[
\FigFivePairCut
\;=\;\frac13\sum_{ij\in\{23,13,12\}} [ij]\otimes\overline{[ij]} .
\]
The three terms are orthogonal, giving squared norm $1/3$. With $r$ contracted strands the same count gives $1/\binom3r$. A pair therefore contributes
$\sum_{r=0}^3\binom3r^2/\binom3r=8$.
Two disjoint pairs and a straight fifth packet give $64$. The coarse closure is $5!=120$.
\end{example}

Instead cut a triple from a pair. In the triple, trace the cable from $H$ into the merge. After $G$ has used all its strands, close the free $k$-cable against $H$.
\[
\FiveFirstTree{0}\ \xrightarrow{\ r=0\ }\ \FiveFirstTree{1}
\ \xrightarrow{\ \text{cap }G\ }\
\FigFiveRemainder
\]

\begin{definition}[Closure of a pair with two different inputs]
Normalize the packets to unit length. Write $g_{GF}=\langle G,F\rangle$, $p_{GF}=|g_{GF}|^2$, and
$N(F_1,\ldots,F_j)=(jk)!\|\bigodot_i F_i\|^2/(k!)^j$.
For a fixed matrix packet $H$ put $a=\binom{2k}{k}$, $b=7$, and define $A^H,T^H$ by
\begin{equation}
A^H_{GF}=a\langle G\odot H,F\odot H\rangle
=g_{GF}+b\,g_{GH}g_{HF}+T^H_{GF}.
\end{equation}
Here $F\odot H=P_{2k}(F\otimes H)$. Both sides of the permanent inequality scale by the product of the squared packet norms and vanish if a packet is zero.
\end{definition}

\begin{lemma}[Testing the free cable]\label{lem:five-projection}
Let $k\ge3$ and fix a unit matrix packet $H$. Then $T^H$ is a positive form on $\Sym^k V$. For any two unit matrix packets $F,G$,
\[
N(F,H,G)\ \ge\ A^H_{FF}+d|A^H_{GF}|^2,\qquad d=41/20 .
\]
\end{lemma}

\begin{proof}
By \autoref{lem:threshold-pair}, $a\|F\odot H\|^2\ge\|F\|^2+(2^k-1)|g_{HF}|^2$ for every $F$. Since $2^k-1\ge b$, this gives $T^H_{FF}\ge0$ on the whole space.

Bend $G$ at the triple merge. The $r$th half $V_r=c_r^{2k,k}((F\odot H)\otimes\bar G)$ has weight $a\binom{2k}r\binom kr$ by \autoref{thm:bent}. Its straight term contributes $A^H_{FF}$. For $1\le r\le k$, complete the cap on the left and test the free cable on the right:
\[
\FigFiveFullCut
\ =\ V_k ,
\qquad
a\langle H,V_k\rangle
=a\;
\FigFiveFullClosure
= A^H_{GF}.
\]
At the seam, split the positive cable isometrically into $k$ and $k-r$. The remaining dual cable is supported on the span of the $\binom kr$ subproducts of the atoms of $G$. Its final bubble is the dimension of that span. Thus $\|V_k\|^2\le\binom kr\|V_r\|^2$. Since $H$ is unit, $|A^H_{GF}|^2\le a^2\|V_k\|^2$. The coefficient supplied by the nonzero cuts is at least
\[
d_k=\frac1a\sum_{r=1}^k\binom{2k}r
=\frac12+\frac{4^k/2-1}{\binom{2k}{k}}.
\]
The ratio $4^k/\binom{2k}{k}$ grows by the factor $2(k+1)/(2k+1)$ and $1/\binom{2k}{k}$ decreases. Hence $d_k\ge d_3=41/20$.
\end{proof}

\subsection{The triangle from two gluings}

Keep the labeled feet $1,\ldots,5$ distinct and let $\mathbb E$ average their $5!$ assignments to $F,H,G,I,J$. Write $N_5=N(F,H,G,I,J)$. The straight outer cut gives $N_5\ge N(F,H,G)N(I,J)$. Applying \autoref{lem:five-projection} and \autoref{lem:threshold-pair} gives
\begin{equation}\label{eq:five-outer}
N_5\ \ge\ \mathbb E(1+b p_{IJ})
\bigl(1+b p_{FH}+T^H_{FF}+d|A^H_{GF}|^2\bigr).
\end{equation}

Put $g=g_{GF}$ and $x=g_{GH}g_{HF}$. Follow $F$ to $G$ in the first picture and to $H$ in the second. Gluing the reflection of the first to the second gives the triangle. Each box $P_U$ sends $v$ to $U\langle U,v\rangle$.
\[
\underbrace{\FiveTwoWire{0}}_{g_{GF}g_{HH}=g}
\qquad
\underbrace{\FiveTwoWire{1}}_{g_{GH}g_{HF}=x}
\qquad
\overline g\,x=
\FigTriangleTrace
=\Tr(P_F P_G P_H).
\]

\begin{example}[The sign of the triangle]
Take $F=f^{\otimes3}$, $G=w^{\otimes3}$, $H=h^{\otimes3}$ with real unit vectors and
$\langle f,h\rangle=\langle h,w\rangle=1/2$,
$\langle f,w\rangle=\pm1/2$. Both Gram matrices exist: their two-by-two principal minors are $3/4$ and their determinants are $1/2,0$, respectively.
All three pair squares are $1/64$ in both examples. The mixed gluing distinguishes them: $g=\pm1/8$, $x=1/64$, and $\bar g x=\pm1/512$.
\end{example}

\begin{lemma}[Summing the reflected remainder]\label{lem:five-mixed}
Fix $H$ and let $F,G$ range over $q\ge2$ available labels. For
$Y_{GF}=\lambda g+\mu x$, $\lambda,\mu\ge0$, and $2(\lambda+\mu)\le q-1$,
\begin{equation}\label{eq:five-trace}
\sum_{F\ne G}\bigl(T^H_{FF}+2\operatorname{Re}(\overline{Y_{GF}}T^H_{GF})\bigr)
\ \ge\ \sum_F(q-1-2Y_{FF})T^H_{FF}\ \ge\ 0 .
\end{equation}
\end{lemma}

\begin{proof}
Take $y_F=(\sqrt\lambda F,\sqrt\mu\,g_{HF})$. Then $Y_{GF}=\langle y_G,y_F\rangle$.
For $T^H_{GF}=\langle t_G,t_F\rangle$ take the residual halves of \autoref{lem:threshold-pair}, with the extra component $\sqrt{2^k-1-b}\,g_{HF}$. Conjugate each output antilinear in $F$. Including the equal labels gives one square:
\[
\left\|
\sum_F
\FigFiveRemainderNorm
\right\|^2
=\left\|\sum_F\bar y_F\otimes t_F\right\|^2
=\sum_{F,G}\overline{Y_{GF}}T^H_{GF}.
\]
Subtract its diagonal. Each $T^H_{FF}$ occurs $q-1$ times on the left of \autoref{eq:five-trace}, and $Y_{FF}=\lambda+\mu p_{FH}\le\lambda+\mu$.
\end{proof}

Write $y=Y_{GF}$. Subtract the square $|\sqrt d A^H_{GF}-y/\sqrt d|^2$ from $d|A^H_{GF}|^2$ in \autoref{eq:five-outer}. By \autoref{lem:five-mixed}, we may drop $T^H_{FF}+2\operatorname{Re}(\bar y T^H_{GF})$ after averaging: fix $H$ for the unweighted term ($q=4$), and $H,I,J$ for the term weighted by $bp_{IJ}$ ($q=2$). The remaining contribution from the off-diagonal entry is
\begin{equation}\label{eq:five-scalar}
(2\lambda-\lambda^2/d)|g|^2+(2b\mu-\mu^2/d)|x|^2+
2(b\lambda+\mu-\lambda\mu/d)\operatorname{Re}(\bar gx).
\end{equation}

\subsection{Closing the cycles}
Set $P=\mathbb E|g|^2$, $U=\mathbb E p_{IJ}|g|^2$ and $V'=\mathbb E p_{IJ}|x|^2$. Let $C_j$ be the mean real weight of an oriented $j$-cycle on distinct labels.

\begin{lemma}[The long cycles]\label{lem:five-long}
$C_4\le U$ and $24C_5\le6U+24V'$.
\end{lemma}

\begin{proof}
Use the straight pairing $l$ and the permuted pairing $h$ below. The solid arrows are the factors of $h$ and the dashed arrows those of $\bar l$. For five labels the extra foot of $l$ contributes $g_{55}=1$.
\[
\begin{array}{c@{\qquad}c@{\qquad}c@{\qquad}cc}
l& h&\bar l h&\mathbb E|l|^2&\mathbb E|h|^2\\[2pt]
\FiveHalfPair{4}{0}&\FiveHalfPair{4}{1}&\FiveMixedCycle{4}&U&U\\[5pt]
\FiveHalfPair{5}{0}&\FiveHalfPair{5}{1}&\FiveMixedCycle{5}&U&V'
\end{array}
\]
In the last row $h=g_{15}g_{32}g_{54}$: edges $15,54$ meet at $5$ and edge $32$ is disjoint. The mixed real entries average to $C_4,C_5$. Average $|l-h|^2$ in the first row and $|\sqrt6\,l-\sqrt{24}\,h|^2$ in the second.
\end{proof}

\begin{proof}[Proof of \autoref{thm:five}]
Put $\Theta=\mathbb E\operatorname{Re}(\bar g x)$ and $\Theta_2=\mathbb E p_{IJ}\operatorname{Re}(\bar g x)$. In the coarse clasp group the permutations by cycle type:
\begin{equation}\label{eq:five-types}
\per(g_{ij})=1+10P+15U+20\Theta+20\Theta_2+30C_4+24C_5 .
\end{equation}
The types, from the left, are $1^5,\,2\,1^3,\,2^2 1,\,3\,1^2,\,3\,2,\,4\,1,\,5$. Each triangle has $\binom53$ choices of labels and two orientations. The other two labels stay fixed or form a pair:
\[
\underbrace{\FiveCycleType{0}}_{g_{12}g_{23}g_{31}}
\qquad\qquad
\underbrace{\FiveCycleType{1}}_{g_{12}g_{23}g_{31}\,p_{45}} .
\]
Reversal gives the conjugate triangle. To match its real part in \autoref{eq:five-scalar}, choose the coefficients in the unweighted and weighted terms as follows:
\[
\begin{array}{c|cc|c|c}
\text{weight}&\lambda&\mu&b\lambda+\mu-\lambda\mu/d
&\lambda+\mu\ \le\ (q-1)/2\\\hline
1&10/7&0&10&10/7<3/2\\
b p_{IJ}&451/2632&1/4&10/b&1109/2632<1/2
\end{array}
\]
In the second row $\lambda=(10/b-\mu)/(b-\mu/d)$. The two triangle coefficients are now $20$ and $b(20/b)=20$.

The straight terms of \autoref{eq:five-outer} give $1+2bP+b^2U$. Add the two contributions \autoref{eq:five-scalar}, subtract \autoref{eq:five-types}, and apply \autoref{lem:five-long}. With row indices $0,1$ in the order above, the coefficients left are
\[
\begin{array}{c|l|r}
P  &2b+2\lambda_0-\lambda_0^2/d-10 &11776/2009\\[2pt]
U  &b^2+b(2\lambda_1-\lambda_1^2/d)-15-30-6 &73895/247408\\[2pt]
V' &b(2b\mu_1-\mu_1^2/d)-24 &47/164 .
\end{array}
\]
Each multiplies a mean of products of pair squares.
\end{proof}

\begin{proof}[Proof of \autoref{thm:matrix}]
For $m=1$ or $k=1$ the two permanents agree. For two and three packets use \autoref{loc:two} and \autoref{thm:three}. At four packets the quadratic case is \autoref{thm:fourquad}. For $k=3$ use \autoref{thm:fourcubic}, and for $k\ge4$ use \autoref{cor:fourthick}. At five packets apply \autoref{thm:five}. The bound for arbitrary packet number is \autoref{thm:threshold}.
\end{proof}

\section{A negative marked closure}\label{sec:limits}

\subsection{Five quadratic packets}

Let $Y_t$ be the transfer of \autoref{eq:iter-mark} with packet $t$ marked. Insert the same product $F_1\otimes\cdots\otimes F_5$ in every diagram. The circles mark the chosen pair of feet.
\[
\LimitsClosure{X_5}{0}
\;=\; \tfrac15 \left(\LimitsClosure{Y_1}{1}+\cdots+\LimitsClosure{Y_5}{5}\right).
\]

\begin{example}[The coefficient of a routing]
Doubled strands represent two-strand cables. Inputs $1,2$ split, as do outputs $2,3$. Count the split ends at the mark:
\[
\FigFivePacketRouting
\qquad
\begin{array}{c|ccccc}
t&1&2&3&4&5\\\hline
\text{split ends}&1&2&1&0&0
\end{array}
\]
Each marked split end has relative coefficient $m/(2b)$ when $b$ inputs split. Here $b=2$ and the five marks give mean $\tfrac15\tfrac54(1+2+1)=1$. In general the $m$ marks count all $2b$ split ends once.
\end{example}

\begin{lemma}
There are five quadratic matrix packets in a plane for which the fifth marked closure is negative.
\end{lemma}

\begin{proof}
For orthonormal $e_0,e_1$, use the basis $e_0^2,e_0\odot e_1,e_1^2$, with squared norms $1,1/2,1$. Take $F_1=F_2=F_3=e_0^2$ and
\[
F_4=(1,1/5,7/400),\qquad F_5=(1,-60,1200).
\]
Both have the form $(e_0+u e_1)\odot(e_0+\bar u e_1)$: take $u=(2+i\sqrt3)/20$ and $u=-30+10i\sqrt3$, respectively.
Evaluating the whole edges by \autoref{sec:certificate} gives
\begin{equation}\label{eq:mark-witness}
\begin{array}{c|ccc}
t&1,2,3&4&5\\\hline
\langle Y_t\rangle&9041984725911/2000&18083516817/4&-565839/2000 .
\end{array}
\end{equation}
The first three entries agree because the packets agree. The mean of all five marked closures is $18083856010197/5000>0$.
\end{proof}

Thus the marked transfers do not admit separate positive decompositions on all matrix packets. The five-block quadratic inequality and the small degrees outside the tests of \autoref{sec:threshold} and \autoref{sec:th-app} remain open here.

\appendix

\section{Coefficients of the reflected cuts}\label{sec:th-app}

\subsection{Keeping a prescribed four-set}
Put $q=|I|$ and $S(I)=U(I)-H(I)$. In $x_i-x_j$ the two exits share the surviving packet $D=g_{ti}F_j-g_{tj}F_i$ at the clasp. Each straight cut can remove a unit spectator while keeping $D$ intact.

\begin{lemma}
For $q\ge4$ the exit squares in \autoref{eq:th-peel} satisfy
\begin{equation}\label{eq:th-faces}
E(I)\ge\frac3q\sum_{|T|=3}E(T),\qquad
E(I)\ge\frac4{q(q-3)}\sum_{|J|=4}E(J).
\end{equation}
\end{lemma}

\begin{proof}
Keep one spectator $h$, then open the last clasp. The double strand contains the other $q-4$ spectators.
\[
(q-2)!\left\|\ThKeep{2}\right\|^2
\;\ge\;2\left\|\ThKeep{1}\right\|^2
\;\ge\;\left\|\ThKeep{0}\right\|^2 .
\]
Both inequalities are straight cuts as in \autoref{eq:th-straight}. Sum the right square with weight $1/q$ per mark $t$, compared with $1/3$ in $E(T)$. For the middle square also average the $q-3$ choices of $h$. Each mark in $E(J)$ has weight $1/4$.
\end{proof}

\begin{lemma}[Four labels and $r$ disjoint pairs]
Put $c_{q,0}=1/(q-3)$ for $q\ge4$. For $r>0$ put
$q c_{q,r}=(q-4-2r)c_{q-1,r}+2r(q-1)c_{q-2,r-1}$, with $c_{q,r}=0$ outside $q\ge4$, $0\le2r\le q-4$. Then
\begin{equation}\label{eq:th-weighted}
S(I)\ge \sum_{r\ge0}c_{q,r}
\sum_{\substack{|J|=4,\; |L|=r\\ L\text{ a matching on }I\setminus J}}p_L S(J) .
\end{equation}
\end{lemma}

\begin{proof}
At $q=4$ both sides equal $S(I)$. Subtract \autoref{eq:th-peel} from the recursion for $U$, leaving $+E(I)$. For fixed $J,L$ there are $q-4-2r$ single deletions outside them and $r$ possible pair deletions in $L$. Their weights give the stated recursion for $r>0$.

At $r=0$ the single deletion gives $1/q$ for each $S(J)$, leaving $3/[q(q-3)]$. On four labels $S(J)=E(J)+\tfrac14\sum_{T\subset J,\,|T|=3}E(T)$ because $S_2=0$ and $S_3=E_3$. Take $1/4$ of the first bound in \autoref{eq:th-faces} and $3/4$ of the second. Each triple belongs to $q-3$ four-sets, giving precisely the missing weight. All remainders are reflected squares.
\end{proof}

\begin{proof}[Proof of \autoref{prop:th-matrix}]
There are $d_{m,r}=\binom m4(m-4)!/[2^r r!(m-4-2r)!]$ choices of $J,L$. Then $h_{m,r}=4d_{m,r}c_{m,r}$ has the claimed initial value and recurrence: for permitted deletions,
\[
\frac{m-4-2r}{m}\,d_{m,r}=d_{m-1,r},
\qquad \frac{2r}{m}\,d_{m,r}=(m-1)d_{m-2,r-1}.
\]
Appending one edge on $J$ extends $L$ in six ways. Appending two edges extends it in three ways. Relabeling makes each extension uniform. Subtract $U(J)=1+7P_J/3+3D_J$ from the right side of \autoref{eq:th-trade}. By \autoref{eq:th-weighted} the resulting change in the upper bound is
\[
-h_{m,r}(1-u_r)Z_{r+1}+h_{m,r}(u_r^{-1}-1)Z_{r+2}
\]
for each $r$ and any $u_r>0$.

The adjacent degrees $j,j+1$ require the two diagonal entries of
$h_{m,j-1}\left(\begin{smallmatrix}u_{j-1}&-1\\-1&1/u_{j-1}\end{smallmatrix}\right)$.
Start at $j=1$. If the current pivot is $a_j$, take $u_{j-1}=a_j/h_{m,j-1}$. The amount used at $j+1$ is then $h_{m,j-1}^2/a_j$, precisely the subtraction in Gaussian elimination. All links are positive. For $M\succeq0$ a pivot before the last cannot vanish, as its row has a nonzero link. The last pivot is nonnegative, giving the comparison at every matching degree.
\end{proof}

\subsection{Cutting an odd circuit}

For a set $Q$ of $q$ labels put $M_j(Q)=\sum_{|L|=j}p_L$, zero when $2j>q$. Let $T_{2r+1}(Q)$ be the sum over oriented $(2r+1)$-cycles in $Q$, with all other packets fixed. Reversing a cycle conjugates its value, so $T_{2r+1}(Q)$ is real.

\begin{lemma}[Odd-cycle cut]
Let $1\le r$, $2r+1\le q$ and $u>0$. Both $B_{q,r}(u)\pm T_{2r+1}(Q)$ are sums of reflected squares, where
\begin{equation}\label{eq:th-oddlocal}
B_{q,r}(u)=\frac{2^{r-1}r!}{2r+1}
\left[\left(u+\frac{q-2r}{u}\right)M_r(Q)
+\frac{2(r+1)}{u}M_{r+1}(Q)\right].
\end{equation}
\end{lemma}

\begin{proof}
Fix ordered distinct labels $x_1,y_1,\ldots,x_r,y_r$. Cut $r$ alternate edges. In the lower half sum the two-step path through all unused $v$.
\begin{center}
\FigAlternatingCuts
\end{center}
The upper half is $A=\prod_jg_{x_jy_j}$. The lower half is $B=(\prod_{j<r}g_{x_{j+1}y_j})\sum_v g_{x_1v}g_{vy_r}$. For $r=2$ the dashed edges come from $\bar A$ and the solid ones from $B$:
\[
\FigAlternatingCycle
\qquad
\bar A B\big|_{v=5}
=g_{21}\,\underbrace{g_{15}g_{54}}_{\text{two-step path}}\,g_{43}g_{32}.
\]
Gluing gives each oriented cycle $2r+1$ times, once for each choice of its middle vertex. Also $\sum |A|^2=2^r r!M_r(Q)$, by ordering and orienting a matching.

Fix the $r-1$ short caps of $B$. There are $2^{r-1}(r-1)!$ choices of order and orientation, and two orders for the ends of the long path. Apply \autoref{lem:th-path} on the $q'=q-2r+2$ available labels. The path adds one or two edges. In the resulting matching these can be marked in $r$ or $\binom{r+1}2$ ways. Hence
\[
\sum |B|^2\le2^r r!\big[(q-2r)M_r(Q)+2(r+1)M_{r+1}(Q)\big].
\]
Sum $|\sqrt u A\pm B/\sqrt u|^2/[2(2r+1)]$ and restore the diagonal remainders from \autoref{lem:th-path}. This gives the two signs.
\end{proof}

\subsection{Disjoint cycles}

Let $\lambda\vdash m$ be a cycle type, with $m_d$ parts of size $d$. Distribute its fixed labels as evenly as possible among the $a$ nontrivial odd cycles, giving $s_i$ to the $i$th. If $a=0$, close each fixed label against itself. Use $z^j$ for the mean closure $Z_j(I)$.

\begin{lemma}
For $u_i>0$ put $L_s(z;u)=\tfrac12(u/(s+1)+1/u)+sz/(2u)$ and define $D_{m,j}$ by
\begin{equation}\label{eq:th-signedsum}
\sum_jD_{m,j}z^j=
\sum_{\lambda\vdash m}\frac{m!}{\prod_{d\ge1}d^{m_d}m_d!}
z^{\sum_{d\in\lambda}\lfloor d/2\rfloor}
\prod_{i=1}^{a} L_{s_i}(z;u_i).
\end{equation}
Then $H(I)\le\sum_jD_{m,j}Z_j(I)$ on unit packets. The block-permanent inequality follows if $D_{m,j}\le\binom nj(2^k-1)^j$ for all $j$.
\end{lemma}

\begin{proof}
On an even cycle $(i_1,\ldots,i_{2r})$ take the two alternating matchings
$A=\prod_j\overline{g_{i_{2j-1},i_{2j}}}$ and $B=\prod_jg_{i_{2j},i_{2j+1}}$ cyclically. Their mixed gluing $\bar AB$ traverses the cycle. The squares $|A\pm B|^2/2$ bound both signs of its real part by half its diagonal sum. Averaging gives the single monomial $z^r$, also for $r=1$.

For an odd cycle on $q=2r+1+s$ labels, divide \autoref{eq:th-oddlocal} by the number $q!/[s!(2r+1)]$ of cycles. There are $q!/[2^j j!(q-2j)!]$ matchings of size $j$ when $2j\le q$. Substitution gives $z^r L_s(z;u)$.

On disjoint sets let $T_i$ be the mean oriented cycle closure and $D_i$ its mean diagonal bound. Each $D_i\pm T_i$ is a sum of reflected squares. Juxtaposing the pictures for $v\ge1$ cycles gives
\[
\prod_{i=1}^v D_i-\prod_{i=1}^v T_i
=2^{1-v}\!\sum_{\substack{\epsilon_i=\pm1\\ \prod_i\epsilon_i=-1}}
\prod_{i=1}^v(D_i+\epsilon_i T_i).
\]
Thus matching degrees add. Uniform relabeling gives the mean of all permutations of this type, also for repeated cycle lengths. The identity gives one. Multiplication by the class sizes and summation give \autoref{eq:th-signedsum}. For the last assertion use \autoref{eq:th-fine}.
\end{proof}

\begin{example}
On six labels take a triangle, a transposition and a fixed point. Give the fixed point to the triangle and put $u=3/2$. The factors for their mean are
\[
\underbrace{
\FigTriangleWeight
}_{z(17/24+z/3)}
\qquad
\underbrace{
\FigPairWeight
}_z .
\]
Their product is $17z^2/24+z^3/3$. The class contains $\binom63\,2\binom32=120$ permutations.
\end{example}

\begin{example}[Ten packets]
Take $u_i$ to be the nearest positive multiple of $1/2$ to $\sqrt{s_i+1}$. Then
\[
(D_{10,j})_{j=0}^5
=\left(1,130,\frac{143339}{15},\frac{555911}{2},2126771,1349585\right).
\]
Each entry is at most $\binom5j31^j$, giving the inequality for ten blocks of size $k\ge5$.
\end{example}

\section{Evaluating closed routes}\label{sec:certificate}

\subsection{Cutting in symmetric words}

\begin{definition}[Occupation words]
Choose an orthonormal basis $e_1,\ldots,e_d$ with $S=\operatorname{span}(e_1,\ldots,e_s)$. For $a=(a_1,\ldots,a_d)$ put
\[
[a]=e_1^{\odot a_1}\odot\cdots\odot e_d^{\odot a_d},\qquad
|a|=\sum_j a_j,\qquad w(a)=\frac{\prod_j a_j!}{|a|!}.
\]
The vector $[a]$ averages $1/w(a)$ distinct orthonormal words. Thus $\|[a]\|^2=w(a)$.
\end{definition}

\begin{example}[The factor at the seam]\label{ex:cert-half}
Take orthonormal $e,f$ and put $S=\C e$. Cap one strand of a $2$-cable against $S^*$. Reflect the cap on $f$ and glue:
\[
\underbrace{\CertTinyCap{0}}_{C^\dagger f}
=(e\odot f)\otimes\bar e,
\qquad
\underbrace{\CertTinyCap{1}}_{\langle f,CC^\dagger f\rangle}
=\|e\odot f\|^2=\tfrac12 .
\]
At the cut $e\odot f=(e\otimes f+f\otimes e)/2$. The projector onto its line is $2|e\odot f\rangle\langle e\odot f|$.
\end{example}

For a route of \autoref{def:inter-cap}, take occupations $a_i=(a_{ij})_j$ on the positive cables, $b$ on the marked dual, and $n$ on the output. Thus $|a_i|=d_i$, $|b|=r$, $|n|=N$ and $b_j=0$ for $j>s$. Put $A=[a_1]\otimes\cdots\otimes[a_h]\otimes\overline{[b]}$.

\begin{lemma}[Occupation formula]
If $\sum_i a_i=n+b$, then $C_\alpha A=D_\alpha(a,b)[n]$ with
\[
D_\alpha(a,b)=w(b)\!
\sum_{\substack{\sum_i b_i=b\\ |b_i|=\alpha_i}}\;
\prod_i\frac{\prod_j\binom{a_{ij}}{b_{ij}}}{\binom{d_i}{\alpha_i}}.
\]
For any $[n]\in Z_\ell$, the entry $G_{\alpha\beta}(\ell)$ of \autoref{prop:routedcaps} is
\[
G_{\alpha\beta}(\ell)= \frac1{w(n)}\sum_A \frac1{\|A\|^2}\;
\FigRoutedGluing
= w(n)\!
\sum_{\substack{b,a_1,\ldots,a_h\\\sum_i a_i=n+b}}
\frac{D_\alpha(a,b) D_\beta(a,b)}{w(b)\prod_i w(a_i)} .
\]
\end{lemma}

\begin{proof}
Split the cable $i$. Exactly $\prod_j\binom{a_{ij}}{b_{ij}}$ of the $\binom{d_i}{\alpha_i}$ choices have occupation $b_i$ at its cap. On the dual side the split has weight $w(b)/\prod_i w(b_i)$. Closing the caps supplies $\prod_i w(b_i)$, and merging the survivors gives $[n]$. This proves the first formula, with an impossible binomial taken as zero.

At the dashed seam insert the orthogonal line projectors $|A\rangle\langle A|/\|A\|^2$. Their denominators are $w(b)\prod_i w(a_i)$. Each outer pairing supplies $w(n)$, and dividing by $\|[n]\|^2$ leaves one factor of $w(n)$. For an additional clasp on packets of equal size, average the route coefficients over the packet permutations before gluing.
\end{proof}

In \autoref{ex:cert-half} the one cap finds $e$ in one of two slots, giving $D=1/2$. The formula returns $(1/2)^2/(1/2)=1/2$.

\subsection{Forcing whole edges}

For the five quadratics in \autoref{sec:limits}, represent
$F_i=a_i e_0^2+b_i(e_0\odot e_1)+c_i e_1^2$
by $p_i(z)=a_i+b_i z+c_i z^2$, and put
$p_I=\prod_{i\in I}p_i$.
Write $[z^j]p$ for the coefficient of $z^j$ in $p$. On $|I|=|J|=n$ packets the mixed fine closure is
\[
U(I,J)=\frac{(2n)!}{2^n}\sum_{j=0}^{2n}
\frac{\overline{[z^j]p_I}\,[z^j]p_J}{\binom{2n}j},
\qquad g_{ij}=U(\{i\},\{j\}),
\]
since $w(2n-j,j)=1/\binom{2n}j$.

\begin{lemma}[Counting whole edges at a mark]\label{lem:cert-mark}
For a partial matching $D$ of $d$ whole edges put
$L(D)=\bigl(\prod_{ij\in D}g_{ij}\bigr)U(I',J')$, where $I',J'$ are the unused bra and ket labels. Let $v_t(D)$ count the ends of $D$ at $t$ on the two boundaries. Set
\[
B_t(z)=\sum_D L(D)\big[2(z-1)^d-v_t(D)z(z-1)^{d-1}\big],
\]
omitting the second term for $d=0$. Then
$\langle Y_t\rangle=\sum_{h=0}^4\frac{5}{2(5-h)}[z^h]B_t(z)$.
\end{lemma}

\begin{proof}
Force a whole edge $j\to i$. Its two thin routings absorb to the same packet cable:
\[
\tfrac12\!\left(\CertForced{0}+\CertForced{1}\right)
=\FigWholePacketEdge
=g_{ij}.
\]
Each forced edge cancels one factor $2$ in the fine normalization. The sum over routings containing $D$ is therefore $L(D)$.

Fix a routing with $h$ whole edges. Summing $(z-1)^{|D|}$ over their subsets gives $z^h$. For each whole end at $t$, summing $z(z-1)^{|D|-1}$ over subsets containing its edge also gives $z^h$. The coefficient is the routing weight times the number of split ends at $t$. For $h=5$ it vanishes. Otherwise $5-h$ inputs split, so \autoref{eq:iter-mark} supplies the factor $5/[2(5-h)]$.
\end{proof}

\begin{example}[The five marked closures]
Take $p_1=p_2=p_3=1$, $p_4=1+z/5+7z^2/400$ and $p_5=1-60z+1200z^2$ in \autoref{lem:cert-mark}. This gives \autoref{eq:mark-witness}. The fine and coarse closures are $14606209031337/4000$ and $695621115897/20000$. Their difference is the mean of the five marks.
\end{example}

\end{document}